\RequirePackage{fix-cm}

\documentclass{svjour3}   
\smartqed  
\usepackage{graphicx}
\usepackage{url}
\usepackage{etex}
\usepackage{amsmath}
\usepackage{adjustbox}
\usepackage{enumitem}
\usepackage{multimedia}
\usepackage{array}
\usepackage{setspace}
\usepackage{color}
\usepackage{textcomp}
\usepackage{gensymb}
\usepackage{float}
\usepackage{multicol}
\usepackage{wrapfig}
\usepackage{amssymb}

\newcommand{\rem}[1]{}

\begin{document}

\title{Trajectory tracing in figure skating}
\author{Meghan Rhodes         \and
        Vakhtang Putkaradze 
}
\institute{M. Rhodes \at
              Department of mathematical and statistical sciences, University of Alberta, 
              Edmonton, AB, T6G 2G1 Canada \\
              \email{mhall@ualberta.ca}          
           \and
           V. Putkaradze \at
              Department of mathematical and statistical sciences, University of Alberta, 
              Edmonton, AB, T6G 2G1 Canada \\ 
              ATCO SpaceLab, 5302 Forand St SW, Calgary, AB, T3E 8B4, Canada\\
              \email{putkarad@ualberta.ca}
}
\maketitle

\begin{abstract}
In this work, we model the movement of a figure skater gliding on ice by the Chaplygin sleigh, a classic pedagogical example of a nonholonomic mechanical system. The Chaplygin sleigh is controlled by a movable added mass, modeling the movable center of mass of the figure skater. The position and velocity of the added mass act as controls that can be used to steer the skater in order to produce prescribed patterns. For any piecewise smooth prescribed curve, this model can be used to determine the controls needed to reproduce that curve by approximating the curve with circular arcs. Tracing of the circular arcs is exact in our control procedure, so the accuracy of the method depends solely on the accuracy of approximation of a trajectory by circular arcs. To reproduce the individual elements of a pattern, we employ an optimization algorithm. We conclude by reproducing a classical ``double flower'' figure skating pattern and compute the resulting controls.
\end{abstract}

\section{Introduction}
\label{intro}
Figure skating is a popular sport worldwide with competitions taking prime time television spots on major networks, with the combination of athleticism and artistic performance making the sport appealing to a broad audience. Underlying the entertaining performances are technical and complex principals from mathematics and physics. Figure skaters are trained to intuitively manipulate physical properties of their bodies such as center of mass, momentum, and inertia, using the relative motion of the limbs, in order to perform complicated movements including footwork patterns, spins, and jumps. 

\begin{figure} 
	\centering
	\includegraphics[width=0.75\textwidth]{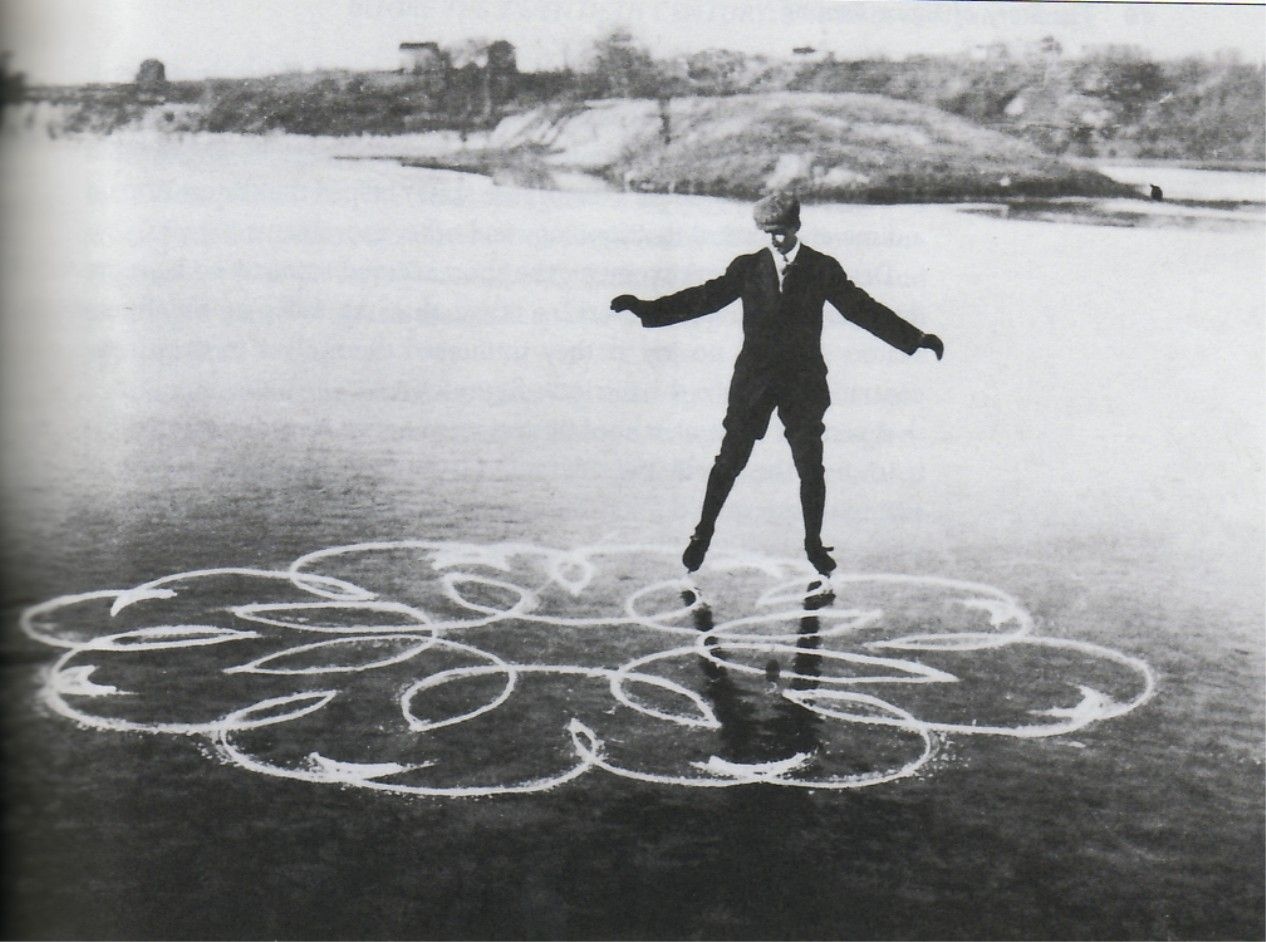}
	\caption{Figure skating pattern on ice.}
	\label{fig:OldPattern}
\end{figure}

Originating in the 19th century, skating derives its name from the patterns, or ``figures'', carved into the ice while skating. Ice skating hobbyists would design intricate patterns and attempt to recreate them on the ice as precisely as possible. An example of such as a design is shown in Figure \ref{fig:OldPattern} \footnote{\url{https://qph.fs.quoracdn.net/main-qimg-b0f3b3d2115f1e1fa1720b7c3733e0bc}}.
As the hobby grew into a sport and ice skate technology advanced, spins and jumps were included in the broadened term of figure skating. The sport of figure skating developed, with categories being made for different forms of figure skating. The name ``compulsory figures'' was given to the process of tracing specific patterns on the ice; ``ice dance'' was used to describe a more dance-like style of skating that excluded jumps and many spins; and the jump- and spin-filled skating style was termed ``free skating''. During testing and competitions, compulsory figures are judged on precision, accuracy, and edge quality (smooth sliding and no grinding), among other metrics. Until 1991, compulsory figures were a required component in competition \cite{hamilton_2019}. Although compulsory figures may be less visually spectacular and less athletically challenging compared to \emph{e.g.} triple and quadruple jumps, performing a precise and complex pattern on the ice requires significant control and technical skill. However, watching competitors slowly trace patterns on the ice was fairly tedious for the mainstream television viewer. Thus, the compulsory skating was excluded from the standard figure skating competitions and  broke off into its own category. Competitions in compulsory figures are held independently and the athletes have their own rankings, independent of the mainstream figure skaters. These competitions are practically never shown on major networks, although plenty of materials are available on video hosting services like \emph{YouTube}. 

Previous mathematical models and descriptions of figure skating have followed this trend in the popularity free style skating, mainly focusing on modeling jumping \cite{hartel2006,king2001,king2005,knoll2005,schaefer2016}. Considering ice skating more generally, much attention has been dedicated to the modeling of a skate's blade gliding on ice with a focus on the friction resulting from ice melting under the blade \cite{Ro2005,LoSzMa2013,PoLeBe2015}. In this paper, we take an idealized approach to the skating and consider an ideal blade that glides in a frictionless fashion on ice, with the gliding motion happening only along the blade's direction. Such assumptions lead to mechanical models lying in the realm of \emph{nonholonomic mechanics} \cite{Bloch2003}, utilizing mechanical models of the skater with affine (in our case, linear) constraints on velocities. 

Several models of skaters have been developed in the context of nonholonomic mechanics. A large portion of the literature has been dedicated to the development of a dynamic model of a skater without control. This direction of work started with the classic development by Chaplygin \cite{chaplygin1911theory}, where a model of the \emph{Chaplygin's sleigh}, or skate, was suggested. Chaplygin's sleigh represents a flat solid body having a certain inertia and center of mass, which is capable of the two-dimensional motion of rotating and sliding on the horizontal ice. The direction of velocity of Chaplygin's sleigh at every point in time must be parallel to a direction specified by a blade attached to the body. In the classical non-controlled Chaplygin's sleigh, the blade's position and orientation with respect to the body is fixed. Notably, the problem of the classical Chaplygin's sleigh is integrable, representing one of the very few known examples of integrability in nonholonomic mechanics \cite{Kozlov1985}. While this model is certainly interesting from the mathematical point of view \cite{Bloch2003,BiBoKoMa2018}, it is much too simple for a realistic description of the skater's motion, which is controlled by the complex movement of the body. For a more realistic description, there have been many generalizations of the classical Chaplygin's sleigh problem in the literature. While we try to keep the literature review relatively compact here, of particular interest to us are the extensions including a movable mass positioned on the sleigh. For the extension of the sleigh's motion without active control, we note  \cite{Ze2003,BlMaZe2009} which considered a movable mass on a spring connected to a given point. As it turns out, that system happens to be integrable, although without the explicit expressions of the integrals of motion allowable by the classical Chaplygin's sleigh. Another extension that has received much interest in literature consists of a sleigh with a mass that is being forced to move in a periodic fashion \cite{bizyaev2017chaplygin}. Such motion leads to Fermi-like acceleration of the sleigh with unbounded energy for certain initial conditions and values of parameters. This work has been further extended in \cite{bizyaev2019chaplygin} with incorporation of viscous friction in the sleigh's motion, and numerical analysis of the zones of parametric resonance. Finally 
\cite{gzenda2019integrability} modeled a static three dimensional figure skater with a model that includes a lean with respect to the vertical direction, but no control or relative change of mechanical properties of the skater. It was shown that for the case when the projection of the center of mass onto the skate's direction coincides with the contact point, the system is integrable, whereas for other configurations the motion is chaotic. An interesting three dimensional  extension of controlled motion of Chaplygin sleigh, not related to skating, was suggested in \cite{fedorov2010hydrodynamic}, in the context of the description of hydrodynamic motion. 

The question of the controls the skater utilizes to obtain desired trajectories on ice is more complicated, especially regarding the mathematical principles behind tracing a desired trajectory on ice. Of course, the detailed mathematical model of the control leading to the actual compulsory figures is currently (and, most likely, will forever be) out of reach for analytical models that are the point of this article. In practice, the skater controls their center of mass by moving the positions of arms, legs, and torso, resulting in trajectory changes. The control is produced in the body frame, while the resulting traces are in the coordinate frame fixed in space (spatial frame). The mechanisms figure skaters use to map the body frame to the spatial frame is very complex, not very well understood, and takes years of practice for skaters to learn precise and accurate control. This lack of current theoretical understanding is even shown in the methods used to teach figure skating. For example, there are many different techniques for each figure skating element, with different body positions considered ``optimal'' in each technique. Often times, when learning a new element, the skater experiments with technique and tries to simply remember the physical feeling of a successful attempt in order to recreate the success of the trajectory on ice. The level to which a skater analyzes their movements and the resulting impact on the physical properties obviously varies skater to skater, but it is not a common practice to significantly examine this relationship. 

In this work, we focus on the controls involved in the field of compulsory figures, based on the simple example of the Chaplygin's sleigh with a moving mass. The mass is moved in such a way that the trajectory of the blade on ice is as close as possible to the desired trajectory.  It seems to us that if we consider this system as a model of a human skater, one should consider the position   of the center of moving mass as the controls, as opposed to the forces exerted on the mass being the controls. Indeed, in everyday life, we can move our hands and legs to the desired position relative to the torso, whereas the forces needed for this position to be obtained are calculated implicitly by the brain to achieve the desired position. 

Compared to the extensive literature dedicated to the dynamics of the sleigh, there have been substantially less work in the literature on the control. In applications relevant to skating, we point out the control of the motion of Chaplygin's sleigh using impulsed force \cite{tallapragada2017steering}, modeling the push off ice. There is also the work on the control of a two-link Chaplygin sleigh similar to the control of two trajectories by changing the relative angle between the skates \cite{fedonyuk2017dynamics}, and the work on the trajectories of the Chaplygin's sleigh forced by the periodic motion of the internal rotor \cite{fedonyuk2018sinusoidal}. However, the actual motion of the skater using trajectory tracing is done gliding on just one skate, with very few additional pushes allowed during the process. This is achieved by changing the position of center of mass and moments of inertia using carefully executed  motion of the body. Thus, the most relevant previous work for our control procedure on trajectory tracing is \cite{osborne2005steering}, where Chaplygin's sleigh was steered based on the prescribed motion of a movable mass. In our work, we show how to develop an optimal control procedure that is capable of tracing a large variety of curves on the ice, possibly including cusps. A short one-page pedagogical review of this work has appeared in SIAM News \cite{RhPu2021}. 

In the following sections, we will present a mathematical description of the dynamics of a blade moving across the ice and describe how controls can be determined to produce a prescribed pattern. As an example, we present a control of the sleigh leading to tracing of a compulsory figure skating pattern illustrated in Figure~\ref{fig:OldPattern}. 

\section{Mathematical Background: Chaplygin's Sleigh}
\label{sect:ModelMvmt}
\begin{figure} 
	\centering
	\includegraphics[scale=0.75]{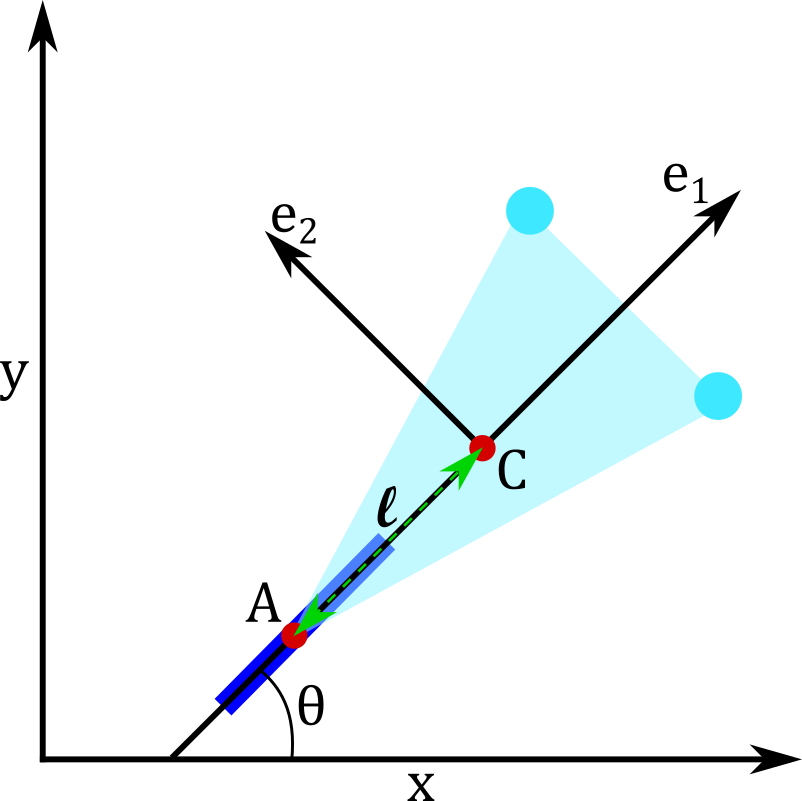}
	\caption{An illustration of the classical Chaplygin sleigh. The solid blue line shows the blade forming the direction of the motion at each point. }
	\label{fig:BasicSleigh}
\end{figure}

The classical Chaplygin sleigh, a model created and analyzed by Chaplygin in 1911 \cite{chaplygin2008theory}, describes a platform of mass $M$, with center of mass at position $C$, which has a knife edge with a single contact point (at $A$) and is supported by two frictionless points that slide freely. A schematic of this mechanical model is depicted in Figure \ref{fig:BasicSleigh}.
The frame attached to the sleigh (the body frame) is defined by the vectors $(\mathbf{e}_1, \mathbf{e}_2)$. The coordinates in the spatial frame are $(x,y)$. The angle between the blade and the $x$-axis of the spatial frame is denoted by $\theta$. The configuration manifold of the Chaplygin sleigh is thus the group of two dimensional rotations and translations $SE(2)$. In the local coordinates, that group can be described by the variables $(x,y,\theta)$.  In Figure~\ref{fig:BasicSleigh}, the blade is aligned with the vector $\mathbf{e}_1$, although in general it does not have to be the case. The Chaplygin sleigh includes a constraint that the direction of movement along must be along the axis of the blade, that is, in the direction of $\mathbf{e}_1$. Mathematically, this constraint is given by 
\begin{align}
\label{constr_nonholonomic}
-\dot x \sin \theta + \dot y \cos \theta =&0\,.
\end{align}
The constraint on the velocity cannot be written exclusively in terms of coordinates on the configuration manifold $(x,y,\theta)$, making the Chaplygin sleigh a nonholonomic system. 

\paragraph{A brief derivation of Chaplygin's sleigh's equation} 
For completeness of the exposition, we show how to derive the Chaplygin sleigh equations using the Lagrange-D'Alembert principle (see \cite{Bloch2003} for full details and alternative derivations).

Denote the coordinates for the sleigh's center of mass as $C=(x_C, y_C)$. If $\ell$ is the distance between the center of mass and the contact point $C$, this distance is along the coordinate $\mathbf{e}_1$, and the coordinates of the contact point $(x,y)$ are 
\begin{align}
x_C=&x + \ell \cos\theta\,, \quad \ y_C=y + \ell \sin\theta\,.
\label{xy_xCyC}
\end{align} 
Then, 
\begin{align}
\dot x_C=&\dot x - \ell \dot \theta \sin\theta\,, \quad  \dot y_C=\dot y+ \ell \dot \theta \cos\theta\,.
\label{dotxy_xCyC}
\end{align}
The Lagrangian for this system is simply the kinetic energy and is given by 
\begin{align}
L=&\frac{1}{2}M(\dot x_C^2+\dot y_C^2)+\frac{1}{2}I\dot \theta ^2\\
=& \frac{1}{2}[M(\dot x^2+\dot y^2)+(I+M\ell^2)\dot \theta^2+2M \ell \dot \theta(-\dot x \sin \theta+\dot y  \cos\theta)]\,,
\end{align}
where $I$ is the moment of inertia. Notice that the last term in the Lagrangian contains the factor that vanishes by the nonholonomic constraint \eqref{constr_nonholonomic}. However, we are not allowed to substitute that constraint in the Lagrangian directly as this leads to \emph{vakonomic} mechanics and different equations of motion. Instead, we take the variations of action using the critical action principle 
\begin{equation} 
\delta \int L \mbox{d} t =0 \, , \quad 
\mbox{on variations} \quad - \delta x \sin \theta + \delta y \cos \theta =0 \, . 
\label{var_action} 
\end{equation} 
If $\lambda$ is the Lagrange multiplier for the constraint on variations in \eqref{var_action}, the equations of motion are
\begin{align} \label{eq:x1}
\ddot x -\ell \dot \theta ^2\cos\theta-\ell \ddot \theta\sin \theta =&\frac{-\lambda\sin\theta}{M}\\
\label{eq:y1}
\ddot y -\ell \dot \theta ^2\sin\theta+\ell \ddot \theta\cos \theta =&\frac{\lambda\cos\theta}{M}\\
\label{eq:theta1}(I+M\ell^2)\ddot\theta+M\ell (\ddot y\cos\theta-\ddot x \sin\theta)=&0\,.
\end{align}
With velocity given by $\label{eq:Vel} v=\dot x \cos \theta +\dot y \sin \theta$, we differentiate the constraint \eqref{constr_nonholonomic} to 
arrive at  
\begin{align}(I+M\ell^2)\ddot\theta+M\ell\dot \theta v=0\,.\end{align}
Taking the time derivative of the velocity and using \eqref{eq:x1} and \eqref{eq:y1}, we get $\dot v= \ell \omega^2$, where we denote the angular velocity as $\omega=\dot \theta$. This gives the final system
\begin{equation}
\begin{aligned} 
\dot v =& \ell \omega^2 \, , \quad 
\ \dot\omega=&-\frac{M \ell}{(I+M\ell^2)}\omega v\,.
\end{aligned}
\label{Chaplygin_standard}
\end{equation} 
There have been many extensions to, and applications of, the Chaplygin sleigh of which we provide a few examples (see  \cite{borisov2016regular,tallapragada2017steering,fedonyuk2017dynamics,fedonyuk2018sinusoidal,kuznetsov2018regular,bizyaev2017chaplygin,jung2016nonholonomic,fedorov2010hydrodynamic,BlMaZe2009,bizyaev2019chaplygin}). 
\begin{figure} 
	\centering
	\includegraphics[scale=0.75]{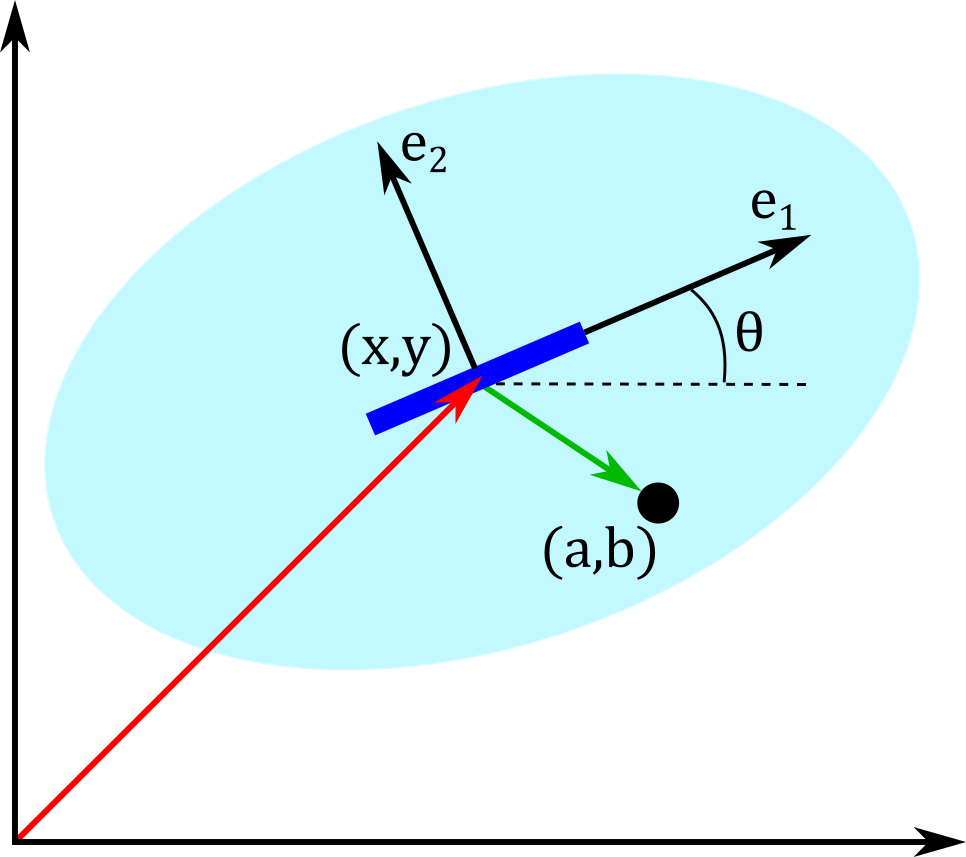}
	\caption{Chaplygin sleigh with added mass.}
	\label{fig:SleighSteer}
\end{figure}
We will employ the Chaplygin sleigh with added moving mass as the simplest model of a figure skater on the ice. The position and velocity of the mass are the controls of the system. The specifications of the Chaplygin sleigh, along with the extension of applying an added mass, make the Chaplygin sleigh a good choice for modeling figure skating, and compulsory figures in particular. First, the controlled Chaplygin sleigh is applicable to figure skating because of the unique profile of the figure skate. Unlike hockey skate and speed skate blades, figure skate blades are curved from front to back. This curve along the length of the blade means that only a small portion of the blade contacts the ice at any given time, corresponding to the single contact point in the Chaplygin sleigh model. Of course, in the full 3D motion, the contact point will move back and forth along the blade, an effect which we neglect in this paper. In contrast, hockey and speed skate blades have large flat portions on the blade, and consequently more of the blade contacts the ice at all times. Second, the nonholonomic constraint on velocity is equivalent to the requirement of smooth gliding in compulsory figures, except for a finite set of predefined points where the skater is expected to turn. The velocity constraint clearly would not hold during other figure skating moves such as jumps where the skate leaves the ice. Finally, the mass added to the classical Chaplygin sleigh can be viewed as moving the skater's center of mass. 

\section{Control of a Skater's Trajectory Using an Added Mass}
\subsection{General Considerations} 
\label{sect:Control}
One of the most basic ways to introduce control of a skater's motion is to introduce an added mass $m$ that can move with respect to the skater. In reality, a skater will use the complex motion of limbs that change both the position of the center of mass, the moment of inertia and, in addition, introduce additional forces and torques from the motion of the limbs. The controlled motion of a single mass can be thought of as a single, two-dimensional version of a multi-limb control. 

There are many works that address the dynamical evolution of the sleigh with an added mass. In particular, we will point out a surprising result \cite{BlMaZe2009} that the problem of a sleigh with added mass $m$ on a spring is integrable. The work \cite{bizyaev2017chaplygin} studied the motion of the Chaplygin sleigh with a prescribed, periodic motion of the added mass. The interesting part about this problem is the existence of trajectories with unlimited growth of energy, and thus unbounded acceleration, which was described as 
``an analog to Fermi acceleration''. However, it is important to point out that not all trajectories lead to indefinite increase of energies, so the periodic motion of the added mass is not a guarantee of indefinite energy growth. 

The problem of making predefined curves on ice, also known as trajectory tracking, using a movable mass has received much less attention in the literature. In \cite{osborne2005steering} it was shown that it is possible to go from one straight line to another using a predefined motion of the moving mass. We are however not aware of works solving the trajectory tracing of complex trajectories as shown on Figure~\ref{fig:OldPattern}. This is precisely the point of this article. We shall present an algorithm that can be used to trace a large class of trajectories using suitable approximations of the curves. 

\subsection{Equations of Motion}
\label{sect:Model}
Let the angular momentum, measured relative to the contact point,  and the linear momentum, measured along the direction of the blade, be denoted $p_1$ and $p_2$, respectively. The spatial coordinates of the blade contact point are given by $(x,y)$, while $\theta$ is the angular orientation of the blade. The configuration space for this system is $SE(2)$. The position of the moving mass is $(a,b)$ in the coordinate system of the sleigh. The equations of motion, generalizing \eqref{Chaplygin_standard} were derived in \cite{osborne2005steering} and are written as
\begin{align}
\dot{p_1} =& -m\eta \xi^2\,, \label{p1_eq}\\
\dot{p_2} =& m\eta \xi^1\,, \label{p2_eq} \\
\dot{\theta} =& \xi^1\,, \label{theta_eq}\\
\dot{x} =&\xi^2\cos \theta\,,\label{x_eq}\\
\dot{y} =& \xi^2\sin \theta\,,\label{y_eq}
\end{align}
where
\begin{align}
\xi^1=&\frac{(M+m)(p_1-ma\dot b)+mb(p_2+M\dot a)}{(M+m)(I+ma^2)+Mmb^2}\,, \label{xi1_eq}\\
\xi^2=&\frac{m[b(p_1-ma\dot b)-(I+ma^2)\dot a]+[I+m(a^2+b^2)]p_2}{(M+m)(I+ma^2)+Mmb^2}\,, \label{xi2_eq}\\
\eta=&\frac{[Mmb^2+I(M+m)]\dot b +a[(M+m)p_1+mb(p_2+M\dot a)]}{(M+m)(I+ma^2)+Mmb^2}\,. \label{eta_eq} 
\end{align} 
For further details of the derivation of the above system, we refer the reader to \cite{osborne2005steering}. The control variables are the position $(a(t), b(t))$ and the velocities $(\dot a, \dot b)$ of the added mass $m$. While $a(t)$ and $b(t)$ can be arbitrary functions with sufficient smoothness properties, for the simulations in this paper we choose the control functions to take the functional form 
\begin{align} 
a(t)=A\cos (\omega t)\, , \quad 
b(t)=A\sin (\omega t)\, , \label{eqn: controls}
\end{align}
where $A$ and $\omega$ are parameters to be optimized to generate a prescribed curve. In other words, the control mass is moving along a circle with a given angular velocity with respect to the sleigh. 

Suppose a trajectory on ice is given. We shall define the trajectory as a piecewise smooth curve $(X(s),Y(s))$ where $s$ is the arc length. The curve tracing problem we are interested does not care about the speed with which the curve is traced. Therefore, we formulate the general trajectory tracing problem as follows. 

In our computations, we shall use the difference between the desired and actual trajectory defined as follows. If $x=x_s(t)$ and $y=y_s(t)$ is the solution curve parametrized by time, we recompute the parametrization of the curve by the arc length $s$, to obtain $x=x_s(s_s)$, $y=y_s(s_s)$. If $S_d$ is the length of the desires solution curve, and $S_s$ is the length of the actual solution curve, we take $S={\rm min}(S_s,S_d)$. Normally, we always have $S_d>S_s$ since the desired solution curve is known analytically, so $S=S_s$. Our control procedure then minimizes the quantity 
\begin{equation} 
C=  \int_0^S \big| x_s(s) - X(s) \big|^p+\big| y_s(s) - Y(s) \big|^p \mbox{d} s  \, , \quad p>1
\label{C_control} 
\end{equation}
\emph{i.e.}, the $p$-th power of $L_p$ norm of the distance between the desired and actual trajectory. 
We then formulate the following general problem.
\begin{problem}[General statement of control procedure]
	\label{prob:general} 
	Suppose a given piecewise plane curve $x=X(s),y=Y(s)$  forms a graph $G$ on $(x,y)$ plane.  Find the initial conditions and controls $(a,b,\dot a, \dot b)$ such that the graph $G_s$  of the solution curve given by (\ref{p1_eq}-\ref{y_eq}) minimizes the deviation defined by \eqref{C_control}. 
\end{problem} 
As it turns out, the problem as formulated above is difficult to solve in the general case, as the mapping from initial conditions and controls to the trajectories is highly nontrivial, and is thus difficult to implement in any robust way. Instead, we develop a procedure that can approximate  any piecewise smooth curve using circular arcs. 

Circular arcs play a special role in the dynamics of the uncontrolled sleigh since they represent some of the exact solutions of the dynamics. Surprisingly, circular trajectories also play a special role in the controlled problem. More precisely, circles yield \emph{algebraic} equations for control variables, which greatly simplifies designing the control procedure. This can be seen as follows. 
\begin{lemma}[On control of circular trajectories] 
	\label{lemma:arcs}
	A controlled solution trajectory is a circular arc of radius $r$ if and only if the controls satisfy 
	\begin{equation} 
	\xi^2 = r \xi^1 \, , 
	\label{arcs_xi1_xi2}
	\end{equation} 
	where $\xi^1$ and $\xi^2$ are given by \eqref{xi1_eq} and \eqref{xi2_eq}, respectively.
	\\
	In the particular case of straight lines, the controls satisfy $\xi^1=0$. \\ 
	Moreover, there is an additional constant of motion $P$ defined as 
	\begin{equation} 
	\mbox{{\rm arcs:}} \quad p_1 + \frac{p_2 }{r}  = \frac{P}{r} = {\rm const}, \quad 
	\mbox{{\rm straight lines:}} \quad p_1 = P = {\rm const}\,.
	\label{p1_p2_const} 
	\end{equation} 
\end{lemma} 
\noindent
{\bf Proof}. From \eqref{x_eq} and \eqref{y_eq} 
we obtain $\mbox{d} s = \xi^2 \mbox{d} t$. Additionally, from \eqref{theta_eq}, we obtain 
\begin{equation} 
\frac{\mbox{d} \theta}{\mbox{d} s } = \frac{\xi^1}{\xi^2}
\label{dtheta_ds}\,.
\end{equation} 
A curve is a circle of radius  $r$ if and only if  $\theta'(s)=1/r$,  giving exactly \eqref{arcs_xi1_xi2}. In the particular case of a straight line, the limiting procedure gives $\xi^1=0$. 

When \eqref{arcs_xi1_xi2} is valid, $p_1$ and $p_2$ are coupled as well from \eqref{p1_eq} and \eqref{p2_eq}, allowing us to notice that
\begin{align}
\dot{p_1} =- m\eta \xi^2 = -m \eta (r \xi^1) \quad 
\implies \dot{p_1}=-r \dot{p_2}
\quad 
\end{align}
and \eqref{p1_p2_const} follows. 
\\[2mm] 
Note that the condition for the sleigh to follow a straight line with constant velocity 
\begin{align}
(M +m)p_1 +m b p_2 = 0\,.
\end{align}
found in \cite{osborne2005steering} follows exactly from $\xi^1=0$ with $\dot a=0$ and $\dot b=0$. 

As the next step in the control procedure, we note the result by Meek and Walton \cite{meek1995approximating} which presented ``an algorithm for finding arbitrarily close arc spline approximation of a smooth curve'' : 
\begin{lemma}[On approximating smooth curves by circular arcs]
	If the bounding circular arcs enclose a given spiral segment of positive curvature $Q(s)$, $s_0\leq s \leq s_1$, and a biarc that matches the same data as the bounding
	circular arcs is found, then the maximum distance between the biarc and the spiral is $\mathcal O (h^3)$, where
	$h=s_1 -s_0$. 
\end{lemma}
Thus, we approximate any smooth part of the trajectory by circular arcs as in \cite{meek1995approximating}, and use the result of Lemma~\ref{lemma:arcs} to design the control procedure, outlined below. To complete the control of the trajectory in a realistic case, as illustrated on Figure~\ref{fig:OldPattern}, we notice that the cusps make the trajectories non-smooth. The cusps are executed by an experienced ice skater by performing a quick turn exactly at the moment when the speed of the skate with respect to the ice vanishes. It is only at these points the finite turn is possible. We shall note that there are also cusps which can be executed without the finite turn, when two arcs touch tangentially \cite{gzenda2019integrability}, but we do not focus on this cases here. Thus, we allow that at the points when the velocity vanishes, a finite turn can be executed, and a new arc can be started with zero velocity. The control algorithm is thus described as follows. 
\begin{enumerate} 
	\item Separate the desired trajectory into smooth parts. 
	\item Approximate every smooth part by circular arcs.
	\item Find a trajectory with control satisfying \eqref{arcs_xi1_xi2} following \emph{exactly} the chosen circular arc. At all the cusp points, the velocity of the skate with respect to the ice must be zero. 
	\item Join circular arcs by allowing a finite turn at the cusp point. 
\end{enumerate} 
Notice that the advantage of this procedure is that the difficulty of finding the controls for an arbitrary trajectory is now substituted by approximating the curve with circular arcs as in \cite{meek1995approximating}. Once the approximation is found, we only need to find the controls enforcing vanishing velocity at the end points, as with the right initial conditions the solution is \emph{guaranteed} to follow the arc. Moreover, there is an additional bonus of simplified solution due to \eqref{p1_p2_const} when using circular arcs or straight lines. 

\subsection{Chaplgyin Sleigh Model Reduced to Circular Trajectories}

Under the conditions that trajectories are circular arcs, the Chaplygin sleigh figure skate model is  
\begin{align}
\dot{p_1}(t) =& -m\eta \xi^2\,,\\
\dot{\theta}(t) =& \xi^1\,,\\
\dot{x}(t) =&\xi^2\cos \theta\,,\\
\dot{y}(t) =& \xi^2\sin \theta\,.
\end{align}
with initial conditions 
\begin{align}
p_1(0) =& \bar p_{1}\,,\\
\theta(0)=& \bar \theta\,,\\
x(0) =&\bar x\,,\\
y(0) =&\bar  y\,,
\end{align}
and the constraints
\begin{align}
\xi^2=&r\xi^1\,,\qquad p_2=\frac{c-p_1}{r}\,. 
\label{eqn: circle constraints}
\end{align}
As described in Section \ref{sect:Control}, zero speed is required at cusp points. The speed with respect to the ice $v(t)$ is given by 
\begin{align}
v(t)=\sqrt{( \dot x)^2+(\dot y)^2}=\xi^2\,.
\end{align}
Therefore, the requirement that the speed be zero at cusp points is equivalent to $\xi^2=0$ at such points.
\subsection{Target Curve Parsing and Optimization}
As a first step in reproducing a given curve, we must define the target trajectory. This means parsing and dividing the trajectory into components of straight lines and circular arcs. In this work, this was accomplished manually with the target curves being divided into approximate curves (eg. semicircle, 3/4 circle, etc.). This process could surely be automated using techniques such as machine learning, but this is beyond of the scope of this work. Possible metrics to define target curves could include arc length, radius, start/endpoints, etc. In the following sections, we use arc length as the metric to describe curve components (a point-based approach is also explored in the appendix). We note that in the examples we present here, all the pieces of the curves are circular arcs. 

To reproduce a given curve component with the Chaplygin figure skate model, we use optimization against the target trajectory to determine the controls $a(t)$ and $b(t)$. In order to enforce the requirement of zero speed at cusps while maintaining smooth control functions, circular arcs are produced by combining a forward-in-time and backward-in-time solutions. 

The cost function considered is based on desired arc length and vanishing speed at the end points of the arc. Although not a simple task for a general curve, as described in \eqref{C_control}, determining the difference between two circular arcs is trivial. In particular, two circular arcs of the same length and radius, starting at the same point and tangent to each other, are going to be exactly identical. Thus, assuming that the solution circle and the desired circle have the same radius, we use the following cost function 
\begin{align} \label{eqn:CLengthandSpeed}
C_{length}=|L_s-L_{d}|^p, \quad p>1, \quad \mbox{with} \quad \xi^2(t_1) = \xi^2 (t_2) =0 \,, 
\end{align}
where $L_s$ is the length of the solution circular arc, $L_d$ is the length of the desired (target), and $t_1, t_2$ are the start and end times of the trajectory, respectively. The arc length of the trajectory between two angles of the normal $\Psi_1$ and $\Psi_2$ is given by $L=r(\Psi_2-\Psi_1)$. 

The control parameters are optimized so that the total trajectory (combined backward- and forward-in-time solutions) reaches an optimal length and such that the speed vanishes at each end. In this case, both the forward- and backward-in-time trajectories are solved within the optimization function. Hence, both the forward- and backward-in-time trajectories are considered when determining the optimal control parameters. 

\section{Numerical Methods}
\label{sect: Numerics}
The simulations of the Chaplygin sleigh were carried out using the Python programming language, version 3.7.6 (Python Software Foundation,\\ https://www.python.org/). 
To simulate a single arc with the Chaplygin sleigh, the work flow is as follows:

\begin{enumerate}
	\item The simulation parameters are defined, including integration time (denoted $T$), values for $I, m, M, r$, ODE initial conditions, and target trajectory metric (ie. $L_{opt}$). 
	
	\item The control functions are defined and values set for an initial guess of control parameters to be used in optimizing the control function.
	
	\item Optimize the control parameters by solving the ODE within an optimizing algorithm, scipy.optimize.minimize\footnote{With default solver (a quasi-Newton method of Broyden, Fletcher, Goldfarb, and Shanno (BFGS))}, with the cost function as previously defined. This optimizer iteratively solves the ODE using the solver scipy.integrate.solve\textunderscore ivp
	\footnote{With solver option ``BDF'', an implicit multi-step variable-order method based on a backward differentiation formula for the derivative approximation.} until an optimal optimal parameter set, $(A, \omega)$, is determined.
	\item The control parameters are optimized so that the total trajectory (combined backwards and forwards trajectories) reaches an optimal length where the speed at each end vanishes. Thus, both the forward and backwards trajectories are solved within the optimization function and are considered when determining the optimal control parameters. The length and speed-based cost function, (\ref{eqn:CLengthandSpeed}), defines the error to be minimized by the optimization function.
	\item The ODE system is then solved with the optimized parameter set. The forwards and backwards trajectories are solved separately to ensure that both ends of the combined arc would have $\xi^2=0$. Note that in both the optimization and when solving the ODEs with the optimal control parameters, an integration termination condition was set so that the integration would stop when $\xi^2$ vanished. 
	
	\item The forward time and backward time solutions are then concatenated to give the full arc.
\end{enumerate}

Individual arcs can then be transformed (via translation and rotation) and placed to produce the final pattern. 

\section{Simulations}

As a proof of concept, we reproduce the compulsory figure shown in Figure \ref{fig:OldPattern}. In order to recreate this pattern, the base components must be identified. In this case, there are two separate curves, each with repeating patterns.

Between the inner and outer curves, the inner curve is simpler with only one arc component that is slightly less than a semi-circle (termed ``arc 1''), which is then repeated eight times. The outermost curve can be decomposed into two different arcs: a short arc (termed ``arc 2'') and a longer arc that is slightly more than a semi-circle (termed ``arc 3''). It can also be noticed that the radius of the inner arc components is larger than that of the outer arc components. The full outer curve can then be represented by sequentially connected triplets consisting of a long arc, followed by a short arc, completed by another long arc, and then repeating this triplet, or ``leaf'', eight times.

Next, we will reproduce the full pattern by generating these component arcs using optimization of the cost function \eqref{eqn:CLengthandSpeed}. The simulation parameters for each arc are given in Table \ref{tab:NumericsParsNewFunctions}. 

\begin{table}
	\caption{Values used in simulations to create arc 1, arc 2, and arc 3. In all cases $m=1$, $M=2$, and $I=3$.}
	\label{tab:NumericsParsNewFunctions}
	\centering
	\begin{tabular}{llllll}
		\hline\noalign{\smallskip}
		Arc  & $T$ & $r$ & $L_{opt}$ &  $(\bar{p_1},\bar{p_2}, \bar{\theta},\bar{x},\bar{y})$ & $(A, \omega)$ \\ 
		\noalign{\smallskip}\hline\noalign{\smallskip}
		1 & $6$ & $1.2$& $1.1r\pi$ & $(2, 3, 0, 0, -1.2)$ & $(1, 1)$\\
		2 & $4.1$ & $0.8$&  $0.2r\pi$ & $(0.15, 0.5, 0, 0, -0.8)$ & $(0.1, 1)$\\
		3 & $4.1$& $1.0$ &  $r\pi$ & $(2,3, 0, 0, -1)$ & $(1, 1)$  \\ 
		\noalign{\smallskip}\hline
	\end{tabular}
\end{table}

The initial control parameter guess, resulting optimized control parameters, and optimization errors are shown in Table \ref{tab:OptimizedParsNewFunctions}. As indicated by the small values for optimization error in Table \ref{tab:OptimizedParsNewFunctions}, the length and speed-based optimization was able to successfully produce arcs that reached the target length. The optimization error could be further minimized by choosing optimal simulation parameters, but for use in this work we were satisfied the given results. Choosing optimal simulation parameters that would yield smaller optimization errors could likely be automated using a machine learning algorithm.

\begin{table}
		\caption{Initial guesses and optimized values for control parameters.}
	\label{tab:OptimizedParsNewFunctions}
	\centering
	\begin{tabular}{llll}
		\hline\noalign{\smallskip}
		Arc  & Guess, $(A, \omega)$& Optimized, $(A^*, \omega^*)$ & Error\\
			\noalign{\smallskip}\hline\noalign{\smallskip}
		1 & $(1, 1)$ & $(0.799, 3.283)$& $1.899 \times 10^{-7}$\\
		2 & $(0.1, 1)$& $(-0.343, 1.369)$ & $6.144 \times 10^{-4}$\\
		3 &  $(1, 1)$ & $(-1.379, 2.352)$ & $9.316 \times 10^{-8}$\\ 
		\noalign{\smallskip}\hline
	\end{tabular}
\end{table}

\begin{figure} 
	\includegraphics[scale=0.4]{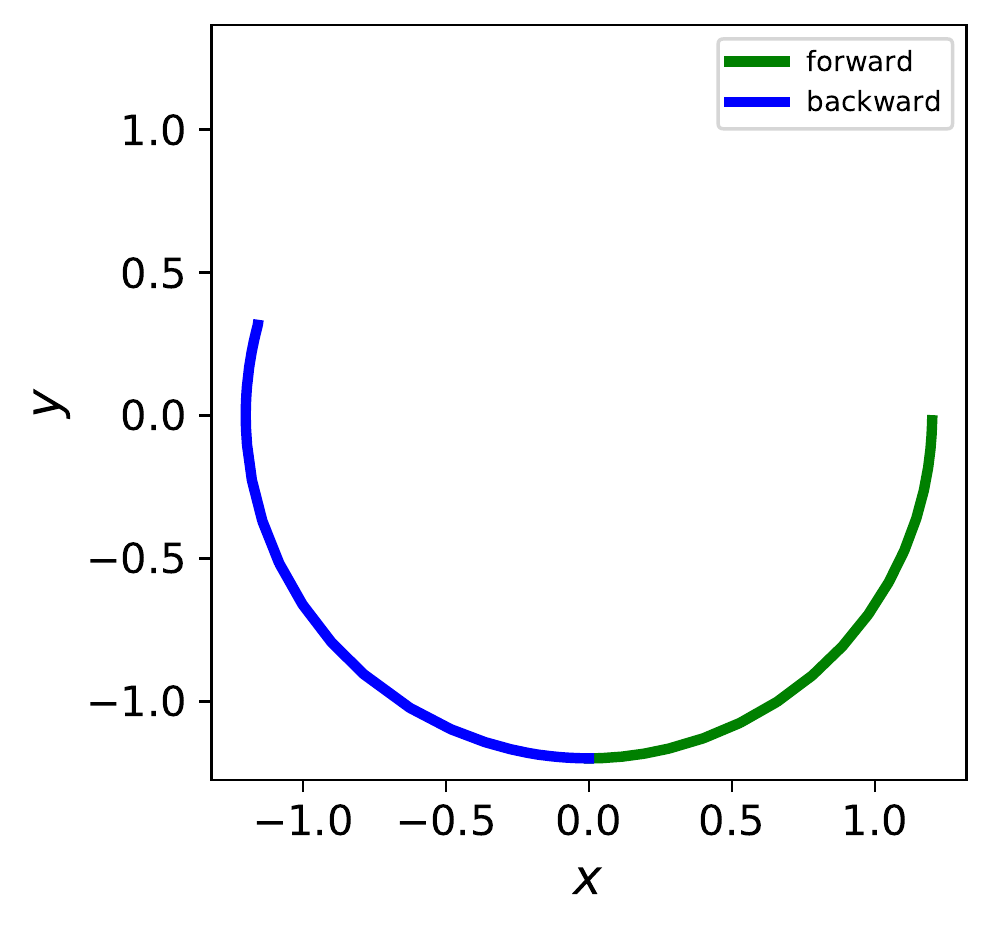} \hfill	\includegraphics[scale=0.4]{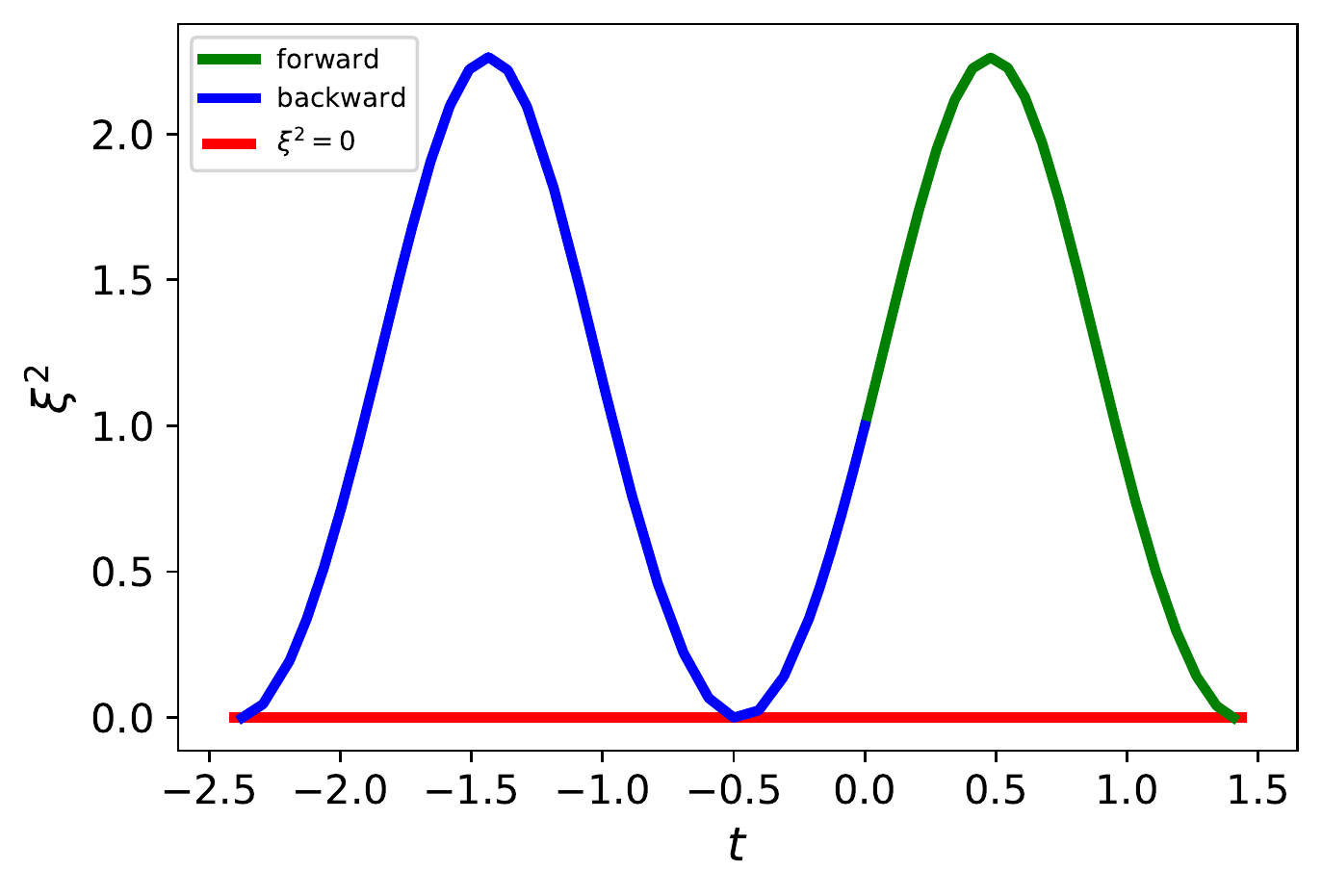}\hfill
	\includegraphics[scale=0.4]{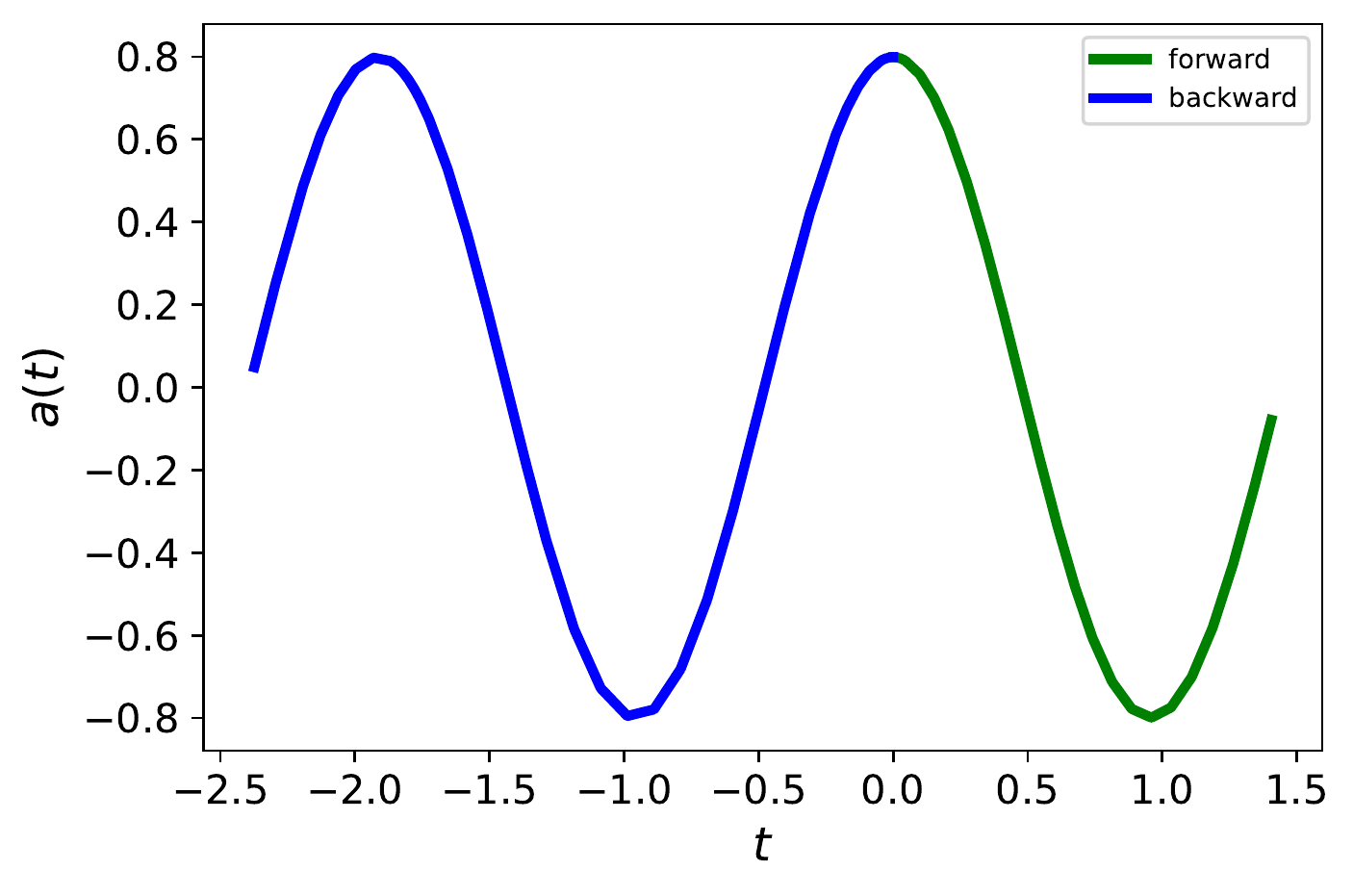}\hfill
	\includegraphics[scale=0.4]{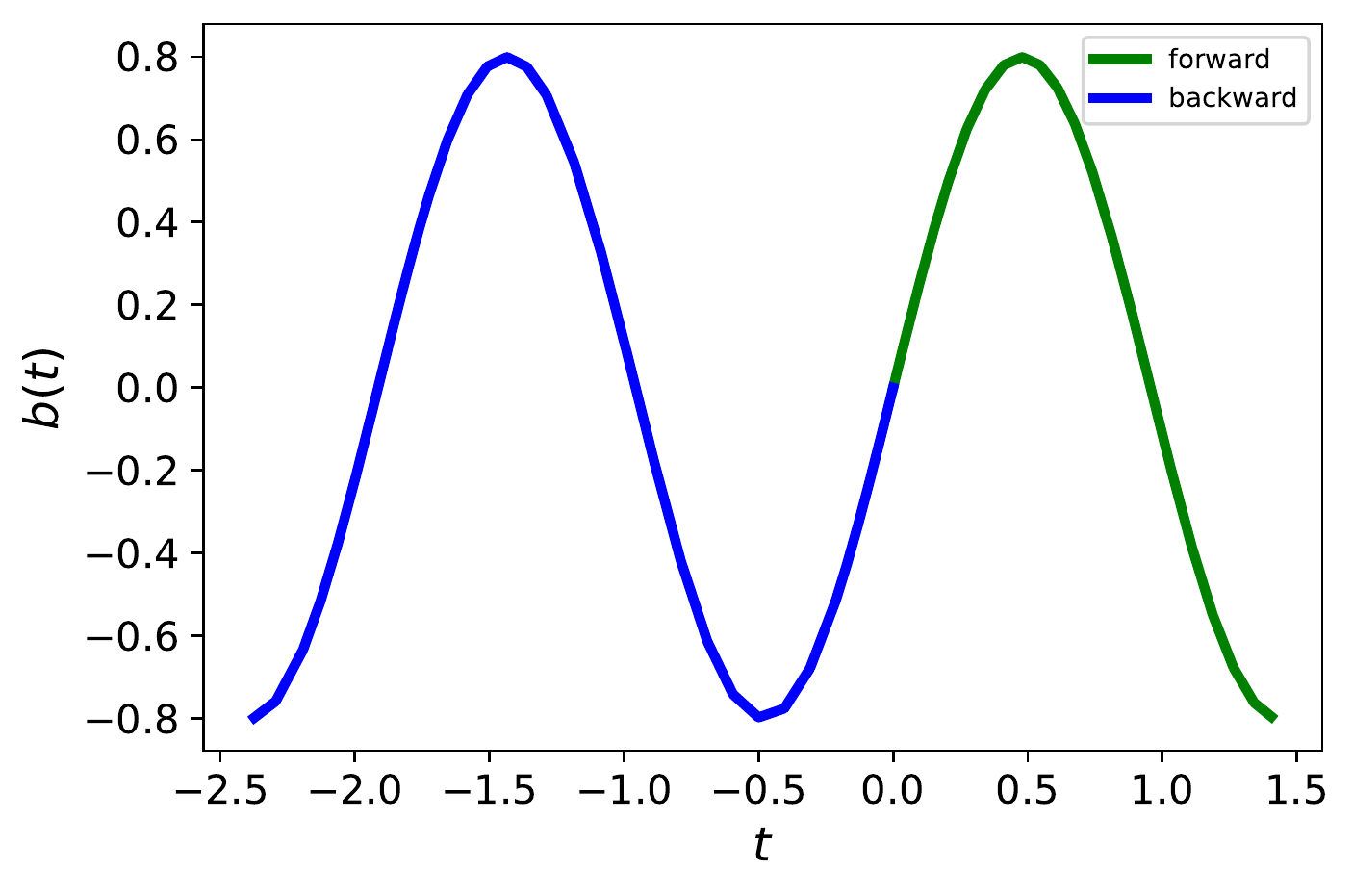}
	\caption{Arc 1 trajectory (top left), $\xi^2$ profile (top right), and optimized control functions, $a(t)$ (bottom left) and $b(t)$ (bottom right). Green and blue indicate the forward- and backward-in-time solutions, respectively.}	\label{fig:Arc1_Length_NewControl}
\end{figure}

The simulated arc 1, as well as the optimized control functions and $\xi^2$ profile, are shown in Figure \ref{fig:Arc1_Length_NewControl}.

The forward- and backward-in-time solutions are represented by the green and blue segments, respectively. As this arc will be will be joined by finite arcs with an execution of a finite turn, there is a physical requirement that the speed, $\xi^2$, be zero at either end of the arc. The top right panel shows that the speed is indeed zero at either end of the arc. The control functions, $a(t)$ and $b(t)$, that result from the optimized control parameters are shown in the bottom two panels. Note that these functions form discontinuities at the cusps, yielding finite impulse allowing for a finite turn at the infinitely short time. Such finite turn in an infinitely short time is an idealized approximation of the motion a skater makes during the actual performance. 

\begin{figure} 
	\includegraphics[scale=0.4]{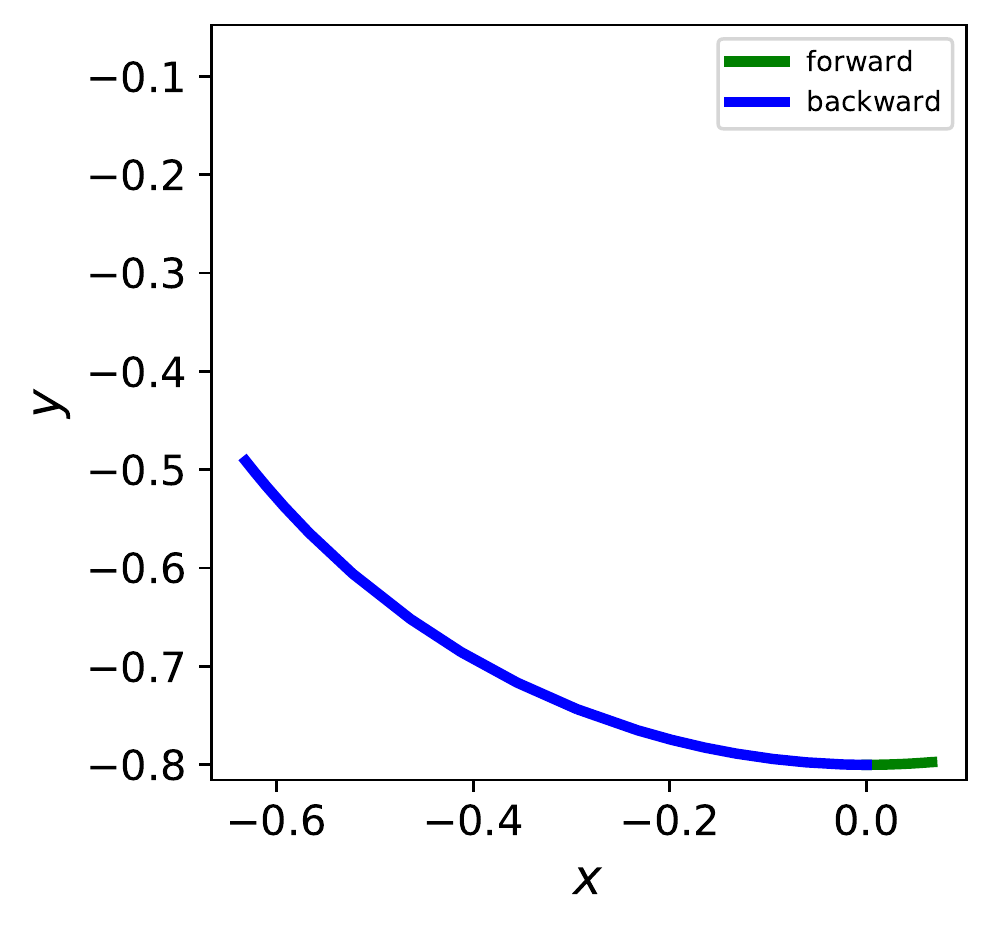}\hfill
	\includegraphics[scale=0.4]{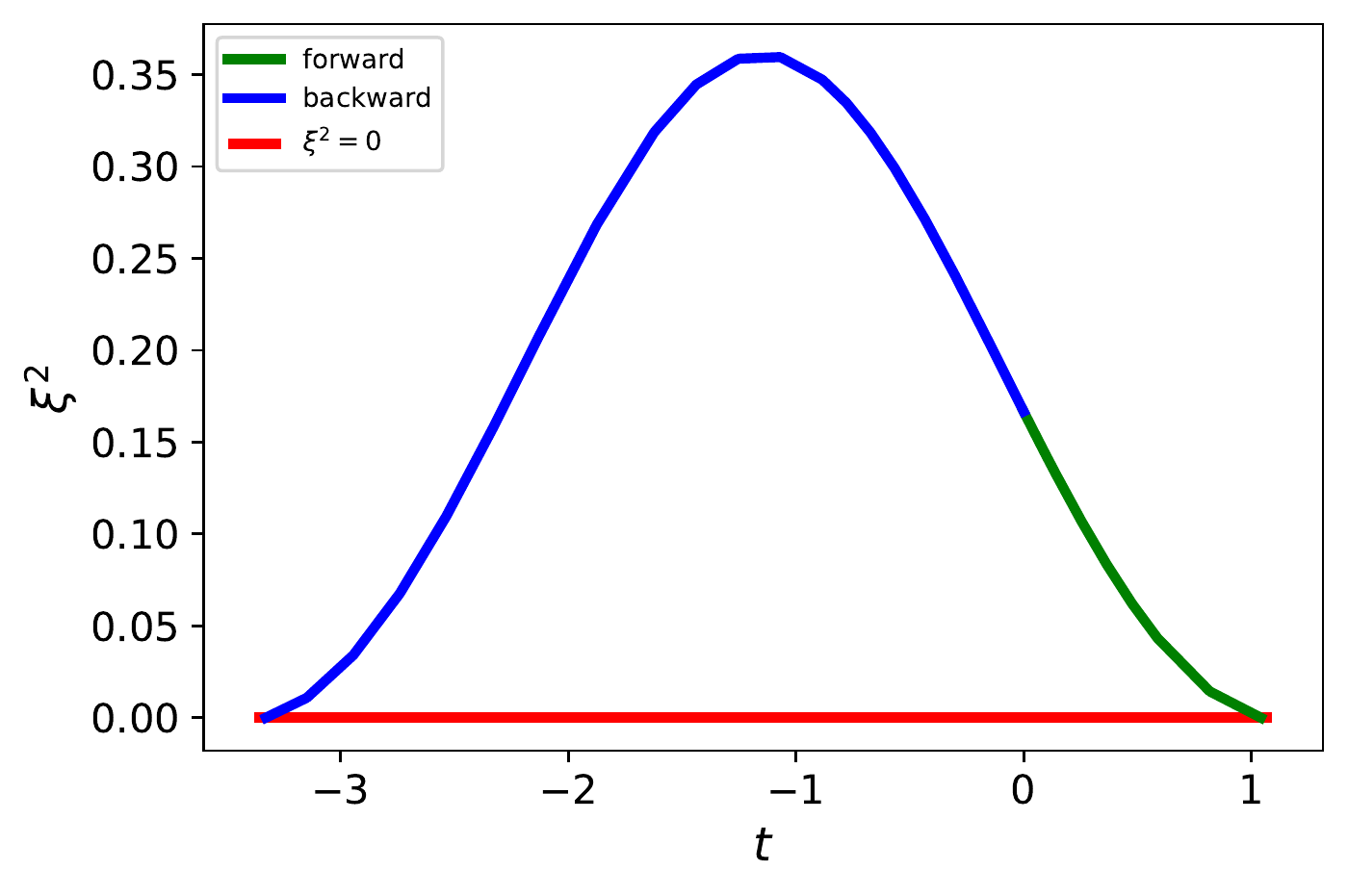}\hfill
	\includegraphics[scale=0.4]{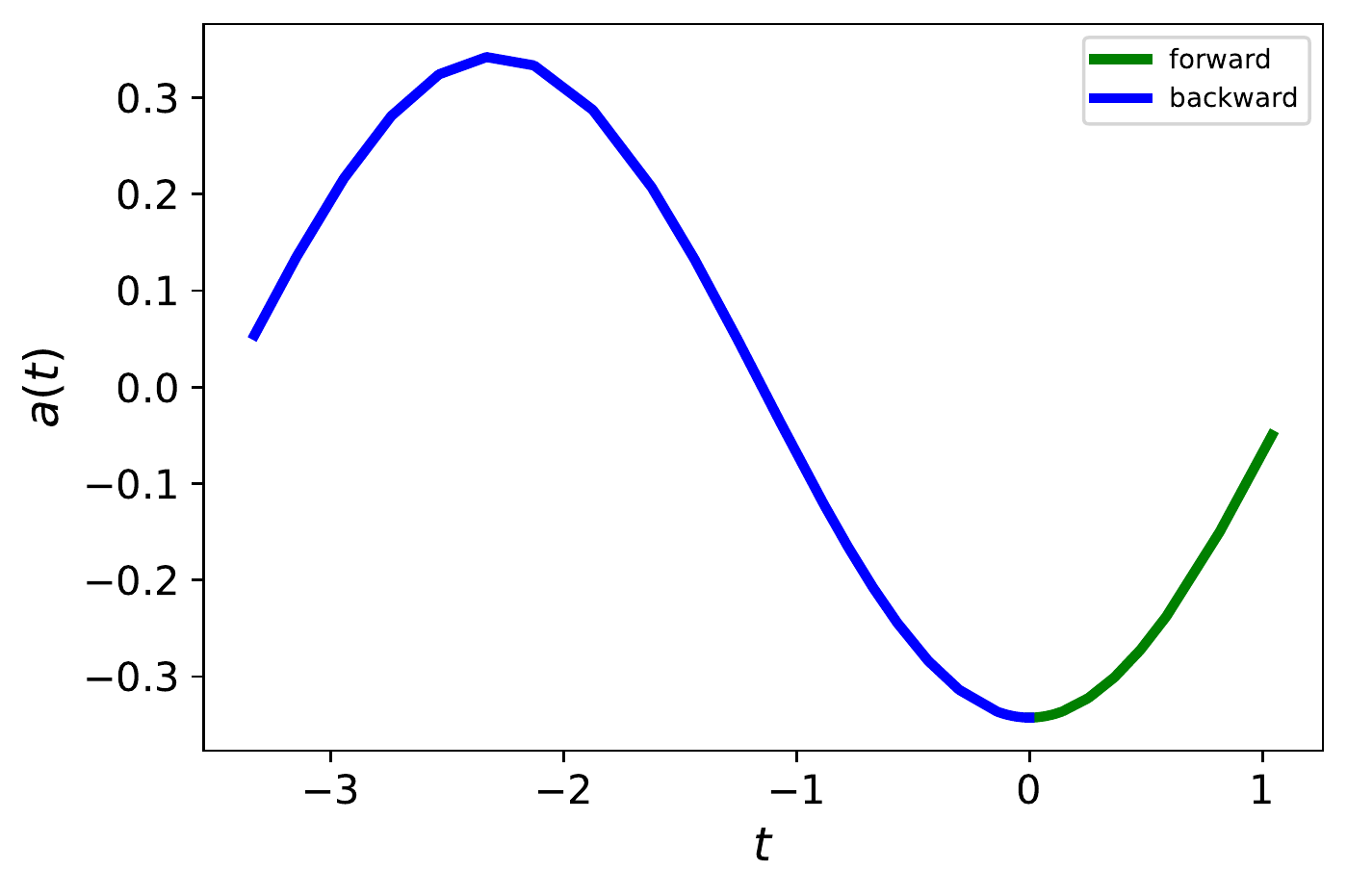}\hfill
	\includegraphics[scale=0.4]{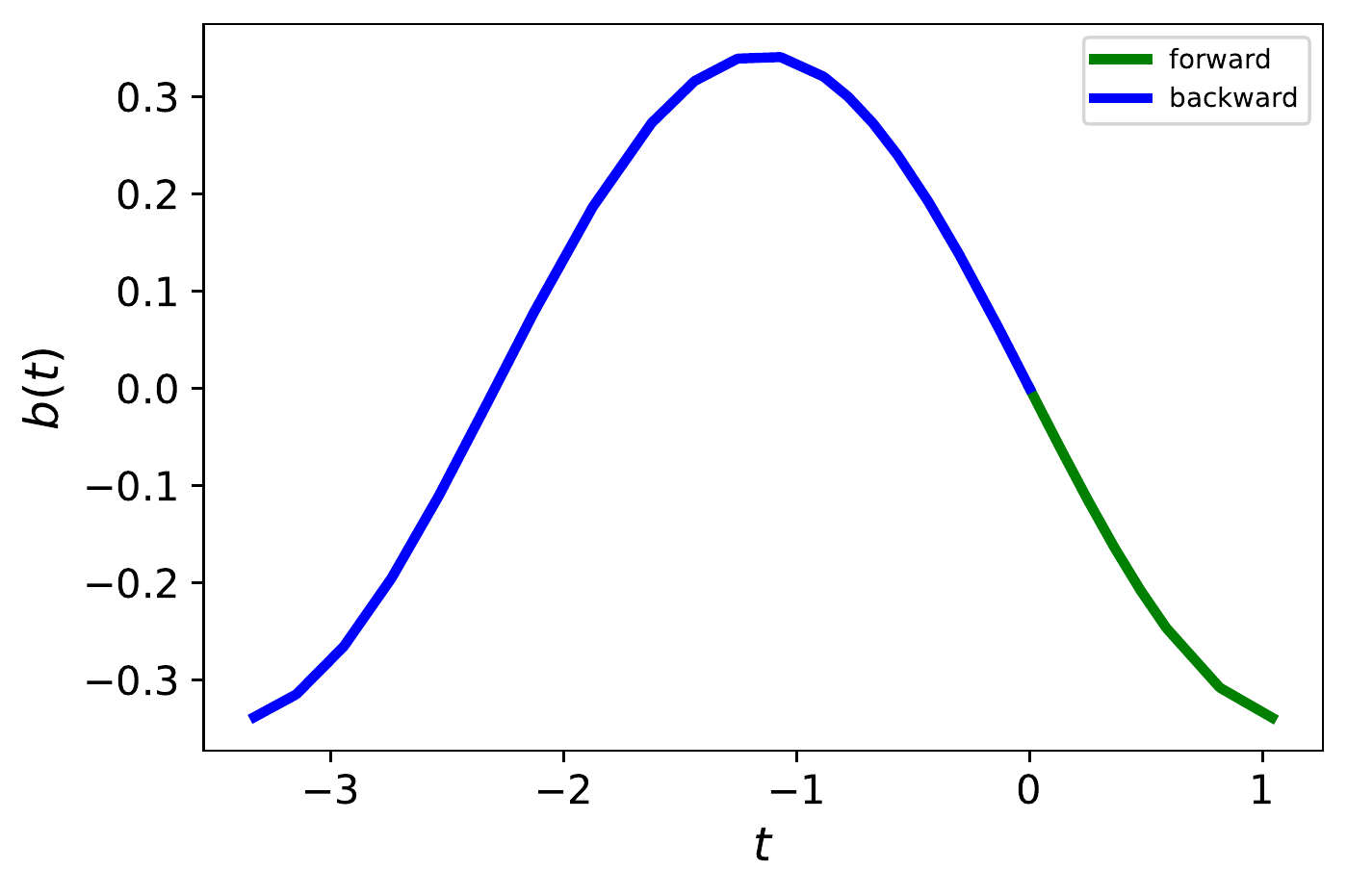}
	\caption{Arc 2 trajectory (top left), $\xi^2$ profile (top right), and optimized control functions, $a(t)$ (bottom left) and $b(t)$ (bottom right). Green and blue indicate the forward- and backward-in-time solutions, respectively.}	\label{fig:Arc2_Length_NewControl}
\end{figure}

\begin{figure} 
	\includegraphics[scale=0.4]{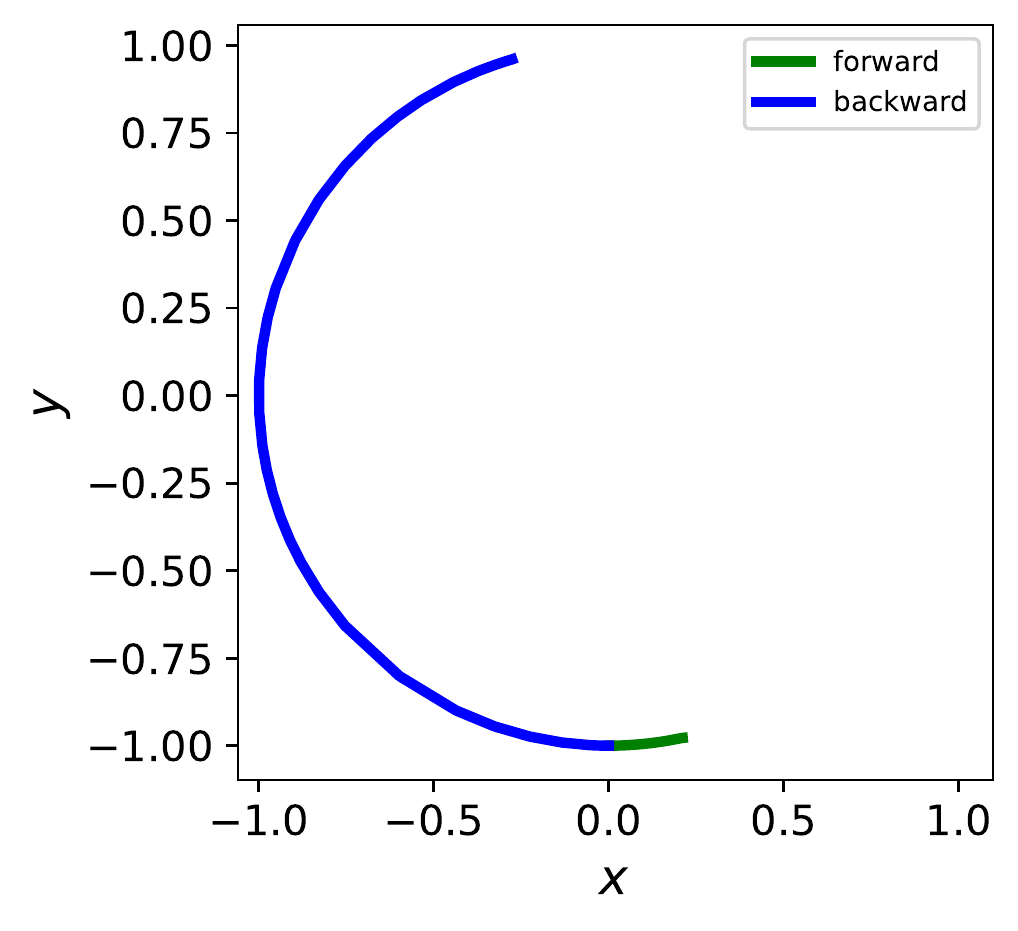}\hfill	\includegraphics[scale=0.4]{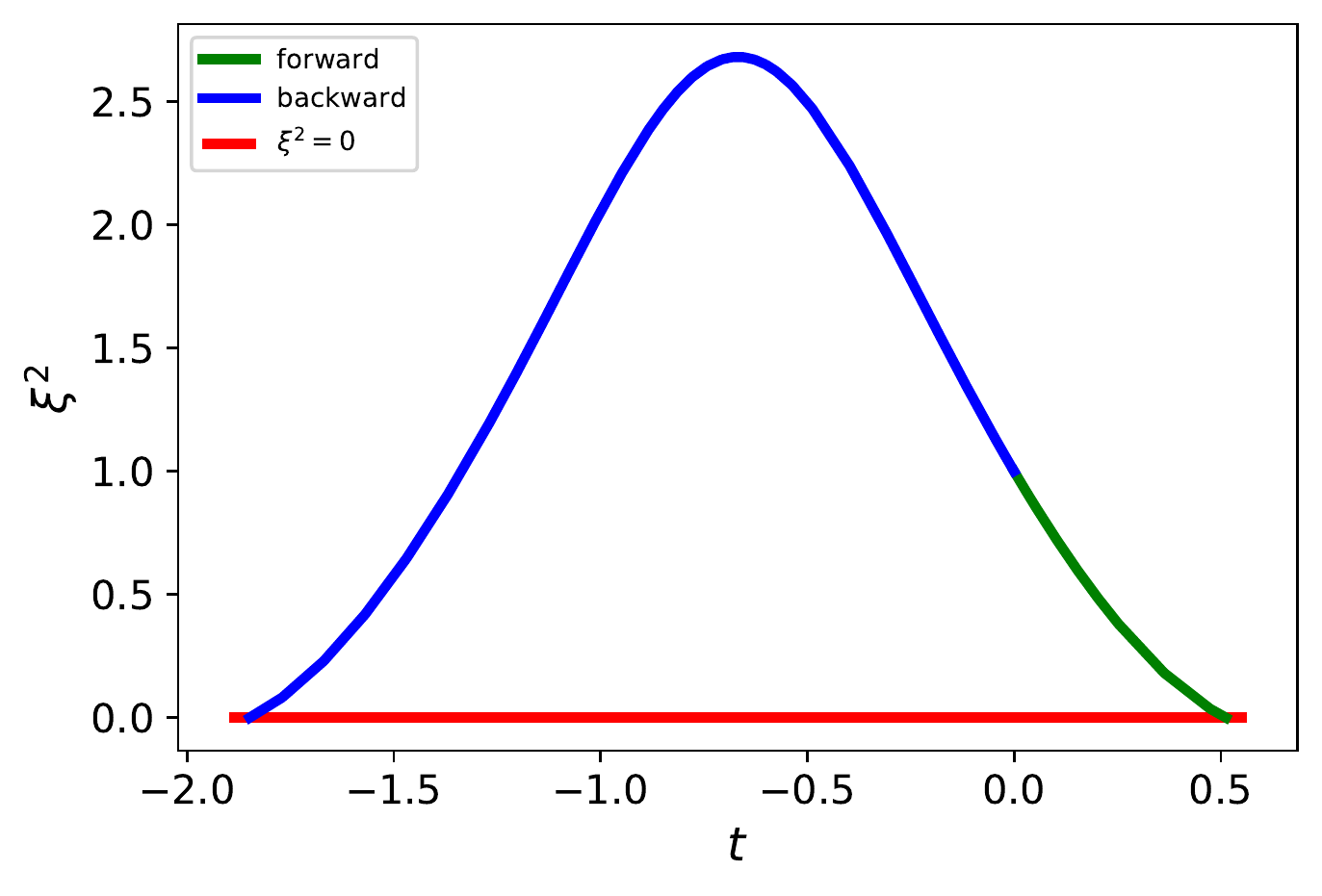}\hfill	\includegraphics[scale=0.4]{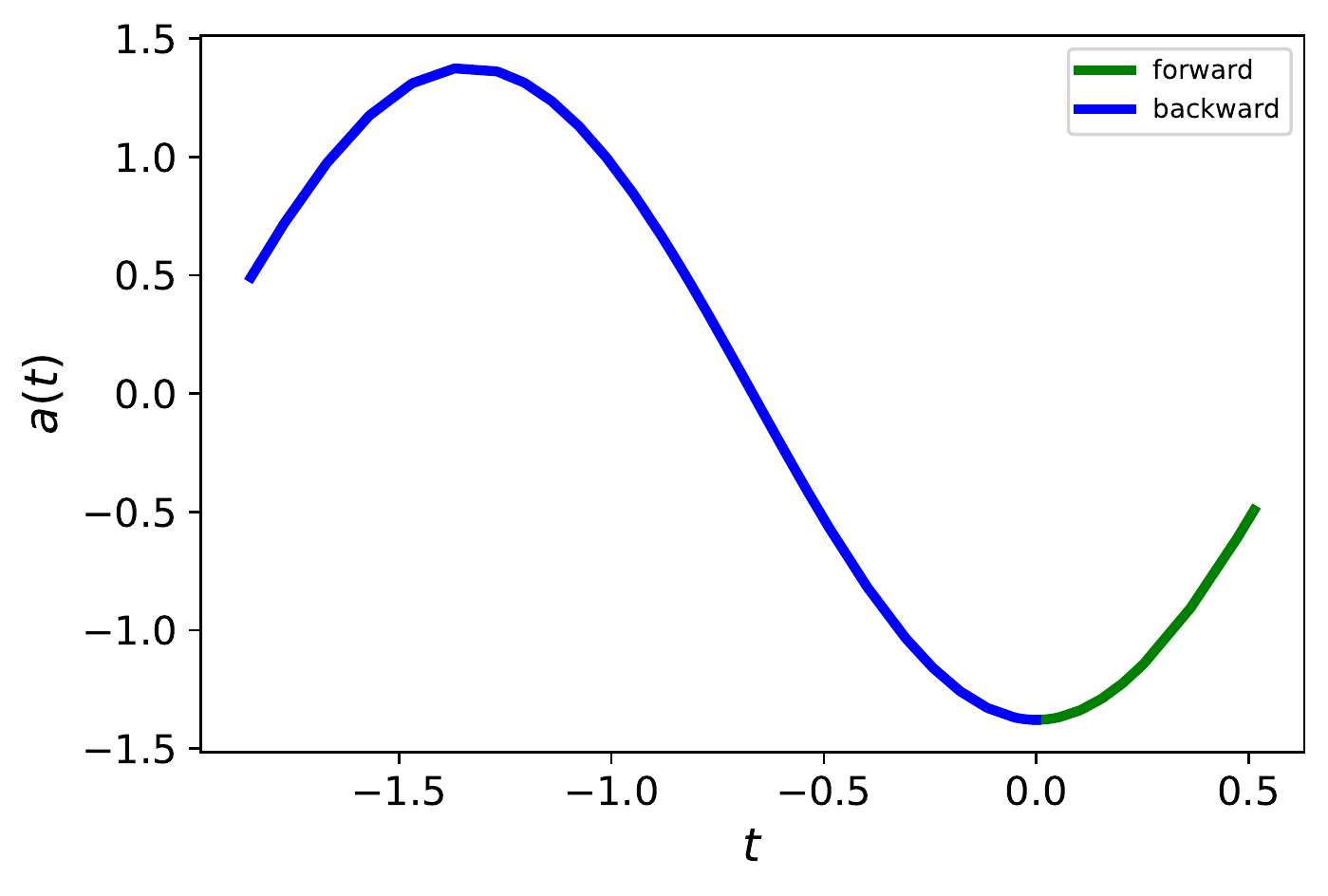}\hfill
	\includegraphics[scale=0.4]{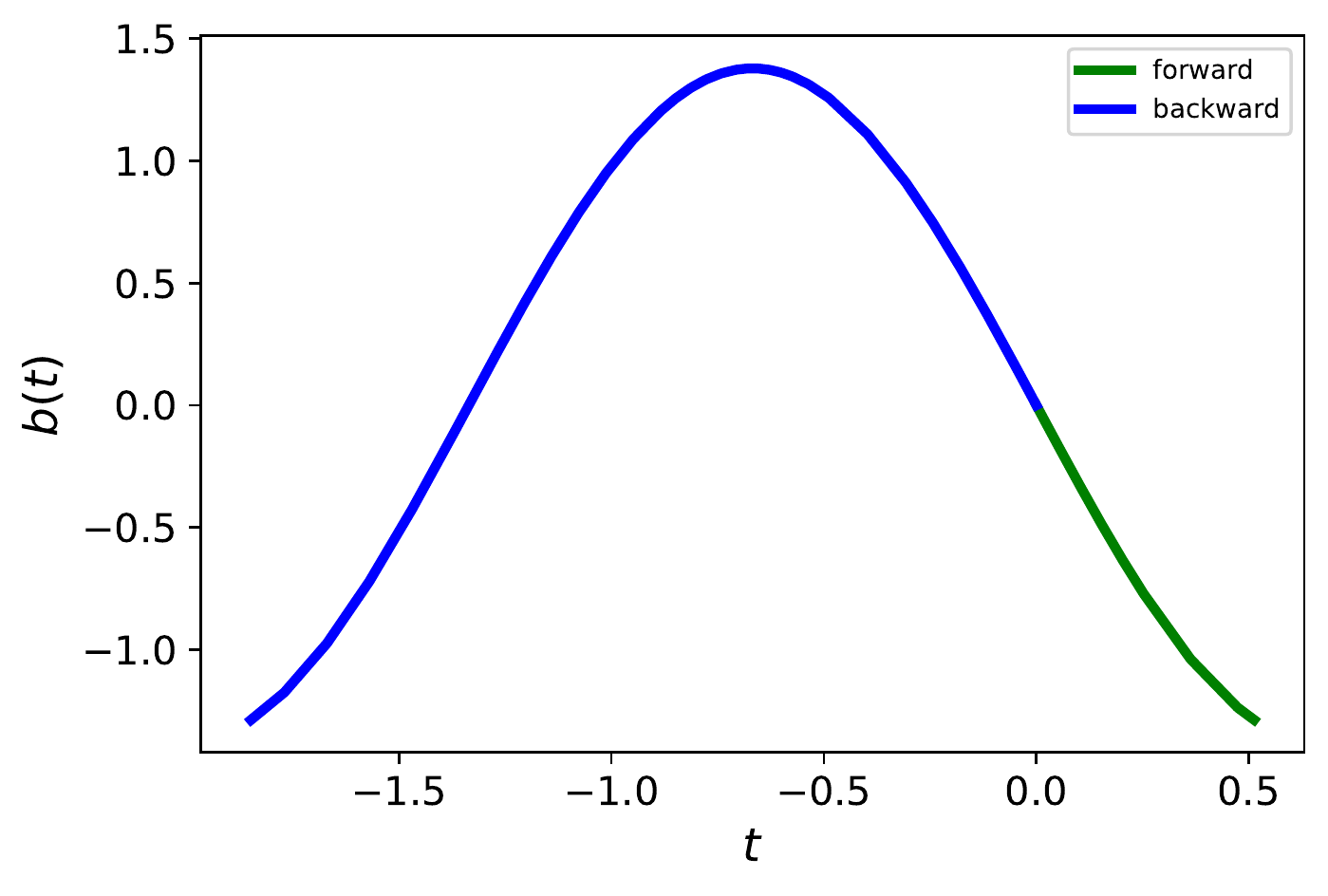}
	\caption{Arc 3 trajectory (top left), $\xi^2$ profile (top right), and optimized control functions, $a(t)$ (bottom left) and $b(t)$ (bottom right). Green and blue indicate the forward- and backward-in-time solutions, respectively.}	\label{fig:Arc3_Length_NewControl}
\end{figure}

\begin{figure} 
	\centering
	\includegraphics[scale=0.4]{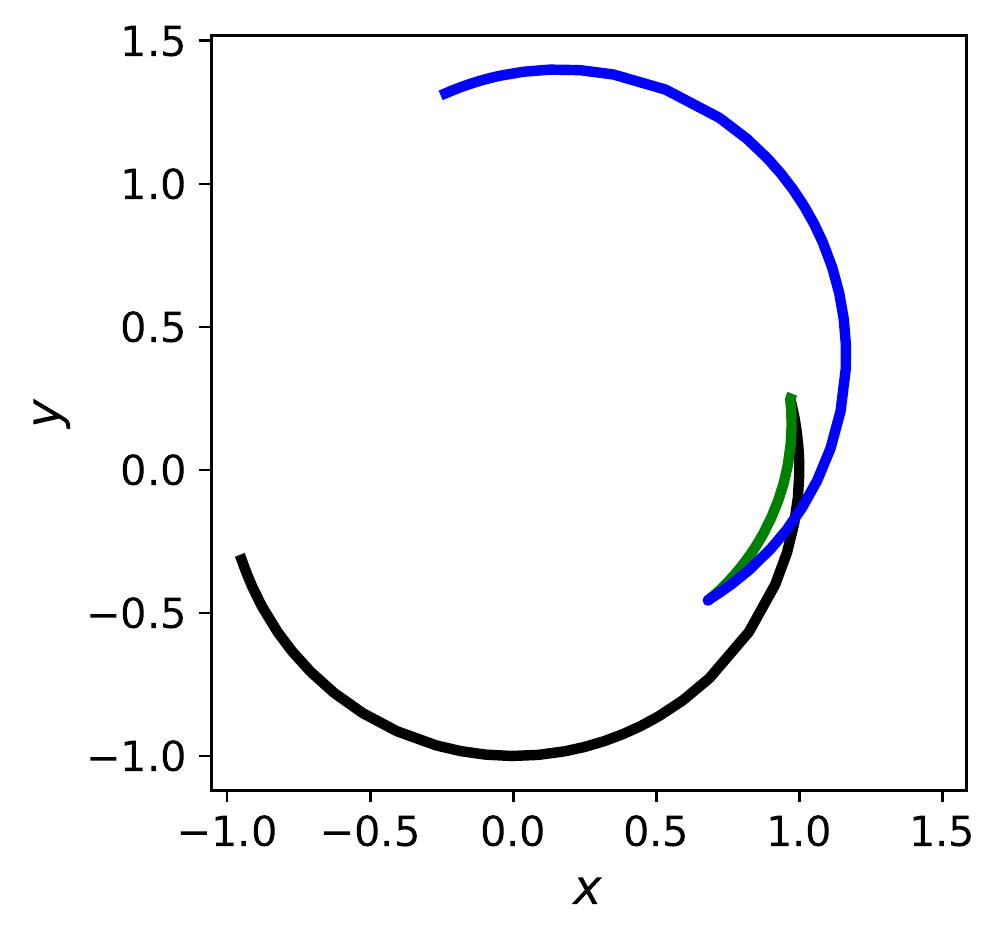}\hfill
	\includegraphics[scale=0.4]{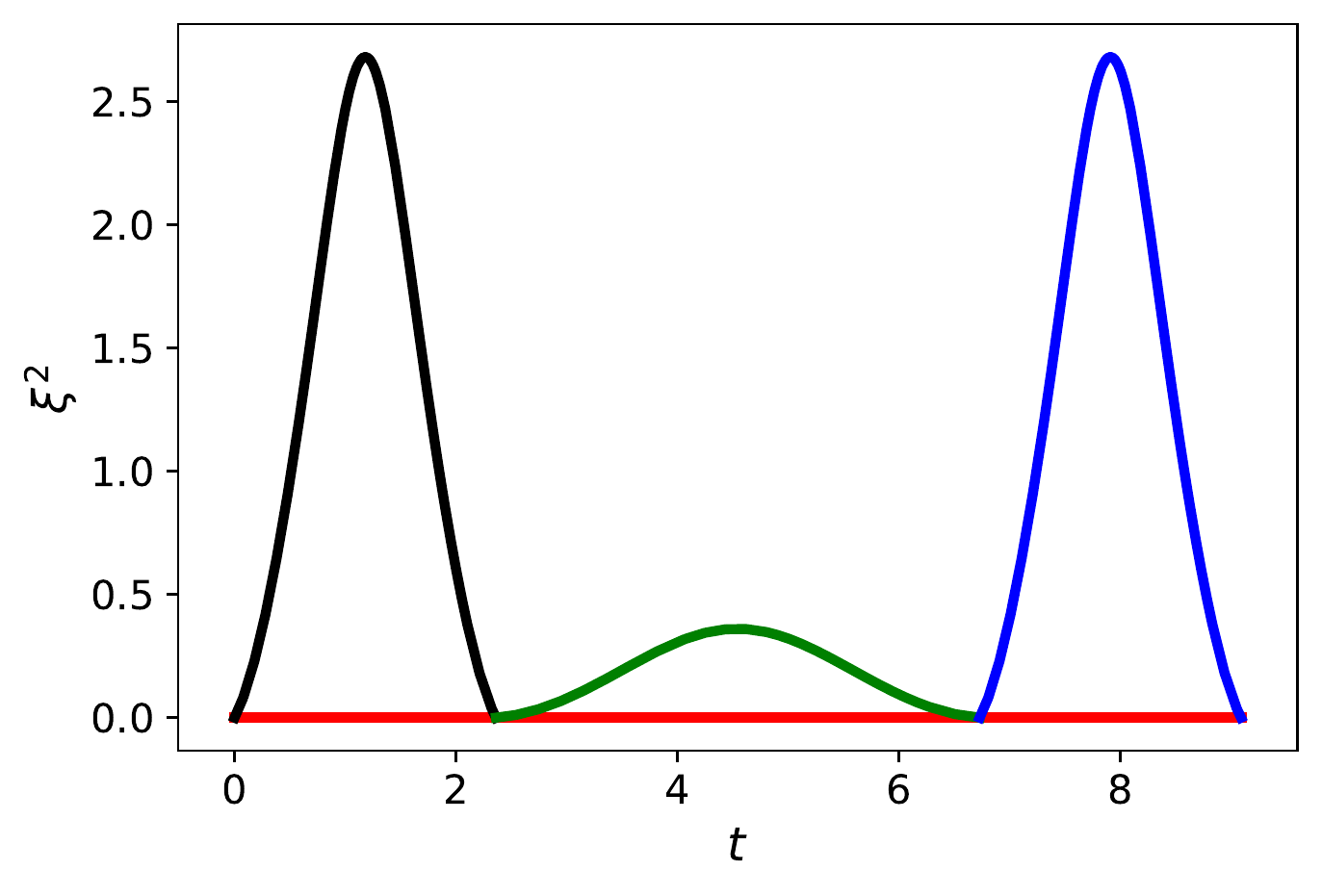}\hfill
	\includegraphics[scale=0.4]{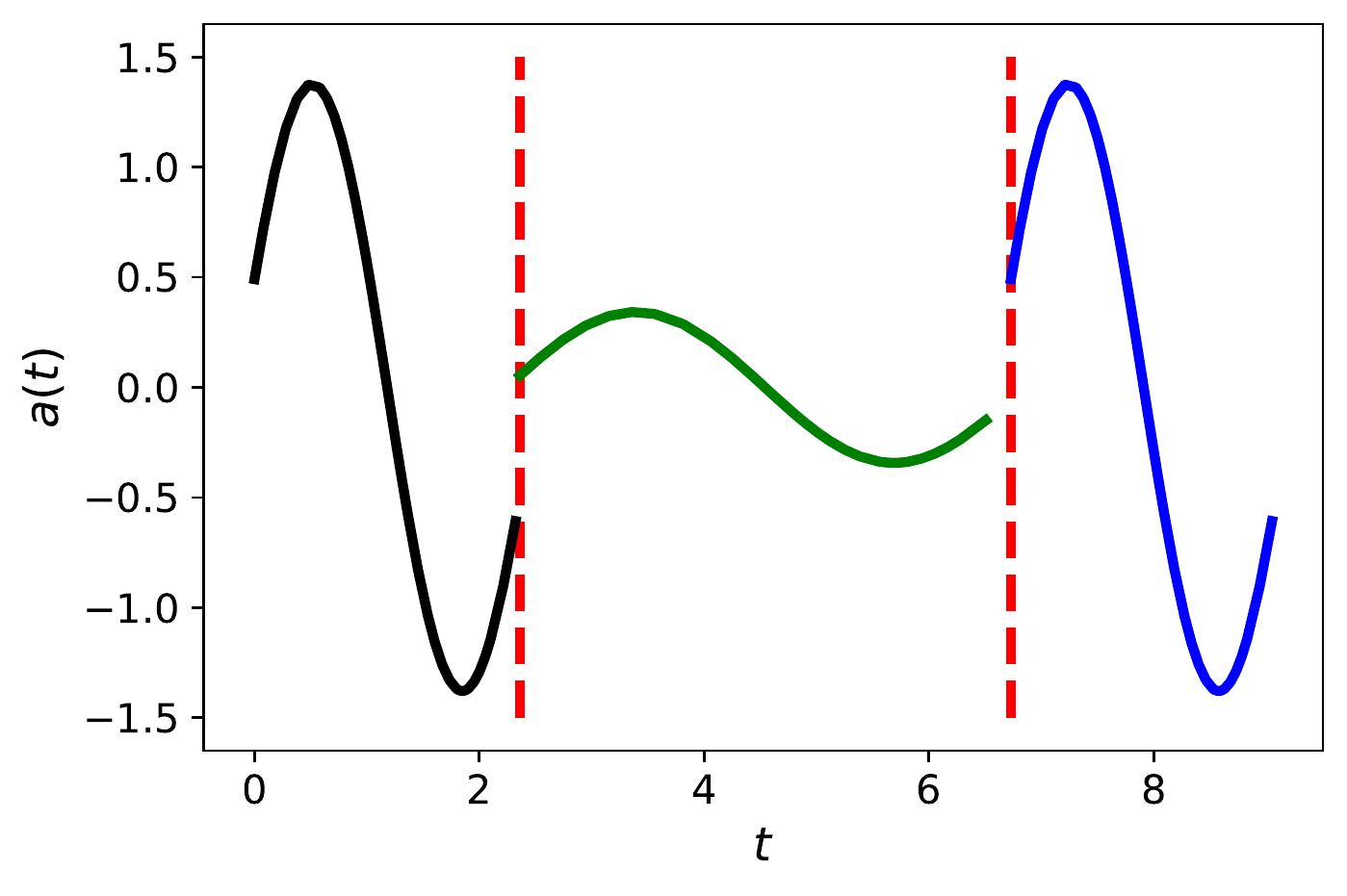}\hfill
	\includegraphics[scale=0.4]{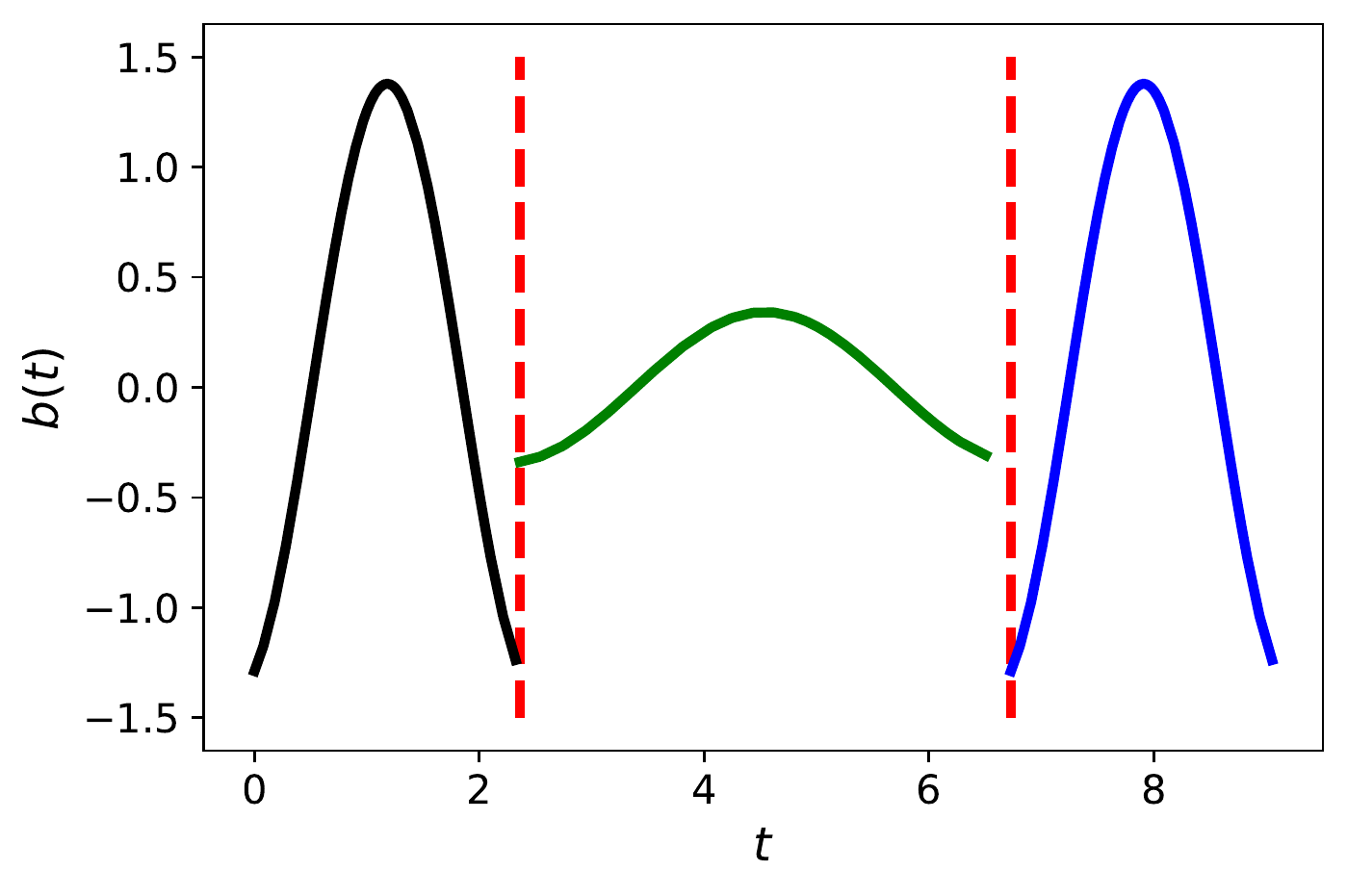}
	\caption{Trajectory of a single leaf (top left), corresponding $\xi^2$ profile (top right) and optimized control functions $a(t)$ (bottom left) and $b(t)$ (bottom right). The black, green, and blue segments distinguish the three arcs that make up the triplet. The vertical red lines indicate the time points where the controls switch between arc types.}
	\label{fig:Triplet_Length_NewControl}
\end{figure}

Similarly, the simulations for arc 2 and arc 3, as well as the optimized control functions and profiles of $\xi^2$, are shown in Figures \ref{fig:Arc2_Length_NewControl} and \ref{fig:Arc3_Length_NewControl}. As noted previously, arc 2 and arc 3 are combined to form the outer portion of the target pattern. Each ``leaf'' of the outer pattern is formed by a triplet consisting of a single arc 2 connected to an arc 3 at both ends. This leaf is then repeated to give the full outer portion. The connections of arc 2 and arc 3, as well as arc 3 to arc 3 (from the connection of leaves/triplets), results in cusps at each of these connection points. Therefore, the velocity at each end of both arc 2 and arc 3 must be zero, as is shown in the top left panels of Figures \ref{fig:Arc2_Length_NewControl} and \ref{fig:Arc3_Length_NewControl}.
The combined triplet simulation, optimized control functions, and corresponding $\xi^2$ profile are shown in Figure \ref{fig:Triplet_Length_NewControl}.

We note that the control functions need not be continuous when switching from arc 2 to arc 3, and vice versa, as sudden shifts in a skater's center of mass is not an unreasonable allowance. 

\begin{figure} 
	\centering
	\includegraphics[scale=0.38]{Figures/June4 Results New Controls Small Files/Arc1}	\includegraphics[scale=0.38]{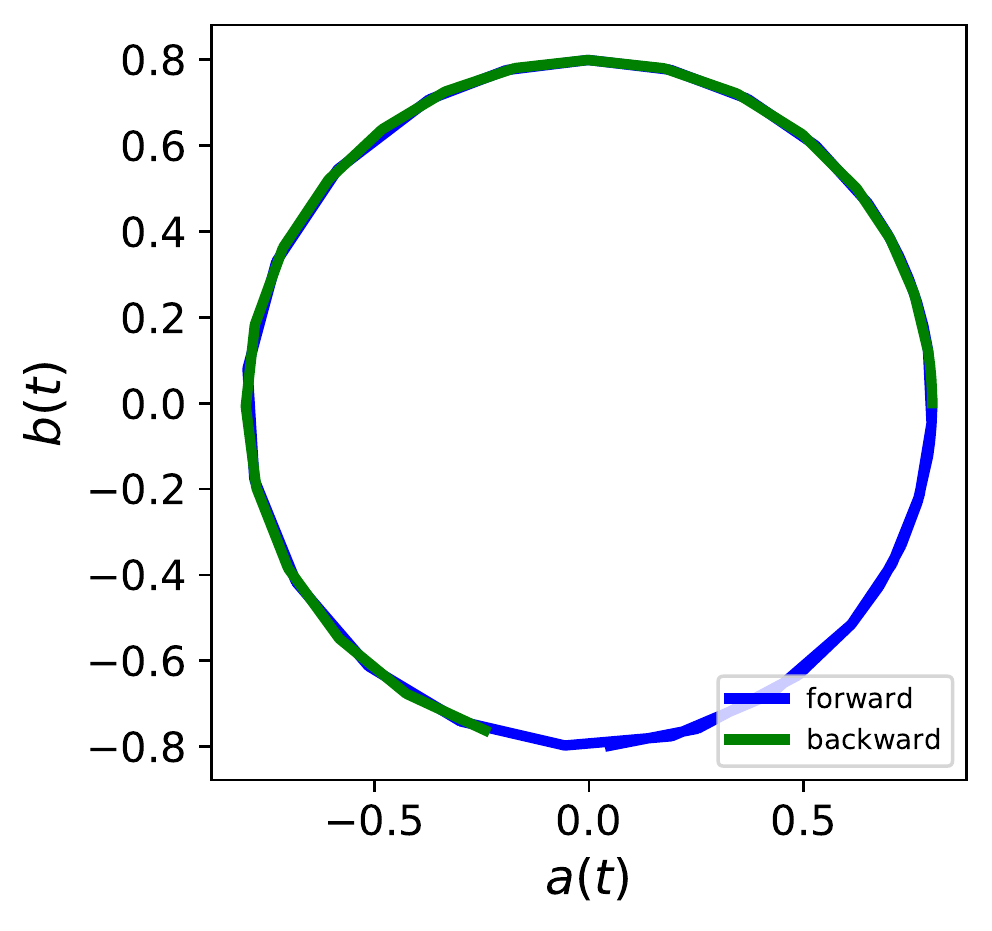} 
	\includegraphics[scale=0.38]{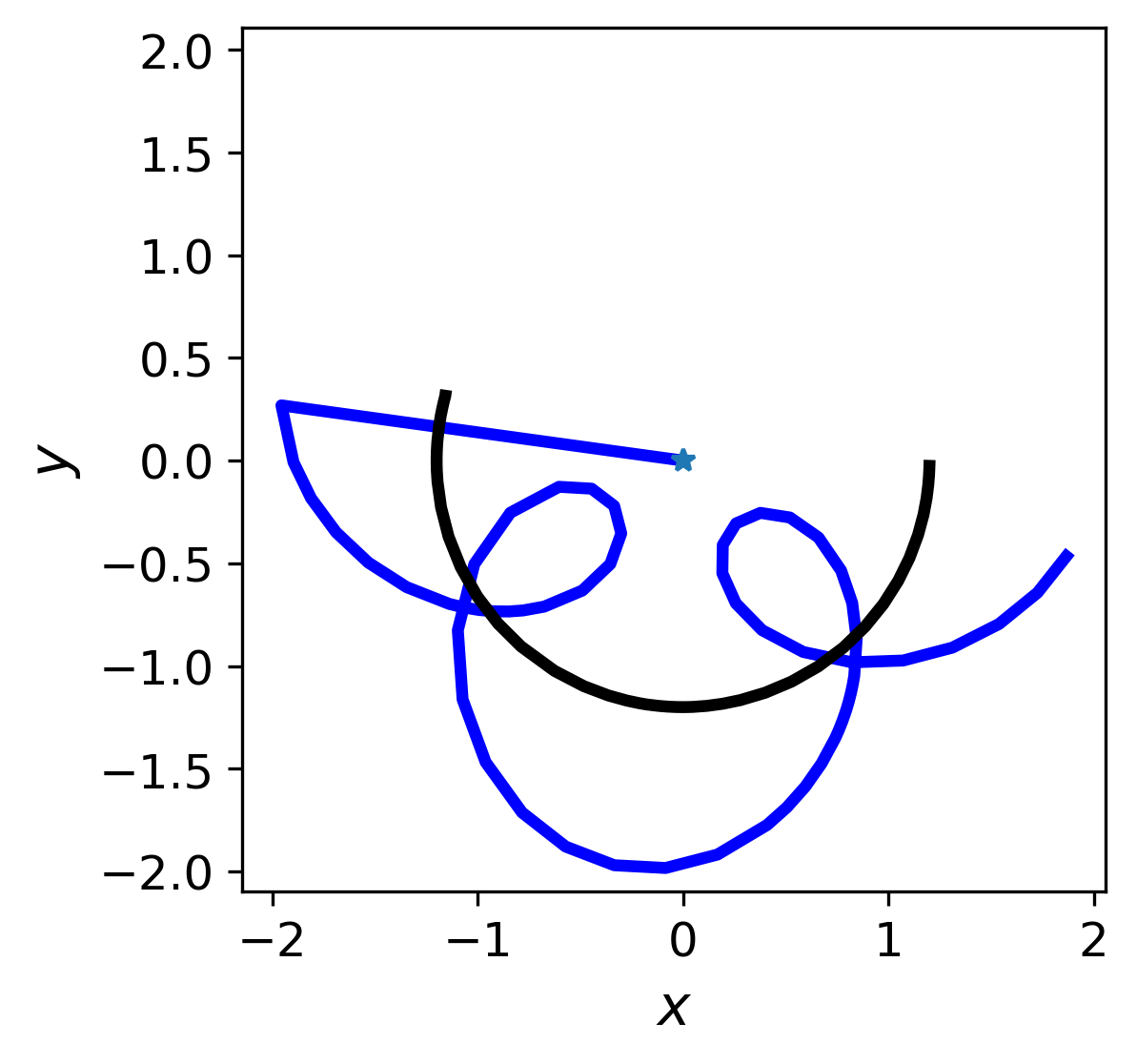}\\
	\includegraphics[scale=0.38]{Figures/June4 Results New Controls Small Files/Arc2}	\includegraphics[scale=0.38]{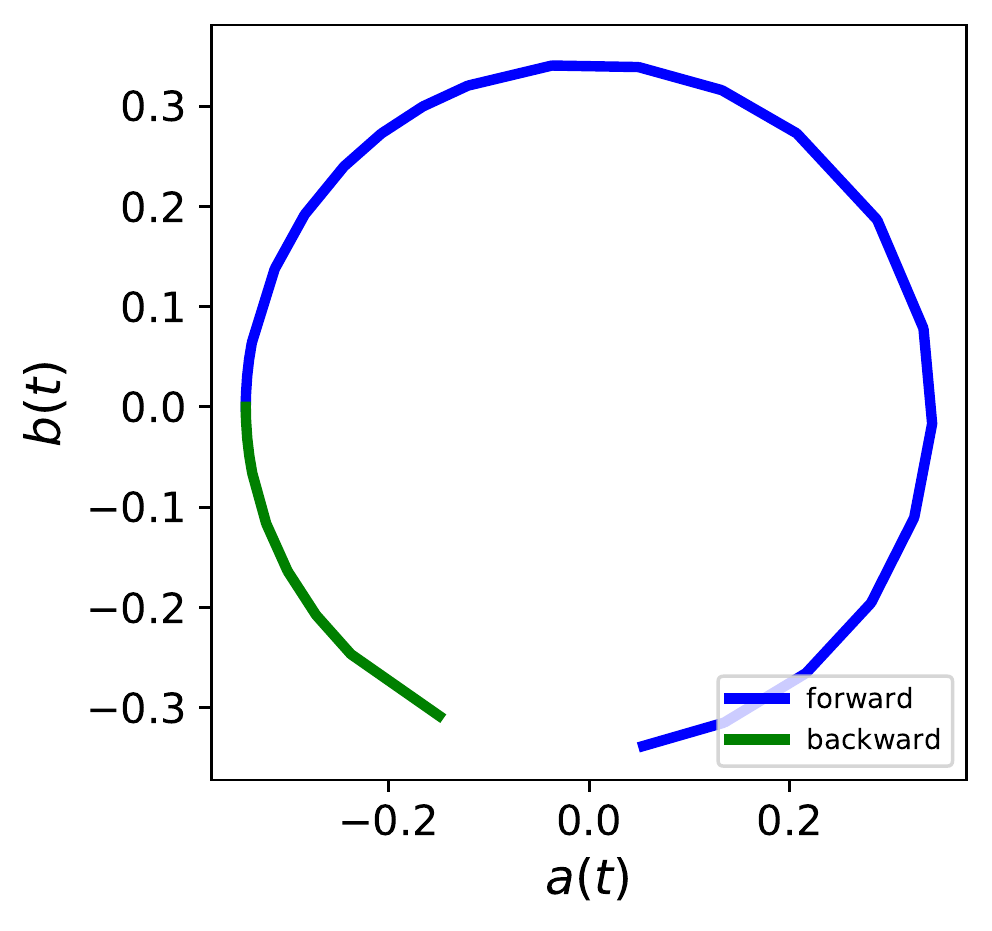} 
	\includegraphics[scale=0.38]{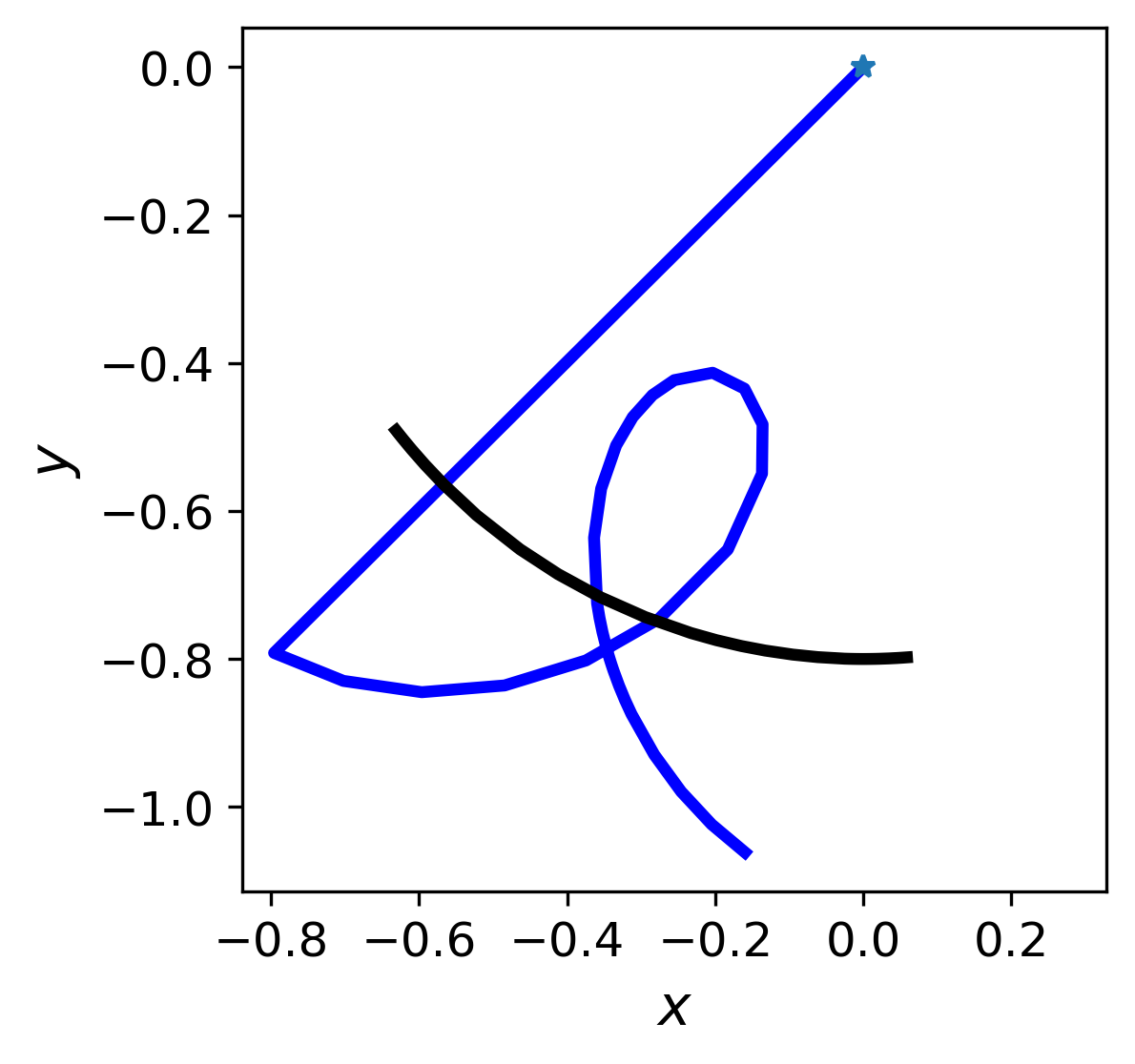}\\
	\includegraphics[scale=0.38]{Figures/June4 Results New Controls Small Files/Arc3}	\includegraphics[scale=0.38]{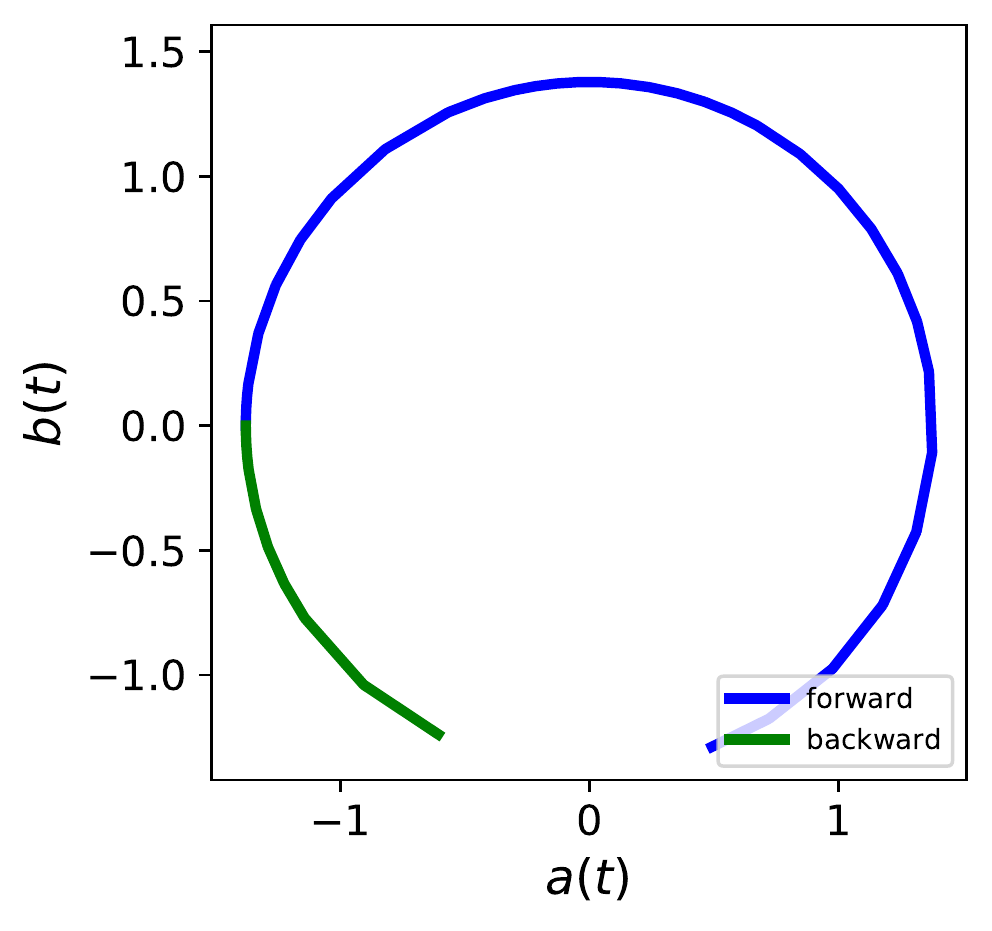} 
	\includegraphics[scale=0.38]{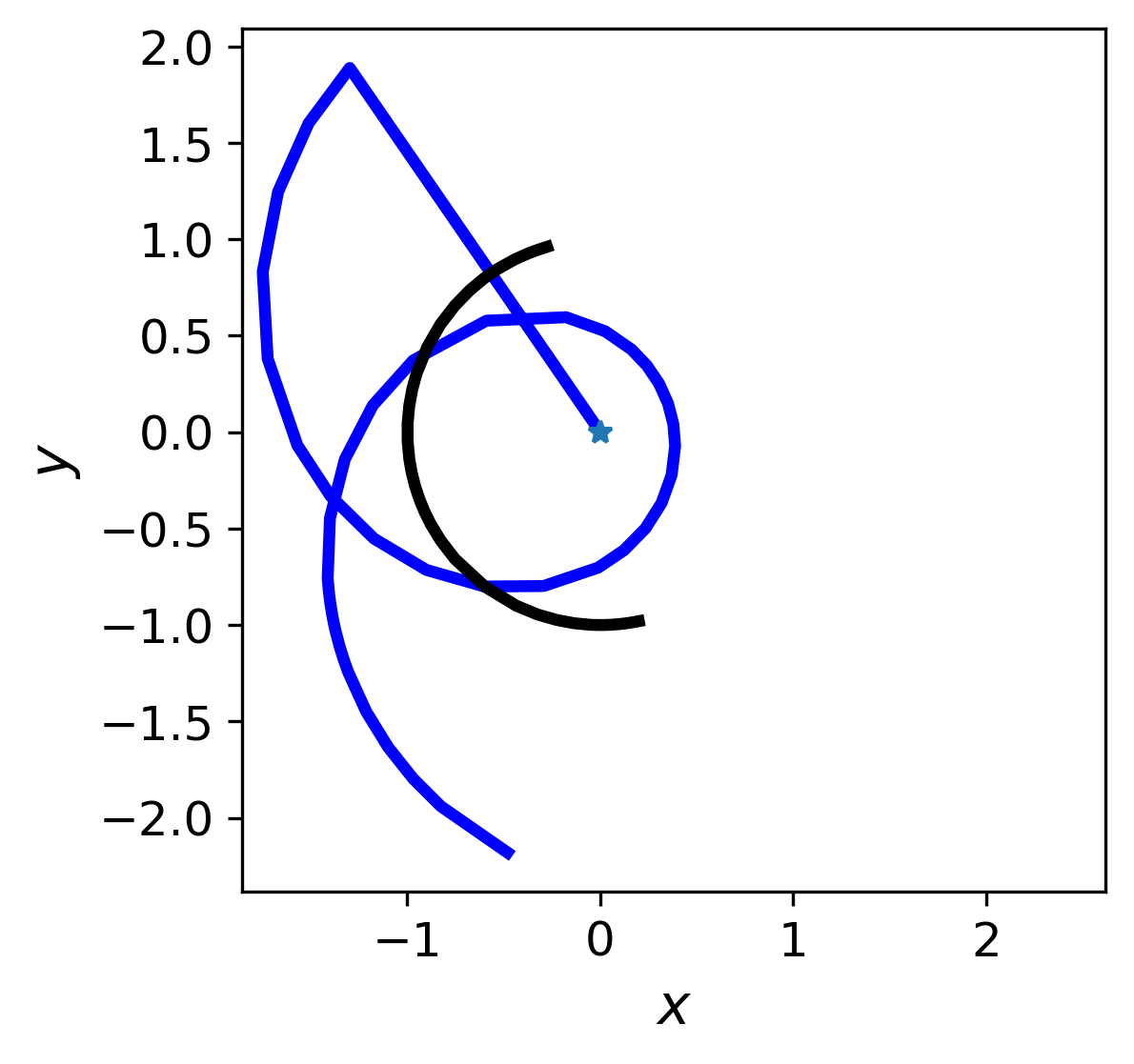}
	\caption{Optimized skate trajectories (left), control mass trajectories (middle), and overlay of skate and control mass trajectories in spatial frame (right) for arc 1 (top), arc 2 (middle), and arc 3 (bottom). Green and blue indicate the forward- and backward-in-time solutions, respectively. In the far right column, the star indicates the initial point of the control mass.}
	\label{fig:ArcControls_Length_NewControl}
\end{figure}

In order to better understand the optimized control functions, Figure \ref{fig:ArcControls_Length_NewControl} shows the trajectory of the position of the added mass $m$ (in body frame) next to the blade trajectory (in spatial frame) for each arc, as well as the skate and mass trajectories together in the spatial frame. Again, the green indicates the forward-in-time solution, while the blue indicates the backward-in-time solution. Thus, the path starts at the beginning of the blue portion of the curve and ends at the end of the green portion of the curve. 

\begin{figure} 
\begin{minipage}{0.49\textwidth}
    \begin{flushright}
        \includegraphics[scale=0.4]{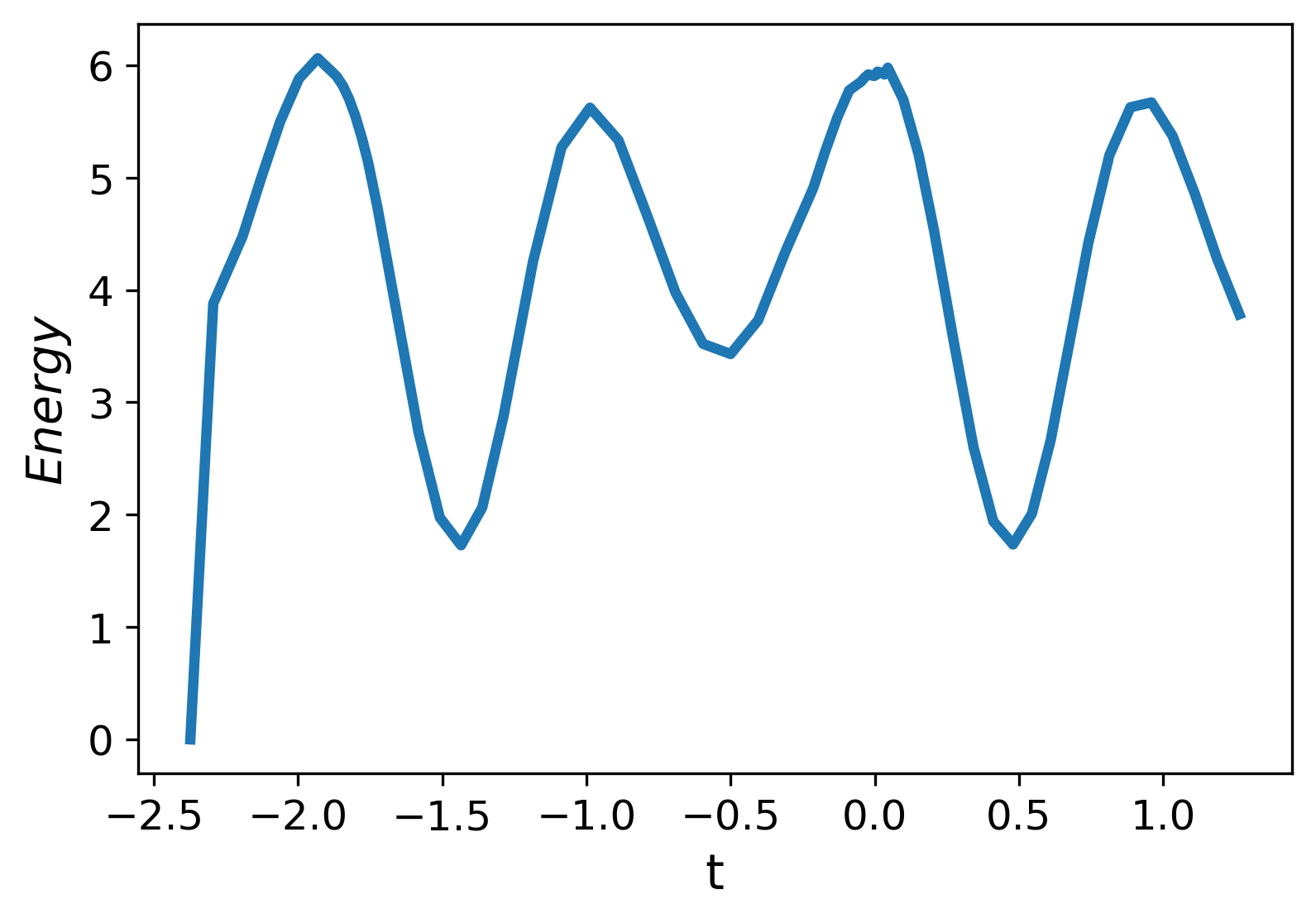}\\	\includegraphics[scale=0.4]{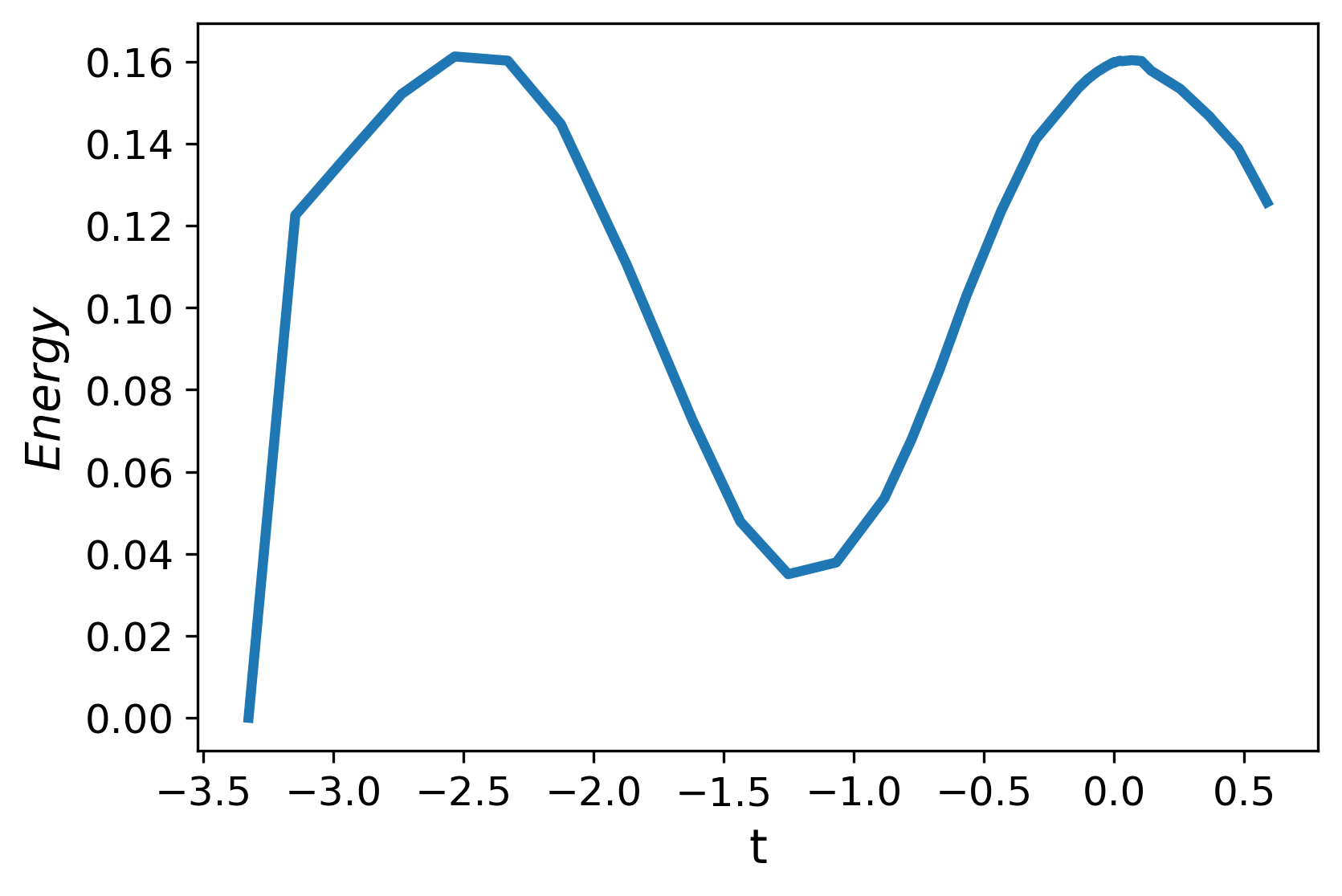}\\	
        \includegraphics[scale=0.4]{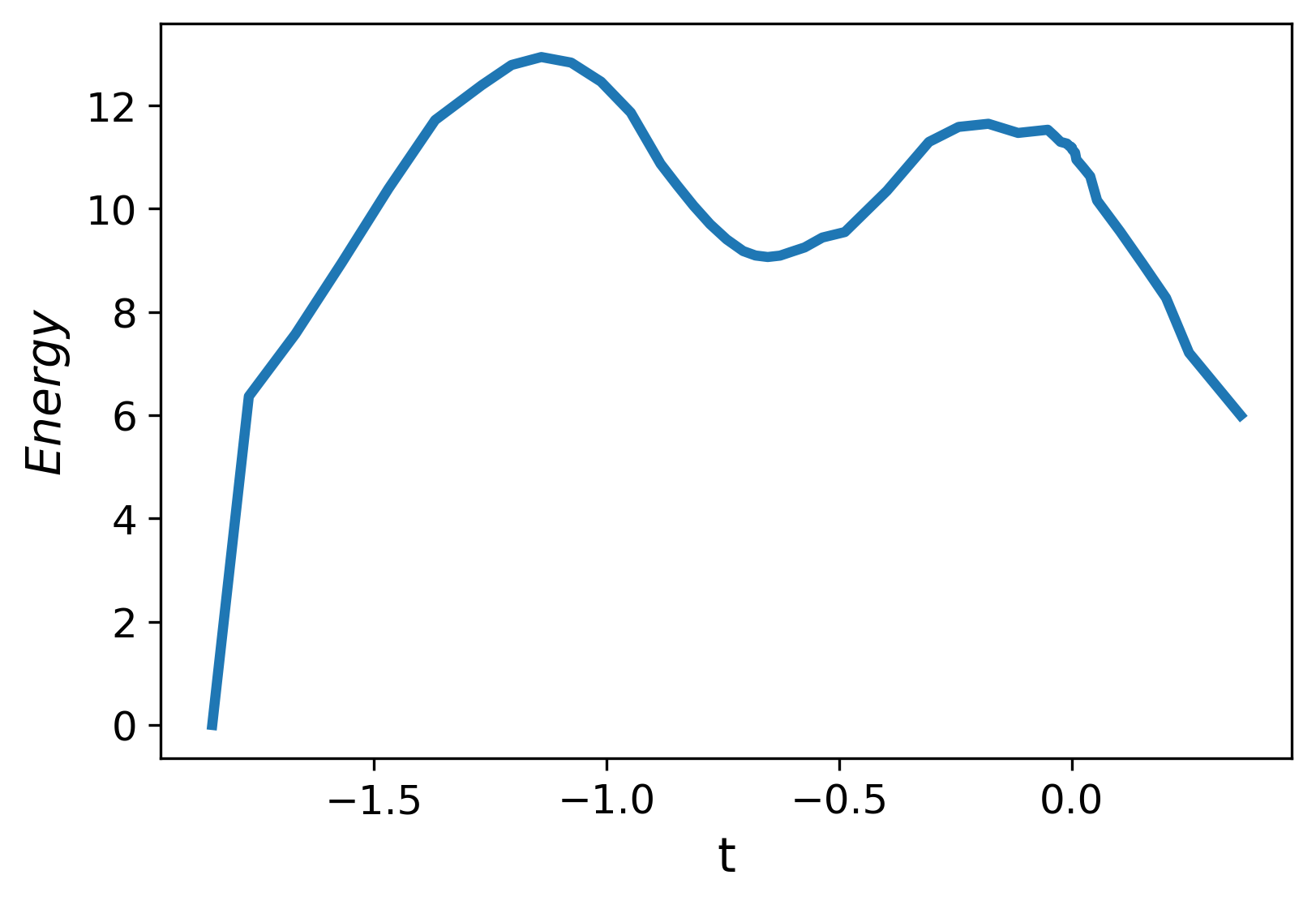}
    \end{flushright} 
    \end{minipage}
    \begin{minipage}{0.49\textwidth}
    \begin{flushright}
    	\includegraphics[scale=0.4]{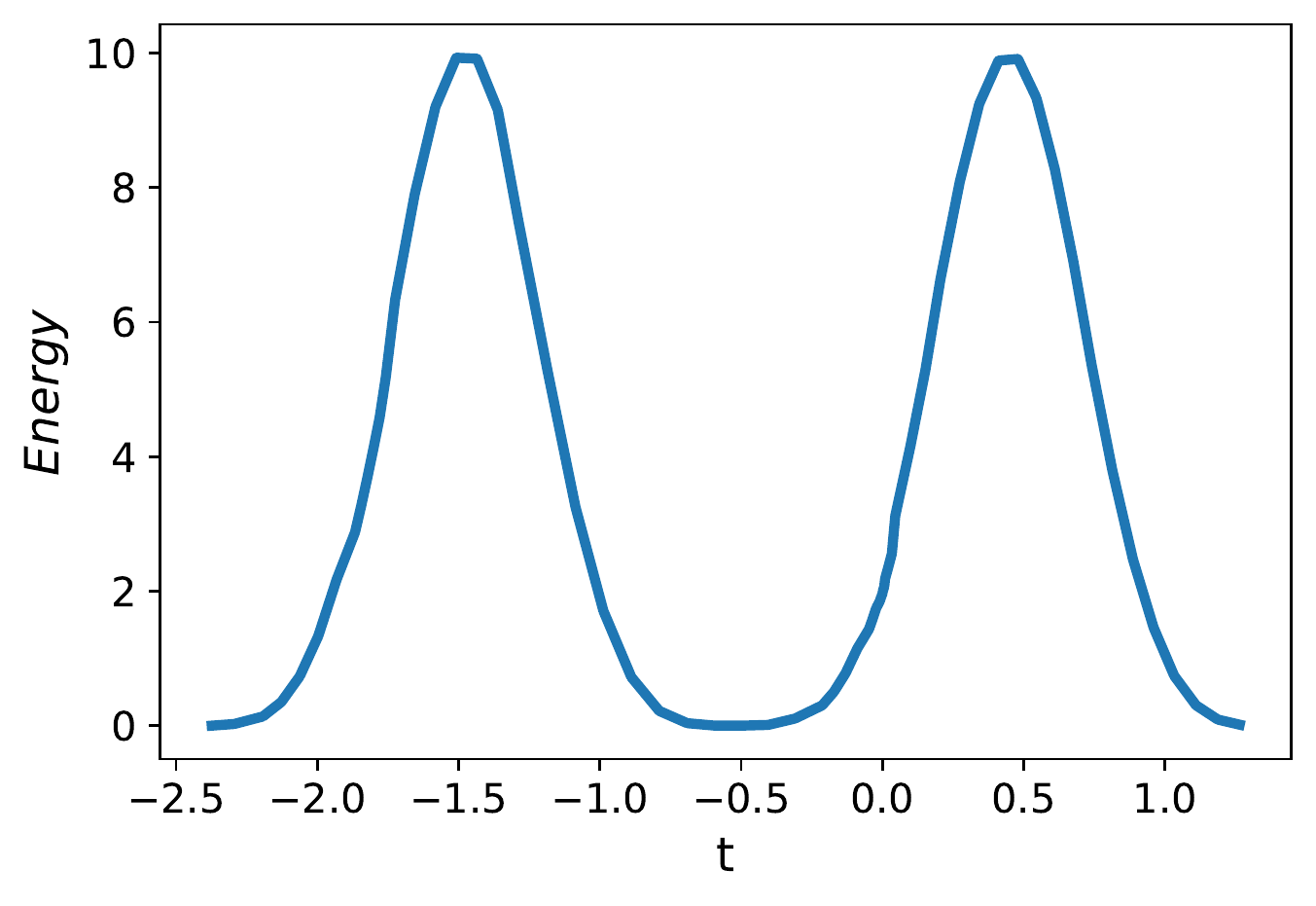}\\
        \includegraphics[scale=0.4]{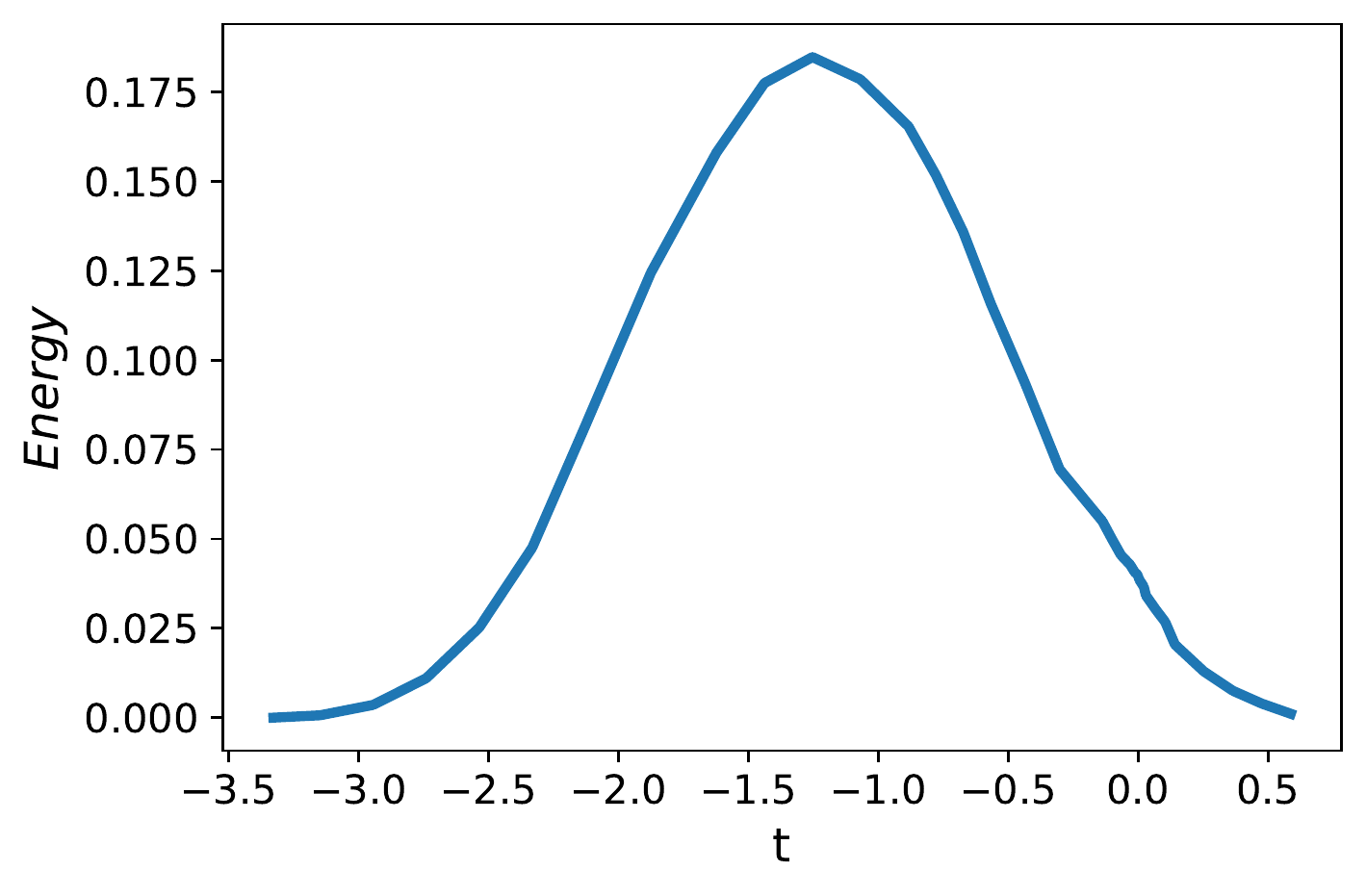}\\	\includegraphics[scale=0.4]{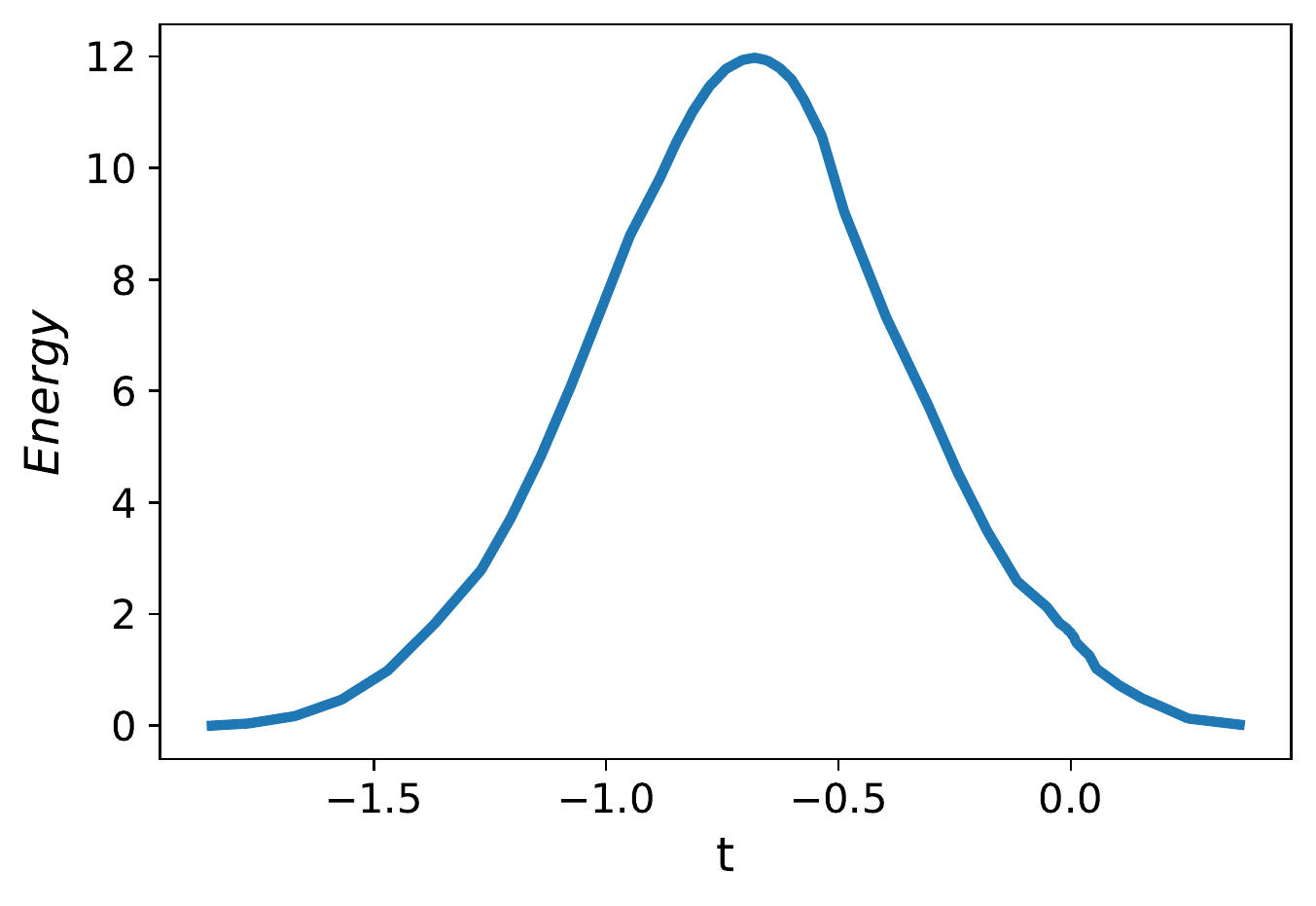}
    \end{flushright} 
    \end{minipage}
	\caption{Energy over time for control mass (left) and skate (right) for arc 1 (top), arc 2 (middle), and arc 3 (bottom).}
	\label{fig: Energy_NewControl}
\end{figure}

It is also useful to compute the energy of the system and confirm that the total energy remains bounded. Thus, we derive the total energy for both the skate and control mass, and compute the total energy for each arc.  
The total energy of the skate is given by 
\begin{align}
KE_{{\rm skate}}=\frac{1}{2}M v^2+\frac{1}{2}I\Omega^2\,,
\end{align}
where $\Omega$ is the angular velocity.
Given a circular trajectory, the angular momentum ($p_1$) can be written as 
\begin{align}
p_1=I\Omega =M v r\,,
\end{align}
where $r$ is the radius of the circular trajectory. Thus, angular velocity ($\Omega$) becomes 
\begin{align}
\Omega=\frac{p_1}{I}=\frac{M v r}{I}\,.
\end{align}
Noting that $v=\xi^2$, we can write the total energy of the skate as
\begin{align}
KE_{{\rm skate}}=&\frac{1}{2}M(\xi^2)^2+\frac{1}{2}I\left(\frac{M\xi^2r}{I}\right)^2
=\frac{1}{2}M(\xi^2)^2\left(1+\frac{M r}{I} \right)\,,
\end{align}
where $M$, $r$, and $I$ are known constants. 
Next, we derive the total energy of the control mass. The spatial coordinates of the control mass, (A,B), are given by 
\begin{align}
(A,B)=&(x,y)+\begin{bmatrix}
\cos(\theta) & -\sin(\theta)\\
\sin(\theta) & \cos(\theta)
\end{bmatrix}(a,b)\\
=& (x+a\cos\theta-b\sin\theta, y+a\sin\theta+b\cos\theta)\,.
\end{align}
The components of the velocity for the control mass in spatial coordinates can then be written as 
\begin{align}
v_A=&\dot x + \dot a \cos\theta-a \dot\theta \sin\theta-\dot b\sin\theta-b\dot\theta \cos\theta\\
v_B=&\dot y +\dot a \sin\theta+a \dot\theta\cos\theta+\dot b \cos\theta-b \dot\theta\sin\theta\,.
\end{align}
The kinetic energy of the control mass $m$ is computed as follows. 
	Suppose $\mathbf{e}_1$ is a vector in the body frame parallel to the skate, $\mathbf{e}_2$ normal to the skate and $\mathbf{e}_3=\mathbf{e}_1\times \mathbf{e}_2$ is normal to the plane. We remind the reader the $\xi^1$ has the physical meaning of the angular velocity of skate's rotation and $\xi^2$ has the physical meaning of the velocity of the skate. Then, the velocity of the control mass $m$ is given by
	\begin{equation} 
	\begin{aligned} 
	\mathbf{v}&=\xi^2 \mathbf{e}_1 + \xi^1 \mathbf{e}_3 \times (a \mathbf{e}_1 + b \mathbf{e}_2) + \dot a \mathbf{e}_1 + \dot b \mathbf{e}_2 
	\\
	&= \mathbf{e}_1 \left( \xi^2 - \xi^1 b+ \dot a \right) +\mathbf{e}_1 \left(  \xi^1 a+ \dot b \right).
	\end{aligned} 
	\label{vec_vm} 
	\end{equation} 
	Thus, the kinetic energy of the control mass is given by 
	\begin{equation} 
	\begin{aligned} 
	\hspace{-3mm}KE_m &= \frac{1}{2} m |\mathbf{v}|^2 = \frac{m}{2} \left[ \left( \xi^2 - \xi^1 b+ \dot a \right)^2 + \left(  \xi^1 a+ \dot b \right)^2 \right] 
	\\
	& = \frac{m}{2}\left[ \dot a ^2+\dot b^2 \!+\!(\xi^1)^2(a^2+b^2)\!+\! (\xi^2)^2\!+\!2\xi^1(a\dot b -\dot a b )+2\xi^2(\dot a-b \xi^1)\right]\,.
	\end{aligned} 
	\end{equation} 
As shown in Figure \ref{fig: Energy_NewControl}, the energy profiles for the optimized control mass and skate for each arc remain bounded as desired.

\begin{figure} 
	\centering
	\includegraphics[width=0.5 \textwidth]{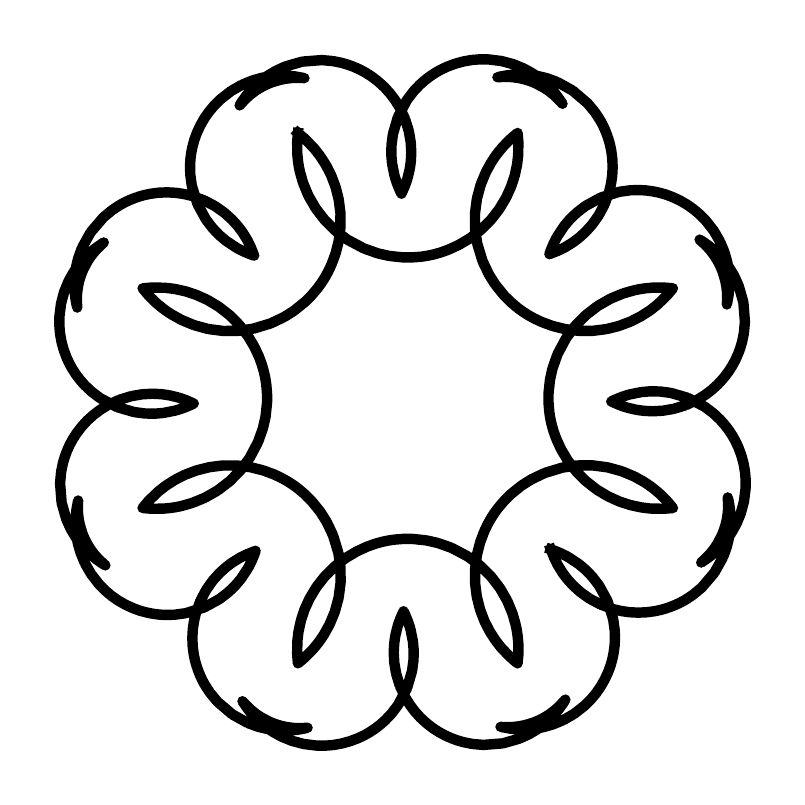} 
	\caption{Full pattern reconstruction using the control procedure. Compare to the target Figure~\ref{fig:OldPattern}.}
	\label{fig:SimulatedPattern_Length_NewControl}
\end{figure}

Finally, the full pattern is produced by transforming (via rotation and translation) and repeating the optimized arcs is shown in Figure \ref{fig:SimulatedPattern_Length_NewControl}. Notice the clear resemblance of the pattern reproduced by our method to the target ``double flower'' pattern of Figure~\ref{fig:OldPattern}. In principle, we could reproduce a large class of piecewise smooth curves, where the smooth parts of the curves are approximated by circular arcs. The limitation on control is then at the points of sharp turns, where the velocity of the skate needs to vanish. Thus, a control to successfully trace trajectories is possible whenever it is possible to create the motion using the control mass while enforcing vanishing velocity at two ends for every smooth part of the curve. We believe that this control mechanism spans a large number of possible trajectories and is substantially simpler than the ``brute force'' trajectory tracing without specifying particular constraints on the motion of the control masses. 

\section{Conclusion}

We have applied the control of the moving mass on a Chaplygin sleigh to model the dynamics of a figure skater tracing a prescribed pattern on the ice. With the added mass representing the controllable center of mass of the figure skater, we determined the controls needed to reproduce an example pattern, which is approximated by circular arcs. The approximation of circular arcs was chosen for the simplicity of algorithms, as in that case the equations for the control simplify considerably: the question of trajectory tracing, which is very complex, simplifies to the matching the beginning and the end of the arc, with the curvature being matched automatically. 

For a more complex trajectory tracing, a machine learning algorithm could be used. For example, deep reinforcement learning algorithms have been used with success for drone control in 3D \cite{koch2019flight}. The problem of tracing skating trajectories is somewhat different as it is done in 2D with a nonholonomic constraint. If a more complex model of human skater than a Chaplygin sleigh is used, this model will also involve substantially more degrees of freedom than a drone. Nevertheless, the application of reinforcement learning to the problem of trajectory tracing looks promising for further research and we hope that the researchers in the field of reinforcement learning will become interested in tackling this challenging problem. 

\section*{Conflict of interest}
The authors declare that they have no conflict of interest.

\section*{Acknowledgements}
MR acknowledges the Pacific Institute for the Mathematical Sciences (PIMS) MathBio Accelerator Award for partial support during this project. The authors are grateful for the infinite patience of, and many discussions with, Prof. D. V. Zenkov. We are grateful to Profs. A. Bloch and D. D. Holm for their constant interest and thoughtful insights in the mechanics and control of nonholonomic systems. We are indebted to V. Gzenda for many discussions related to the mechanics of figure skating robots. We are grateful to Dr. V. Fedonyuk for pointing out some recent references in the control and dynamics of Chaplygin's sleigh. We are appreciative of the work by SIAM News staff, and in particular Lina Sorg, for their interest in showcasing this work in \cite{RhPu2021}.

\bibliographystyle{spmpsci}      
\bibliography{FigureSkatingBib}   

\appendix
\label{app_sec:alt_control} 
\section{Alternative methods of trajectory tracing} 

\subsection{Simulations with Generalized Control Functions}
Here, we explore a more general control function in place of the controls given by \eqref{eqn: controls}. In addition to the length-based cost function \eqref{eqn:CLengthandSpeed}, we will employ a point-based cost function to examine the applicability of optimizing trajectories based on a desired end position.   

Restriction of the control to circular trajectories through $\xi^2=r\xi^1$ also introduces a constraint on the functions $a(t)$ and $b(t)$ and their derivatives. Thus, we can define one control, say $(a(t),\dot a(t))$, and solve a differential equation to determine $b (t)$. To allow for a wide range of control functions while maintaining a simple form, the control function $a(t)$ is taken to be 
\begin{align}
a(t)=&a_1+a_2\sin  \omega t +a_3 \cos  \omega t  \quad \Rightarrow \quad \dot a = \omega a_2 \cos \omega t - \omega a_3 \sin \omega t
\label{eqn: general control function}
\end{align}
where $a_i, i=1,2,3$, and $\omega$ are parameters to be optimized to generate a prescribed curve. The equation  for $\dot b(t)$ is then derived from \eqref{arcs_xi1_xi2} using \eqref{xi1_eq} and \eqref{xi2_eq}, and solving for $\dot b(t)$, leading to: 
\begin{align}
\dot b =&\frac{mbp_1-m\dot a (I+ma^2)+p_2[I+m(a^2+b^2)]-r(M+m)p_1-rmb(p_2+M\dot a)}{m^2ab-r(M+m)ma}\,,
\label{b_diff_eq}
\end{align}
where $ a $ and $ \dot a $ are known. This differential equation has to be solved together with the dynamic equations for $p_1$ and $p_2$. Unfortunately, as one can easily see from \eqref{b_diff_eq}, the equation for $\dot b$ is singular whenever $a=0$, so the velocity of the control mass, and hence the energy necessary for control, diverges at $a=0$. The condition for solutions for $b(t)$ to remain regular is 
\begin{equation} 
a_1^2 > a_2^2 +a_3^2 \, .  
\label{sing_cond}
\end{equation} 
When the singularity condition \eqref{sing_cond} is violated, the velocity of the controls may diverge at a finite set of time points. This singularity presents some problems in the physical interpretation of results in the case when \eqref{sing_cond} is violated, as it leads to an infinite control force. This is in contrast with the method presented in our main paper, where the control forces always remain finite. In the following text, we present two cases, one when \eqref{sing_cond} is satisfied (point-based optimization, Sec.~\ref{sect: PointBased}) and one when \eqref{sing_cond} is violated (length-based optimization, Sec.~\ref{sect: LengthBased}). 

Under the conditions that trajectories are circular arcs and the controls are coupled through \eqref{b_diff_eq}, the Chaplygin sleigh model with the control \eqref{eqn: general control function} is  
\begin{equation} 
\begin{aligned}
\dot{p_1}(t) &= -m\eta \xi^2\,,\\
\dot{\theta}(t) & = \xi^1\,,\\
\dot{x}(t) & =\xi^2\cos \theta\,,\\
\dot{y}(t) &= \xi^2\sin \theta\,,\\
\dot b(t)   & \mbox{ given by \eqref{b_diff_eq}}\,.
\end{aligned}
\label{all_dynamics_eqs} 
\end{equation} 
with initial conditions 
\begin{equation} 
\begin{aligned}
p_1(0) =& \bar p_{1}\,,\\
\theta(0)=& \bar \theta\,,\\
x(0) =&\bar x\,,\\
y(0) =&\bar  y\,,\\
b(0) =&\bar  b\,,
\end{aligned}
\end{equation} 
and the constraints given by \eqref{eqn: circle constraints}.

With the control function $a(t)$ and $\dot a(t)$ given by \eqref{eqn: general control function}, we consider two cost functions to optimize. The first cost function is a point-based optimization tracking the trajectory, while enforcing vanishing speed at the end of each smooth trajectory. The second cost function considered is the length-based cost function with vanishing speed at the trajectory end points given by \eqref{eqn:CLengthandSpeed}. The numerical methods follow those outlined in Section \ref{sect: Numerics}, where the cost functions involved in step 4 are replaced with the appropriate cost functions (\eqref{eqn:Cpoint} for point-based and \eqref{eqn:CLengthandSpeed} for length-based). Under point-based optimization, a further step is required to produce the full trajectory which is described in the following section.

\subsection{Point-based Optimization}
\label{sect: PointBased}
In this section, the circular arc is chosen so it reaches a given end point in space, given vanishing velocities at the end points. Denoting the target endpoint of the forward-in-time trajectory by $(X_{opt}, Y_{opt})$, the cost function is 
\begin{align} \label{eqn:Cpoint}
C_{point}=& \left[(X(t_2)-X_{opt})^2+(Y(t_2)-Y_{opt})^2\right]^p\,, p>1\,, \text{ with } \xi^2(t_2)=0\,,
\end{align}
where $t_2$ is the end time of the trajectory. 

The optimization function solves only the forward-in-time trajectory to obtain the optimal control parameters, with the backward-in-time trajectory not involved in optimization. The optimized control parameters are then used to generate the backwards-in-time trajectory. In the numerical methods, the generation of the backwards-in-time trajectory is an additional step that would occur in step 4 of the numerical methods outline in Section \ref{sect: Numerics}.

The simulation parameters for each arc are given in Table \ref{tab:NumericsPars}. 
\begin{table} 
	\caption{Values used in simulations to create arc 1, arc 2, and arc 3 with point-based optimization. In all cases $m=1$, $M=2$, and $I=3$.}
	\label{tab:NumericsPars}	
	\centering
	\begin{tabular}{llllll}
	\hline\noalign{\smallskip}
		Arc  & $T$ & $r$ & $(X_{opt}, Y_{opt})$ &  $(\bar{p_1},\bar{p_2}, \bar{\theta},\bar{x},\bar{y},\bar{b})$ & $(a_1, a_2, a_3, \omega)$ \\ 
		\noalign{\smallskip}\hline\noalign{\smallskip}
		1 & $8$ & $1.2$& $(\sqrt{27}/5, 3/5)$ & $(1.5, 2.5, 0, 0, -1.2, 1.1)$ & $(1, 1, 1, 1)$\\
		2 & $1.5$ & $1.0$&  $(1/5,-\sqrt{24}/5)$ & $(0.015, 0.025,  0, 0, -1, 0.015)$ & $(0.1, 0.1, 0.1, 2)$\\
		3 & $8$& $1.0$ &  $(1,0)$ & $(1.5, 2.5, 0, 0, -1, 0.75)$ & $(1, 1, 1, 1)$  \\ 
		\noalign{\smallskip}\hline
	\end{tabular}
\end{table}

The initial guess for control parameters and resulting optimized parameters used in optimizing the Chaplygin sleigh solution for each arc are shown in Table \ref{tab:OptimizedPars}.
\begin{table} 
	\centering
	\caption{Initial guesses and optimized values for control parameters.}
	\label{tab:OptimizedPars}
	\begin{tabular}{lll}
	\hline\noalign{\smallskip}
		Arc  & Guess, $(a_1, a_2, a_3, \omega)$& Optimized, $(a_1^*, a_2^*, a_3^*, \omega^*)$ \\ 
		\noalign{\smallskip}\hline\noalign{\smallskip}
		1 & $(1, 1, 1, 1)$ & $(1.005, 1.013, 1.004, 1.018)$\\
		2 & $(0.1, 0.1, 0.1, 2)$& $(0.139, -0.509,  0.323, 2.011)$\\
		3 &  $(1, 1, 1, 1)$ & $(1.028, 1.082, 0.954, 1.099)$ \\ 
		\noalign{\smallskip}\hline
	\end{tabular}
\end{table}

\begin{figure} 
	\includegraphics[scale=0.39, trim={0.25cm 0 0 0 }, clip]{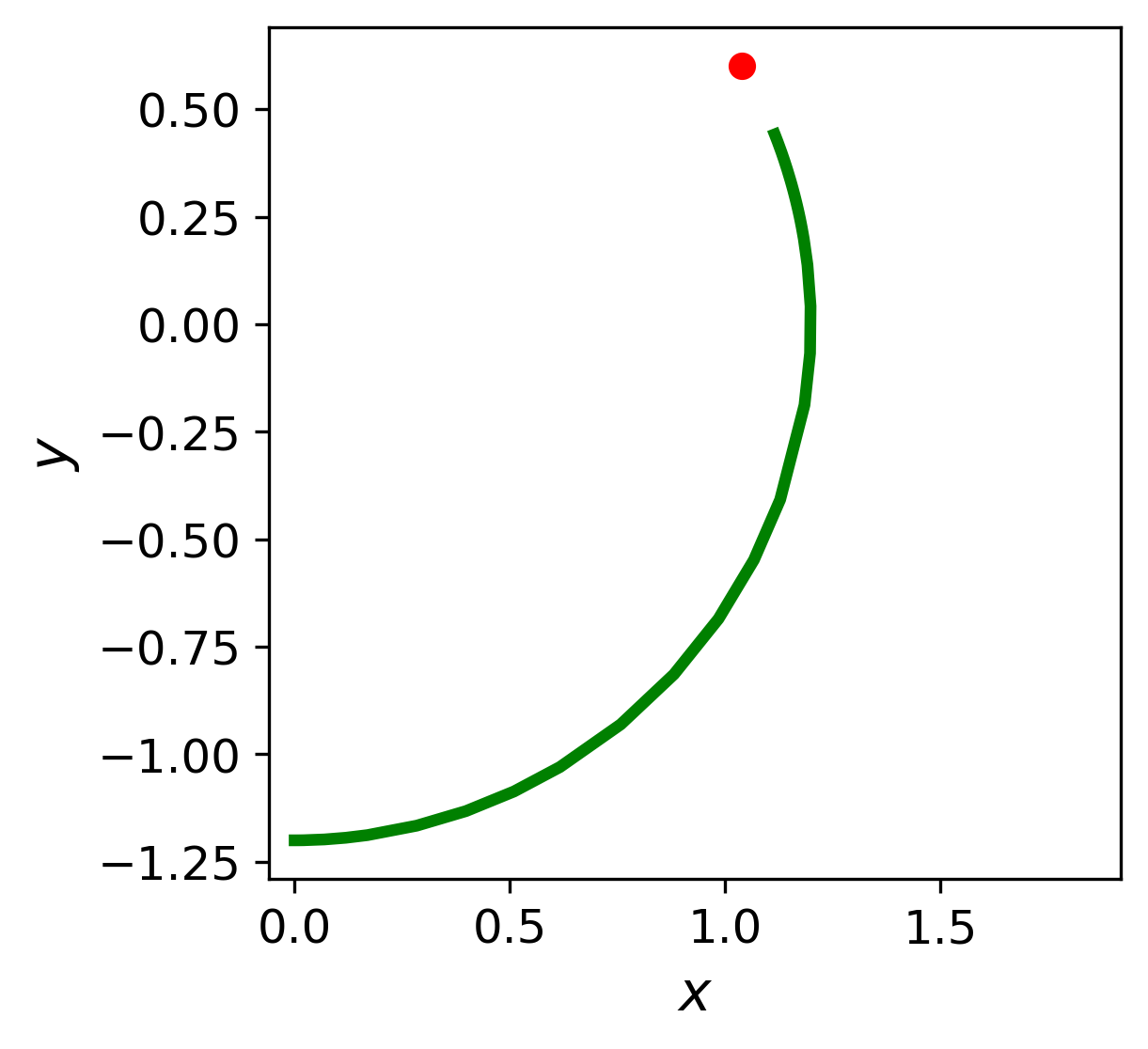}\hspace*{-0.1cm} \includegraphics[scale=0.39, trim={0.25cm 0 0 0 }, clip]{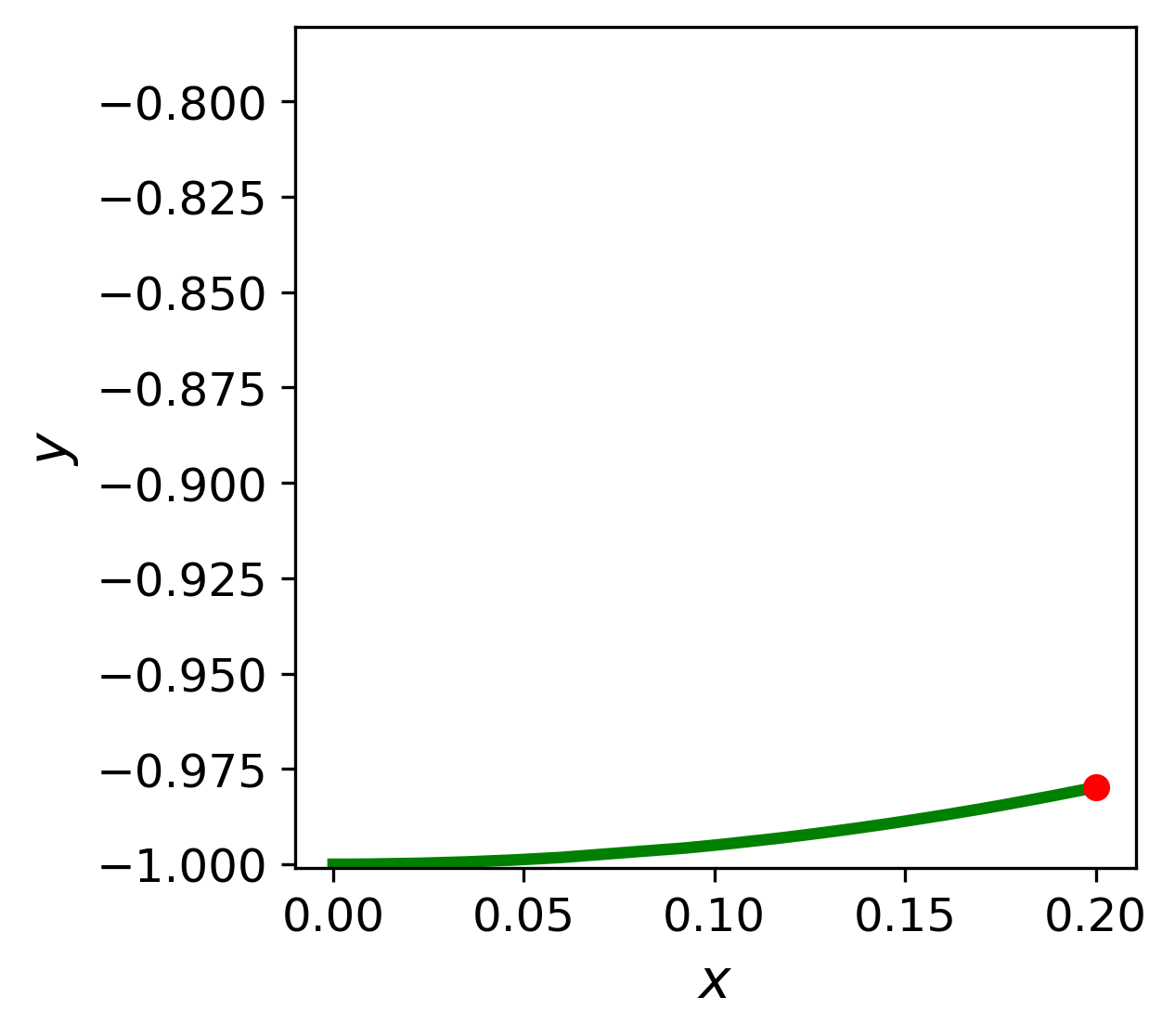}\hspace*{-0.1cm}	\includegraphics[scale=0.39, trim={0.25cm 0 0 0 }, clip]{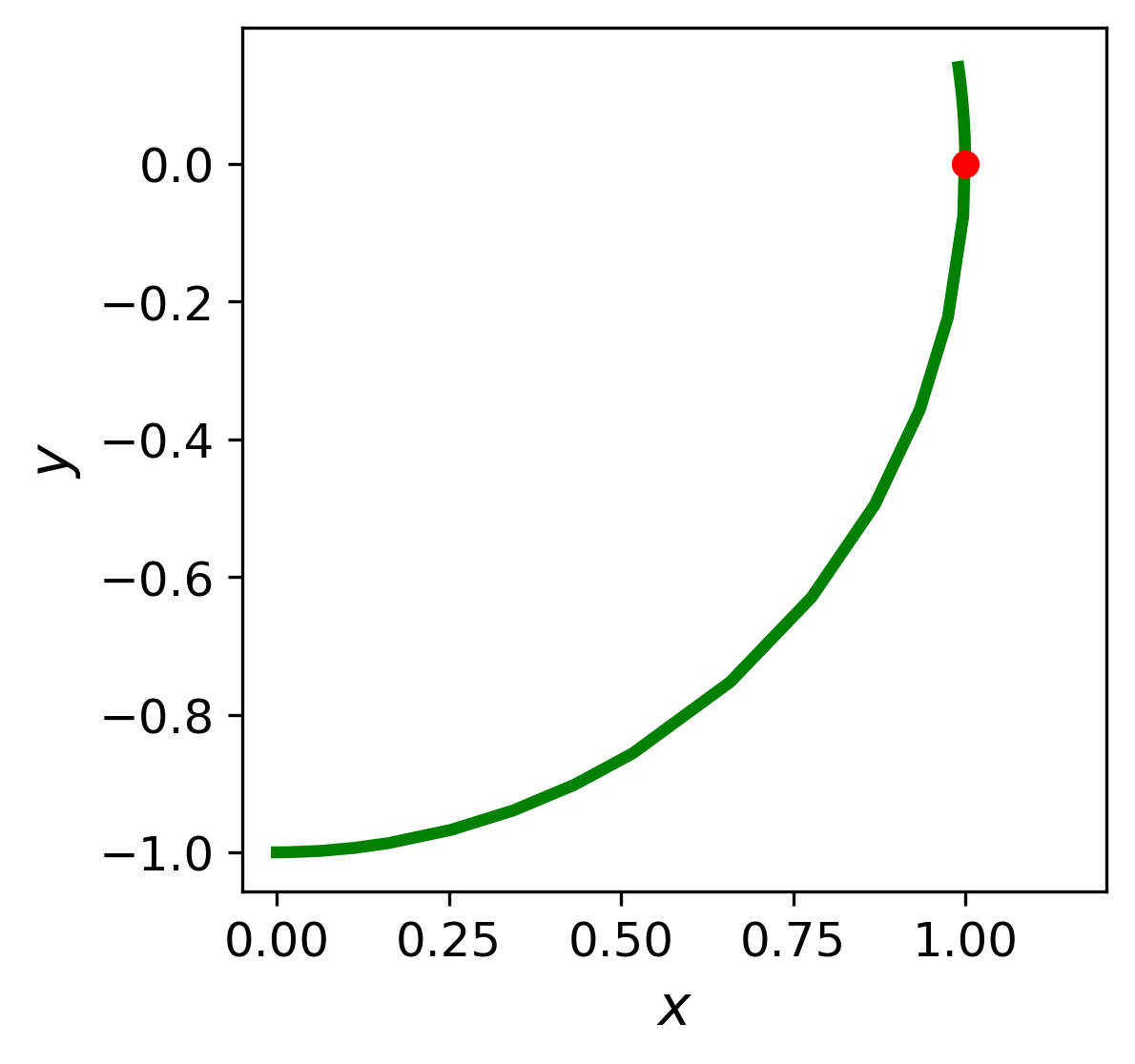}
	\caption{Arc trajectories and optimal end points (red marker) for arc 1 (top), arc 2 (bottom left), and arc 3 (bottom right). Green and blue indicate the forward- and backward-in-time solutions, respectively.}
	\label{fig:ArcEndPoints}
\end{figure}

\begin{figure} 
	\includegraphics[scale=0.4]{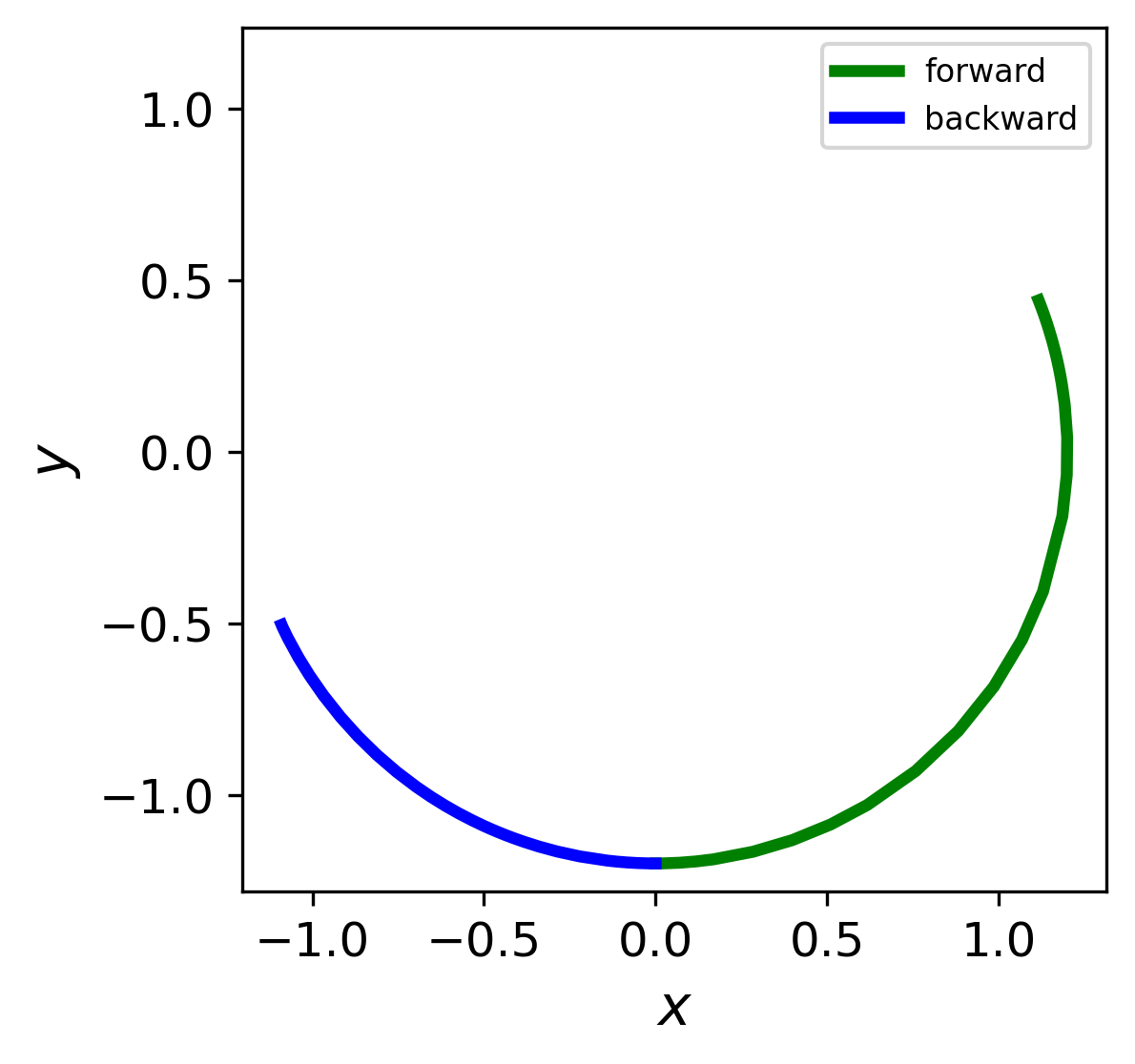} \hfill	\includegraphics[scale=0.4]{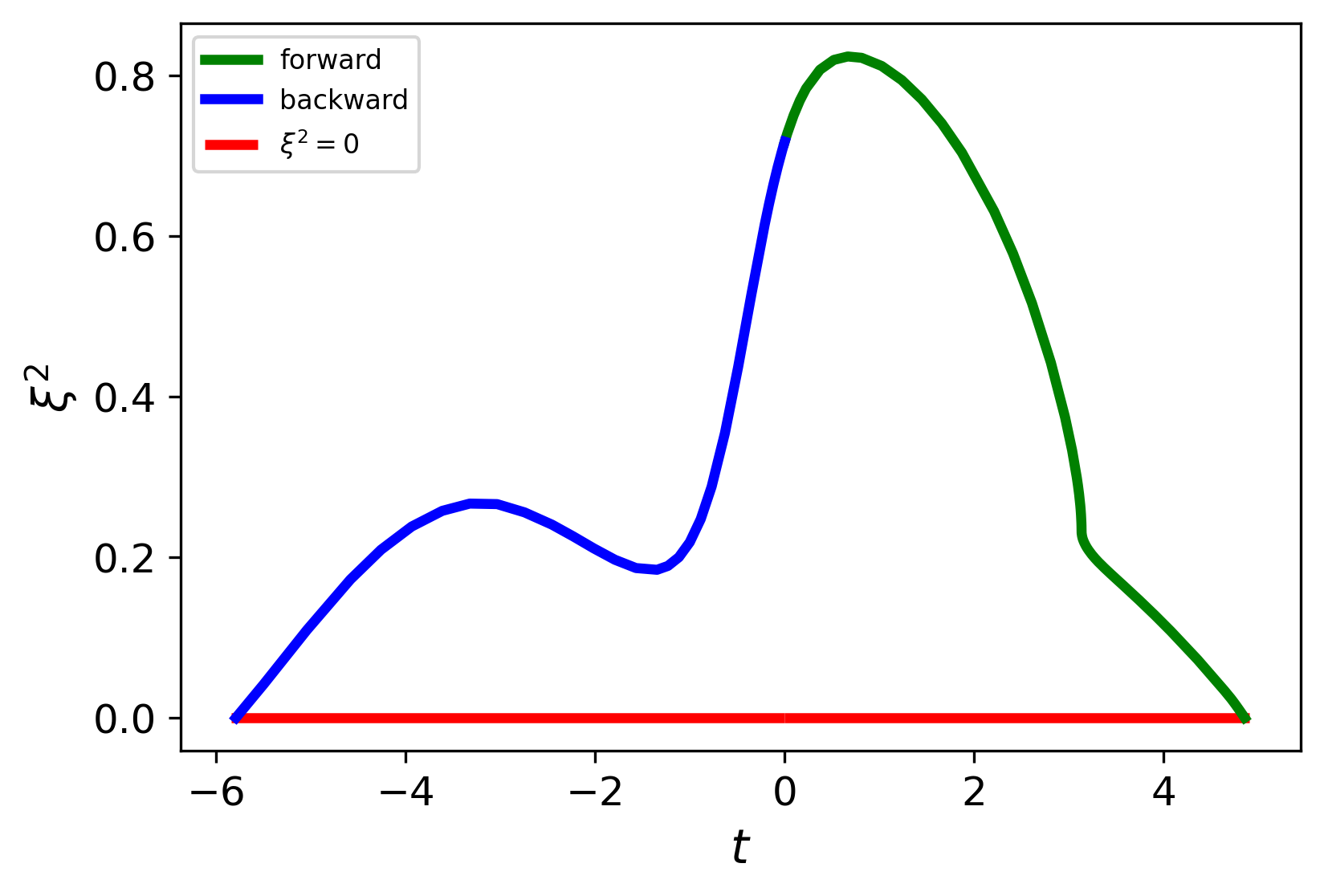}\hfill
	\includegraphics[scale=0.4]{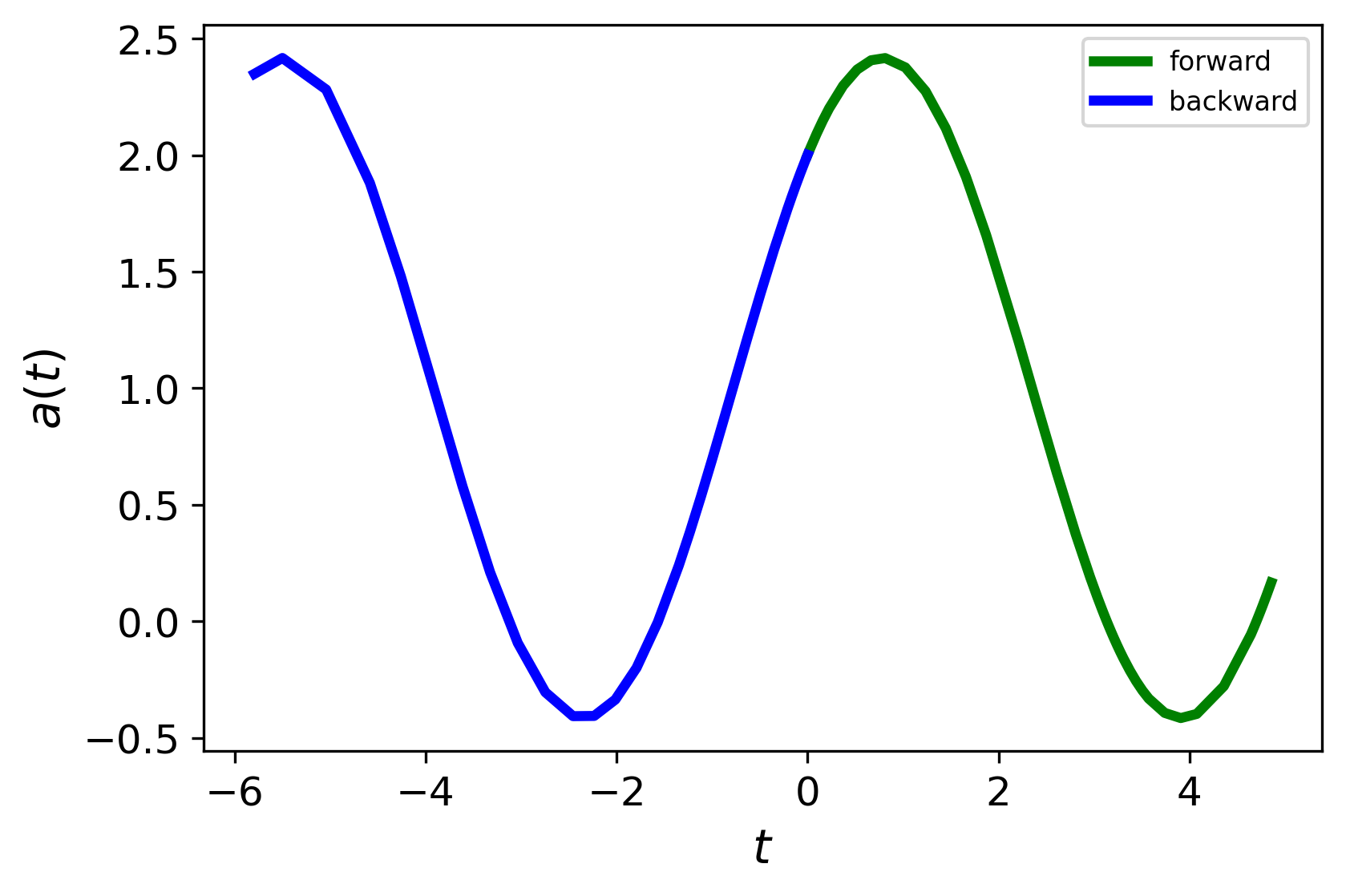}\hfill
	\includegraphics[scale=0.4]{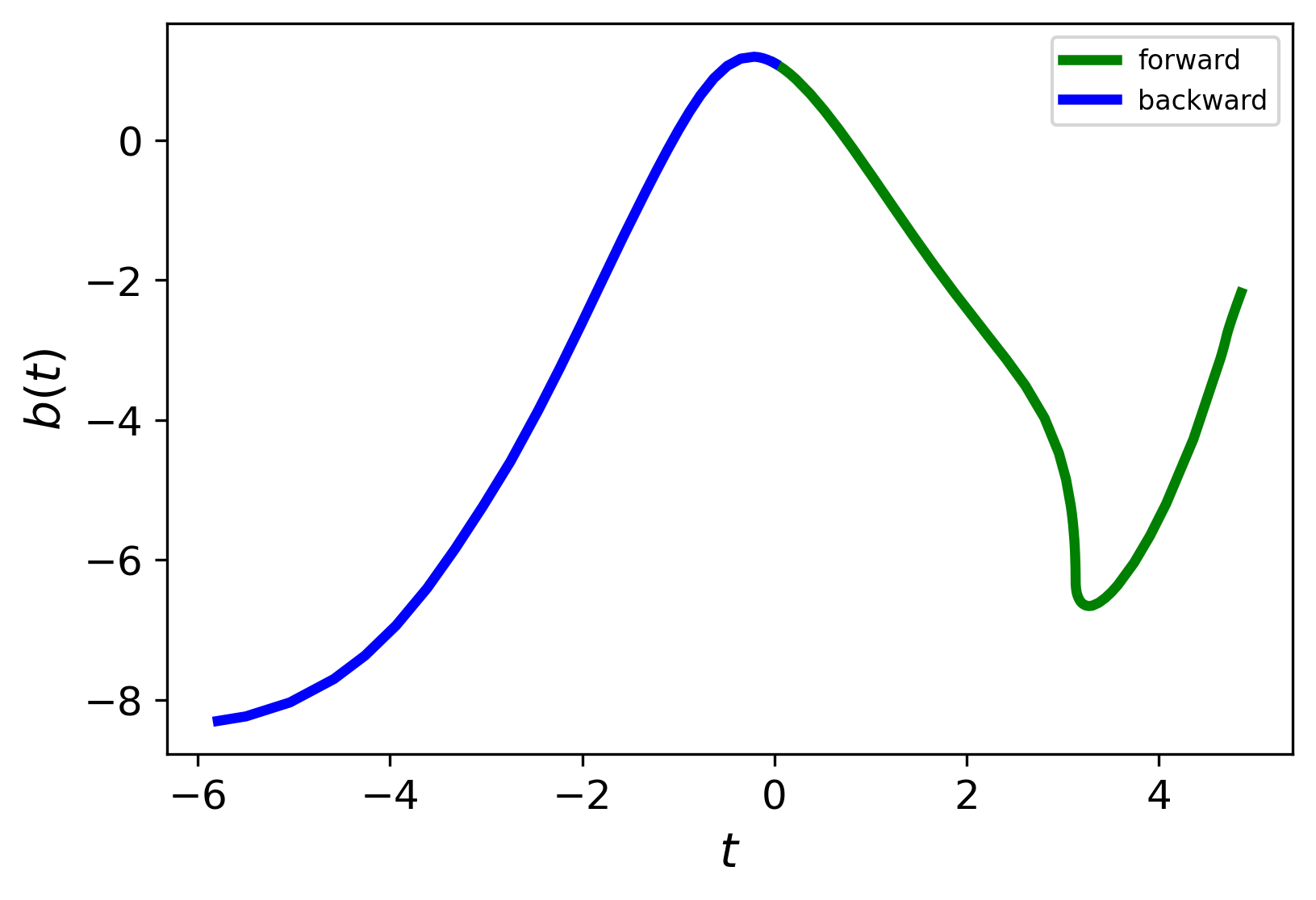}
	\caption{Arc 1 trajectory (top left), $\xi^2$ profile (top right), and optimized control functions, $a(t)$ (bottom left) and $b(t)$ (bottom right). Green and blue indicate the forward- and backward-in-time solutions, respectively.}
	\label{fig:Arc1}
\end{figure}

\begin{figure} 
    \includegraphics[scale=0.4]{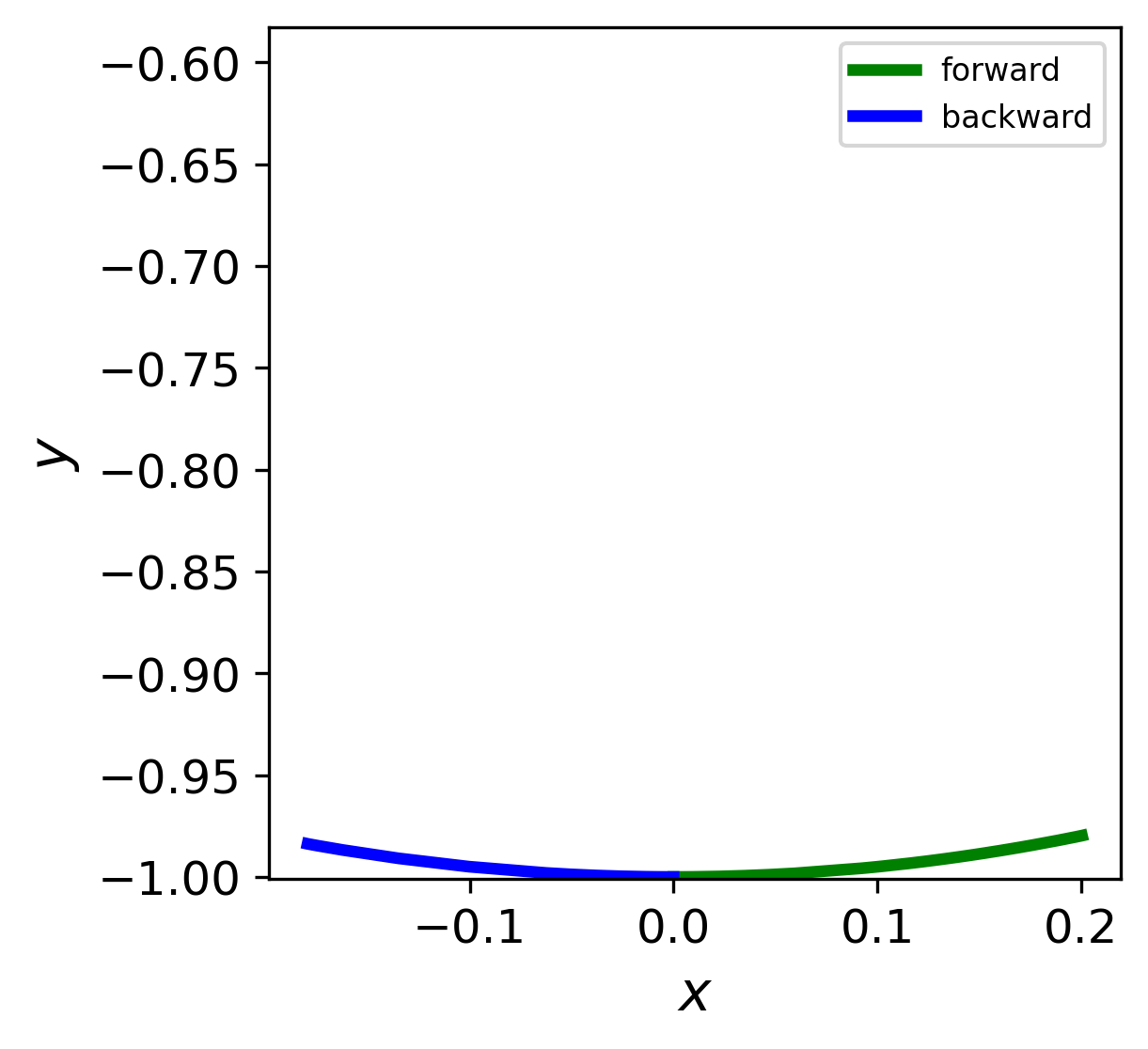}\hfill
    \includegraphics[scale=0.4]{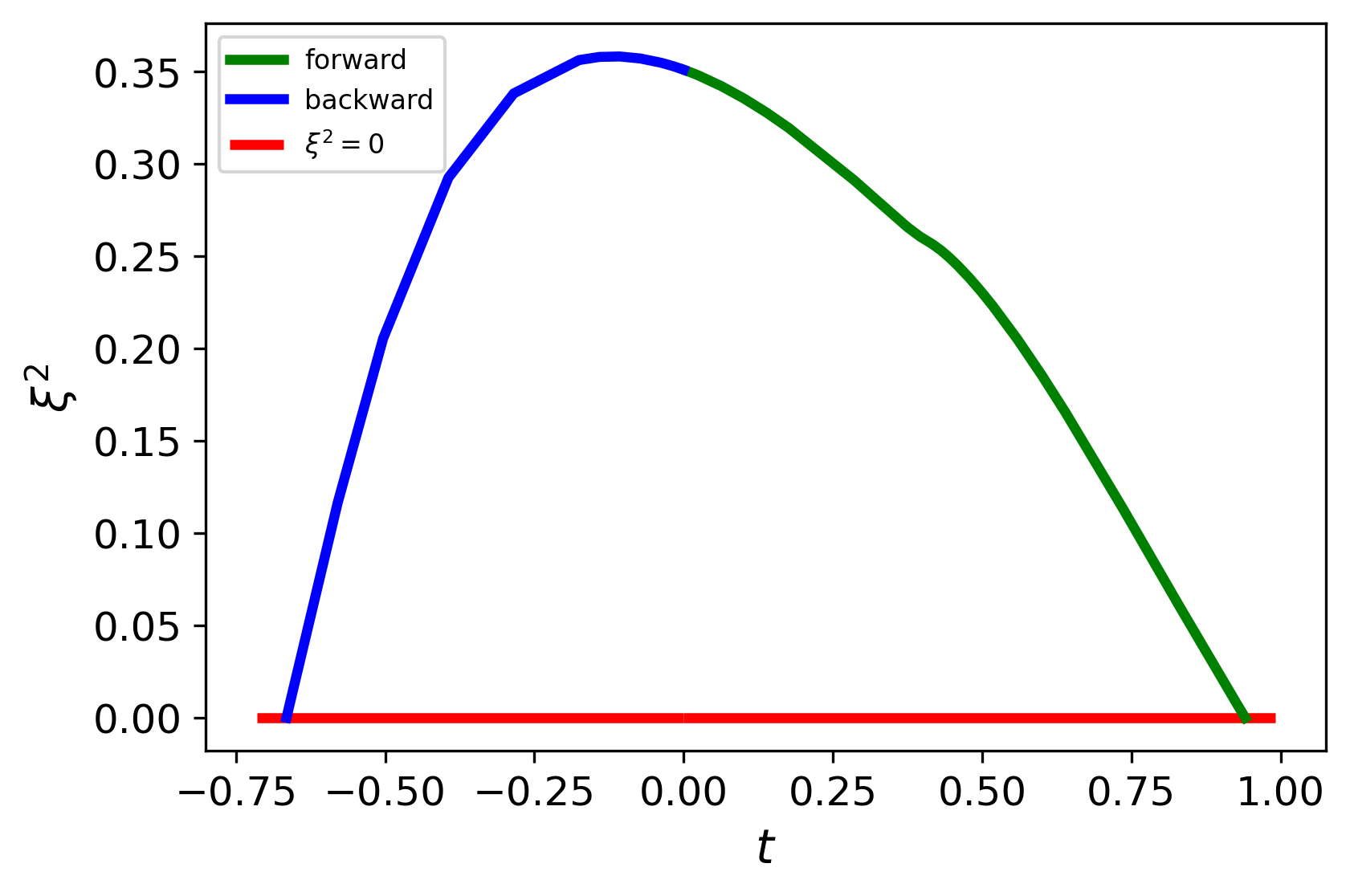}\hfill
    \includegraphics[scale=0.4]{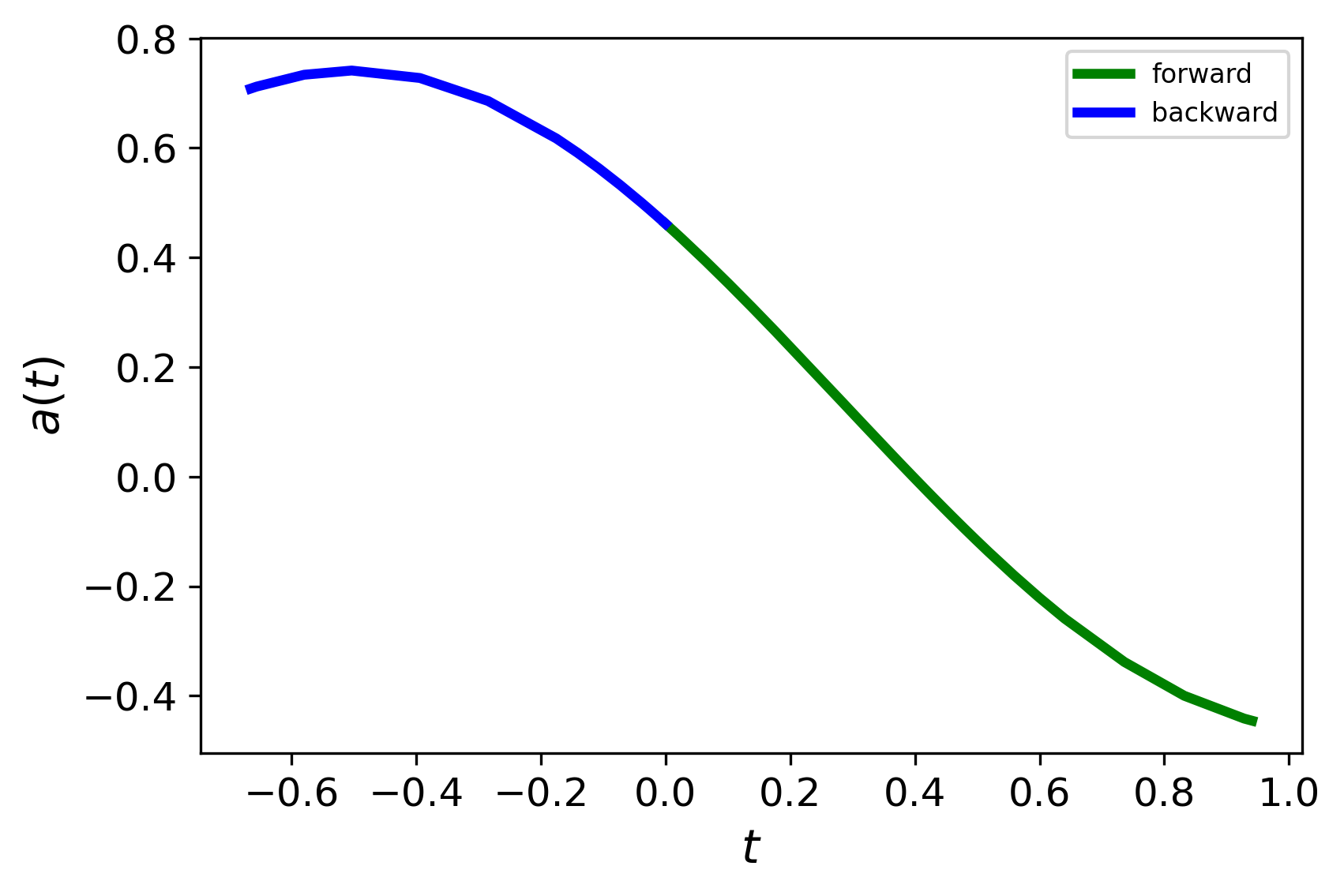}\hfill	\includegraphics[scale=0.4]{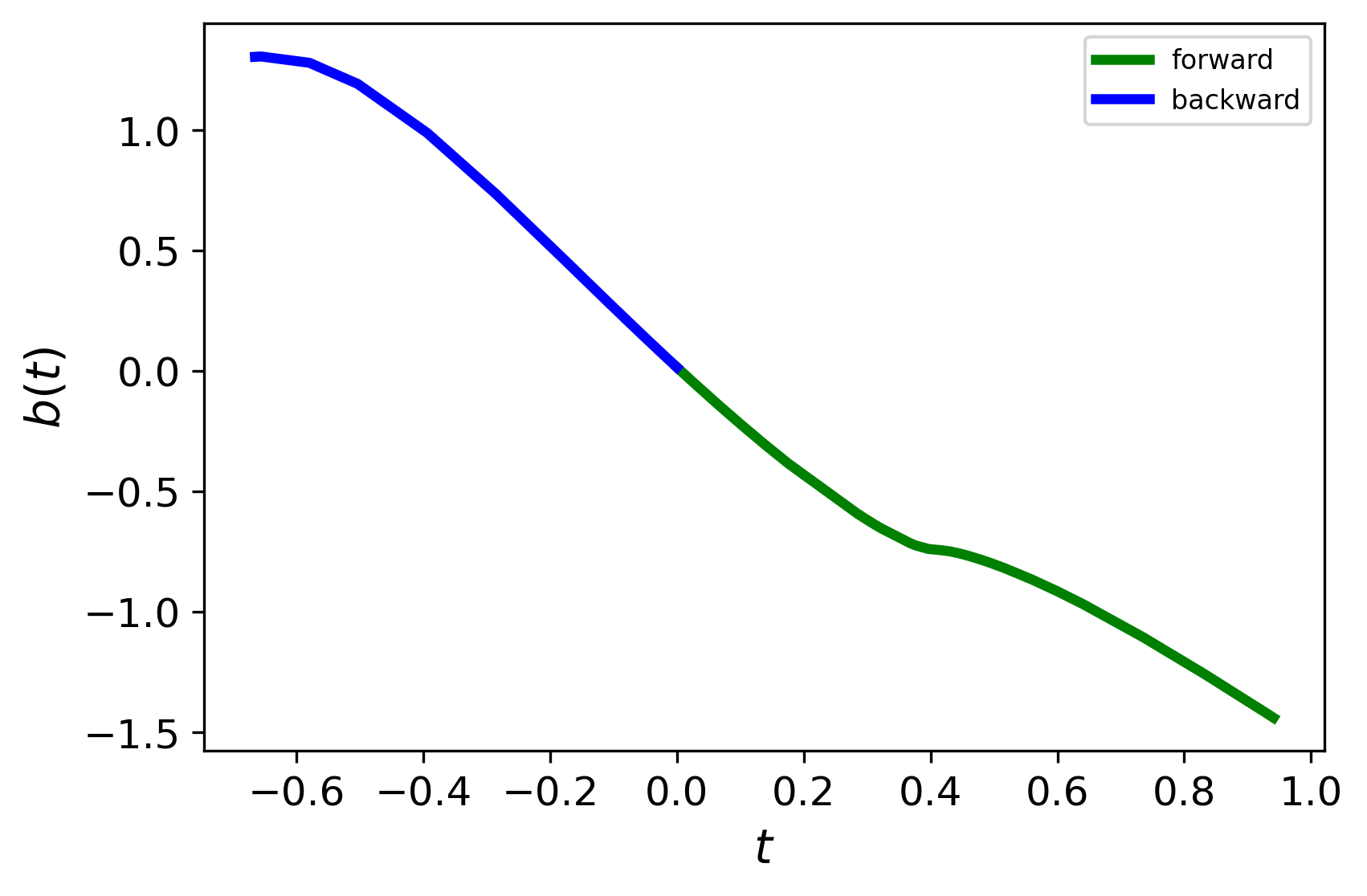}
	\caption{Arc 2 trajectory (top left), $\xi^2$ profile (top right), and optimized control functions, $a(t)$ (bottom left) and $b(t)$ (bottom right). Green and blue indicate the forward- and backward-in-time solutions, respectively.}
	\label{fig:Arc2}
\end{figure}

\begin{figure} 
	\includegraphics[scale=0.4]{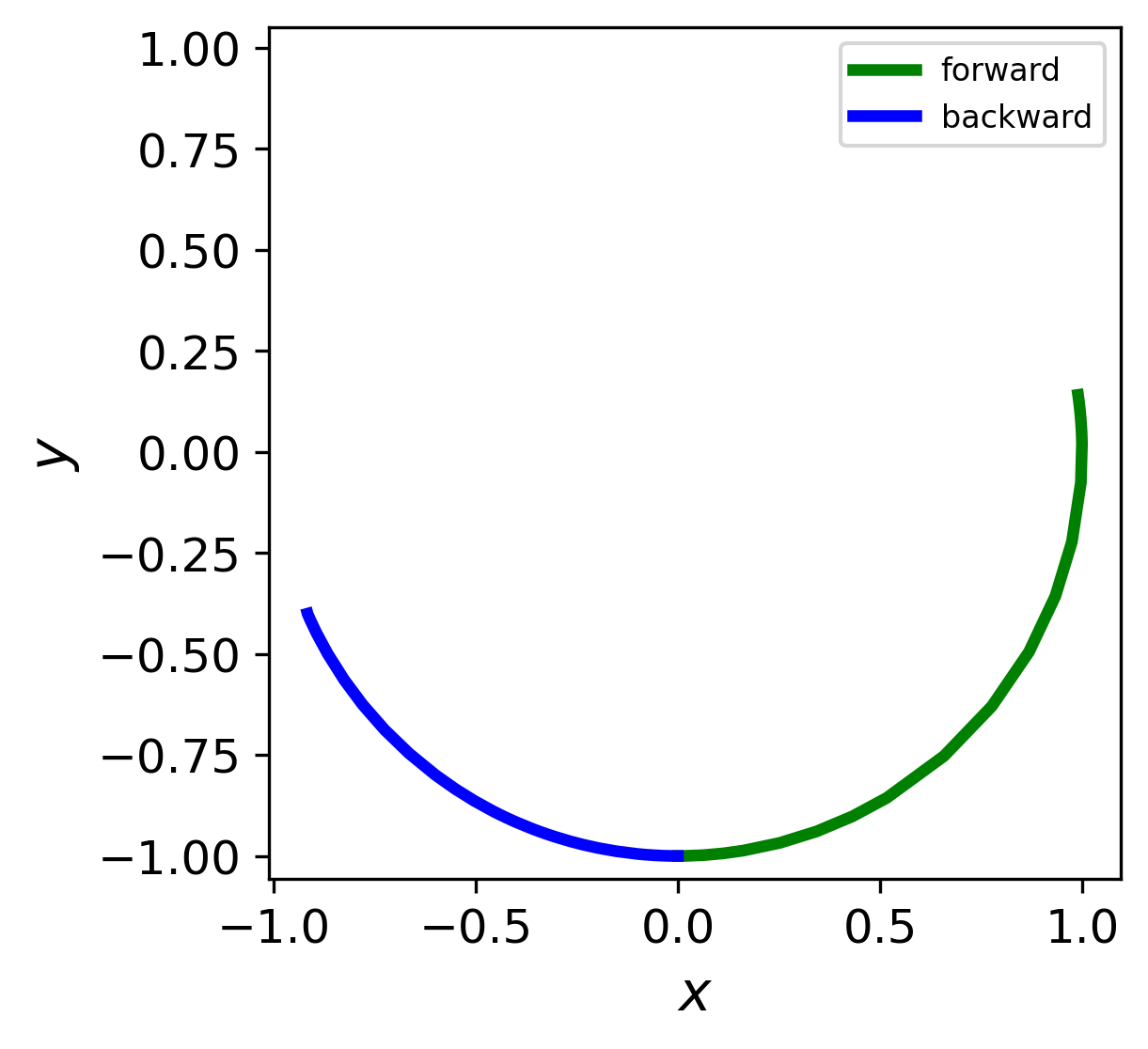}\hfill	\includegraphics[scale=0.4]{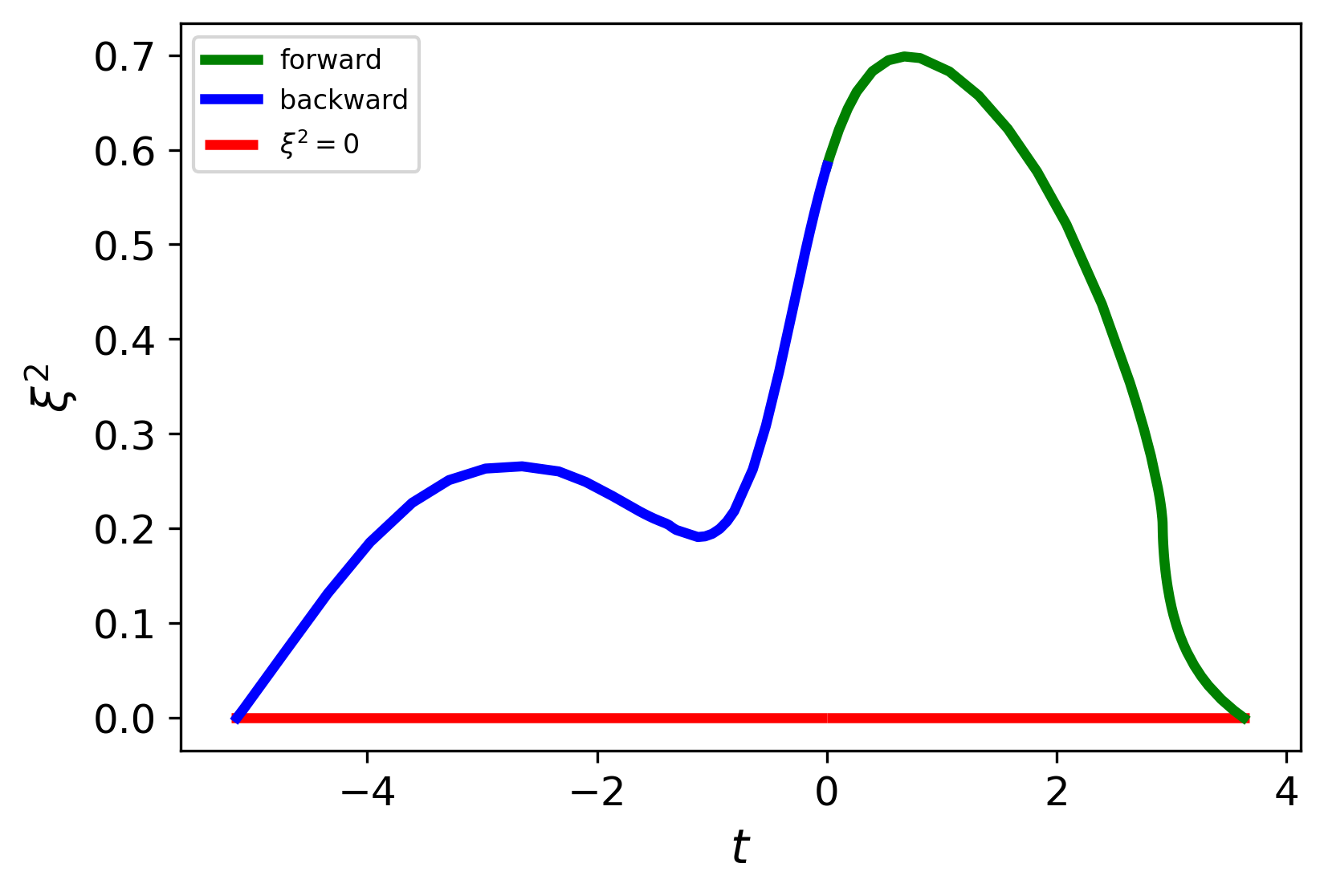}\hfill	\includegraphics[scale=0.4]{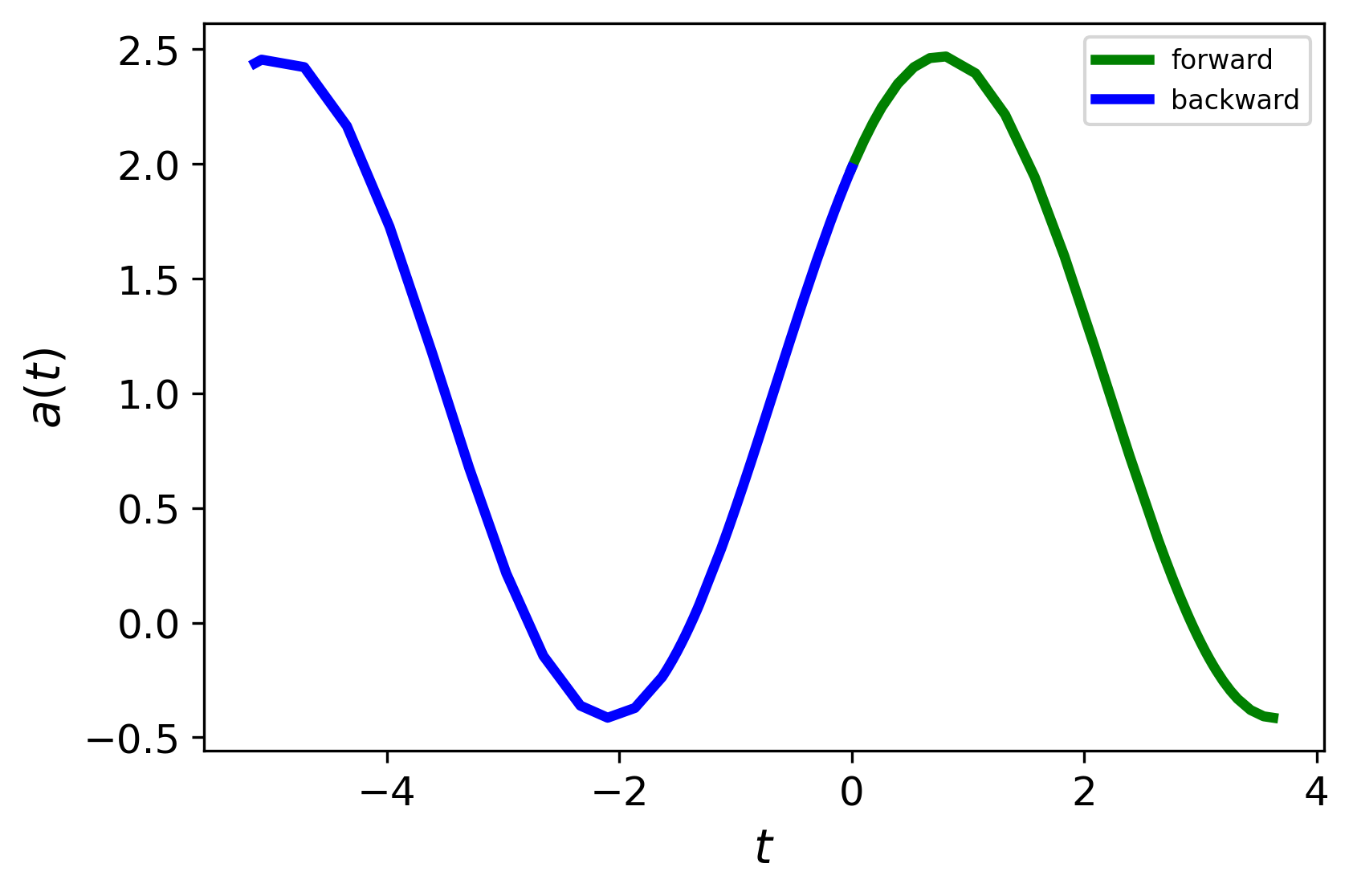}\hfill
	\includegraphics[scale=0.4]{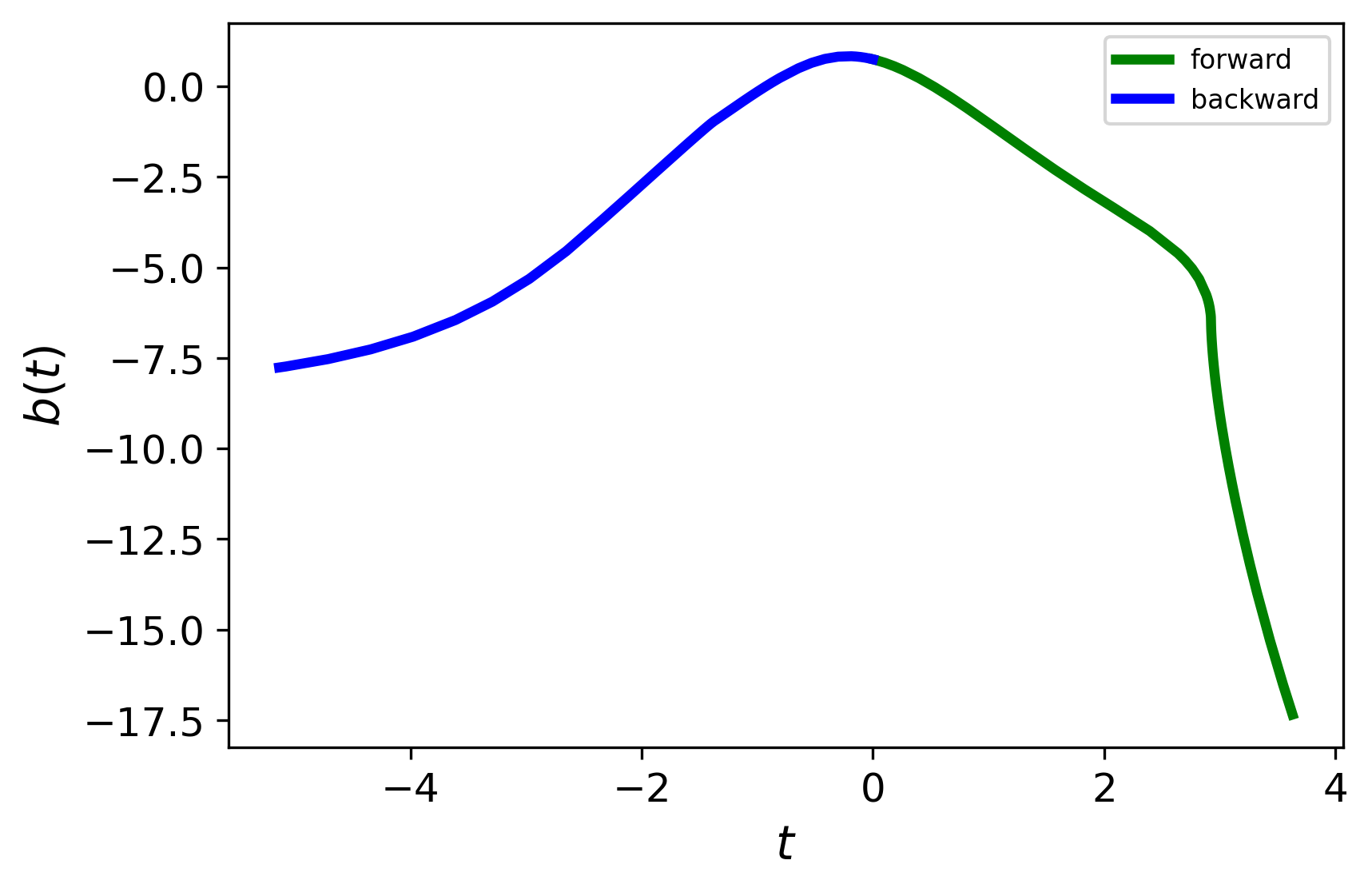}
	\caption{Arc 3 trajectory (top left), $\xi^2$ profile (top right), and optimized control functions, $a(t)$ (bottom left) and $b(t)$ (bottom right). Green and blue indicate the forward- and backward-in-time solutions, respectively.}
	\label{fig:Arc3}
\end{figure}

\begin{figure} 
\centering
	\includegraphics[scale=0.4]{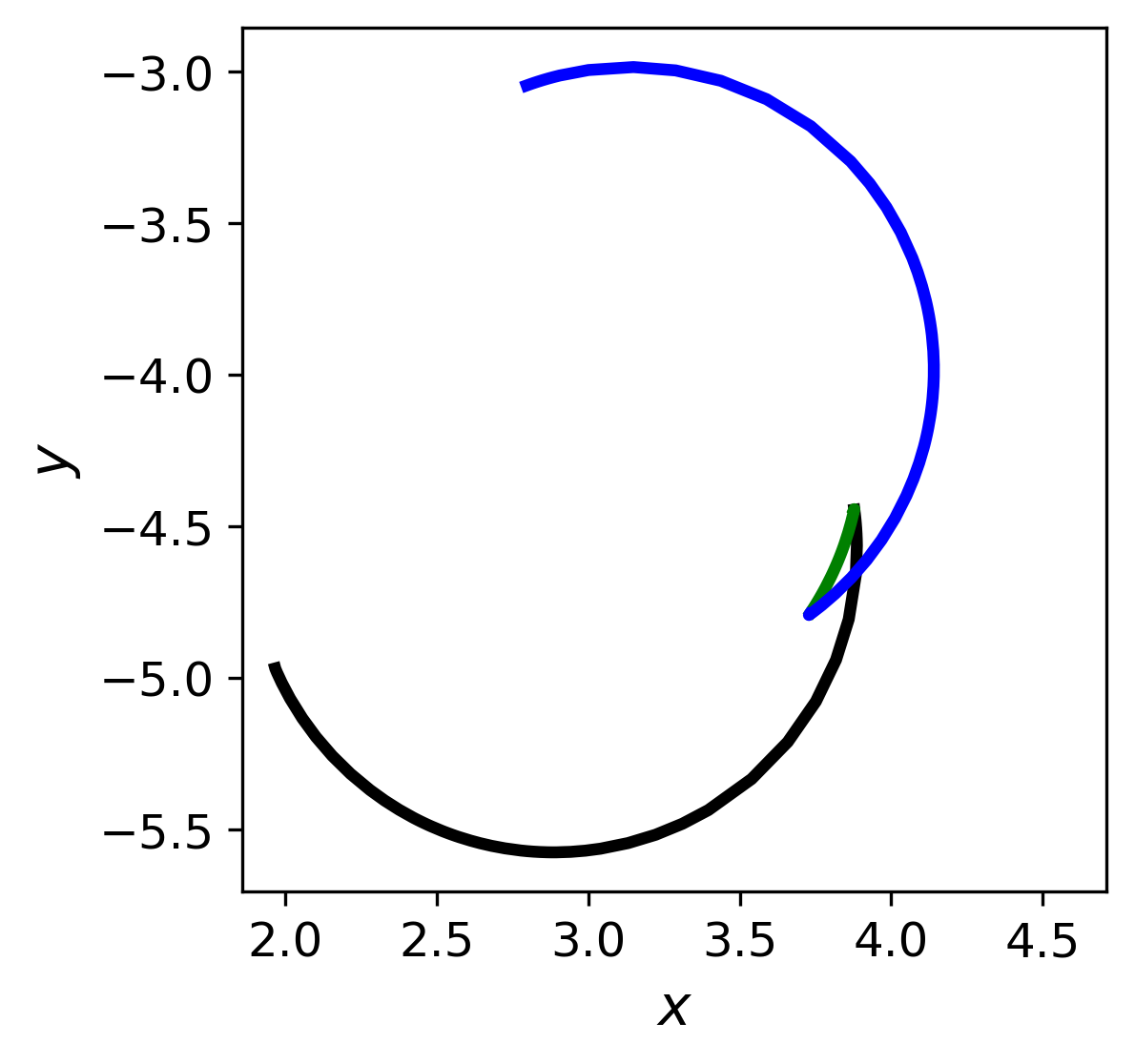}\hfill\includegraphics[scale=0.4]{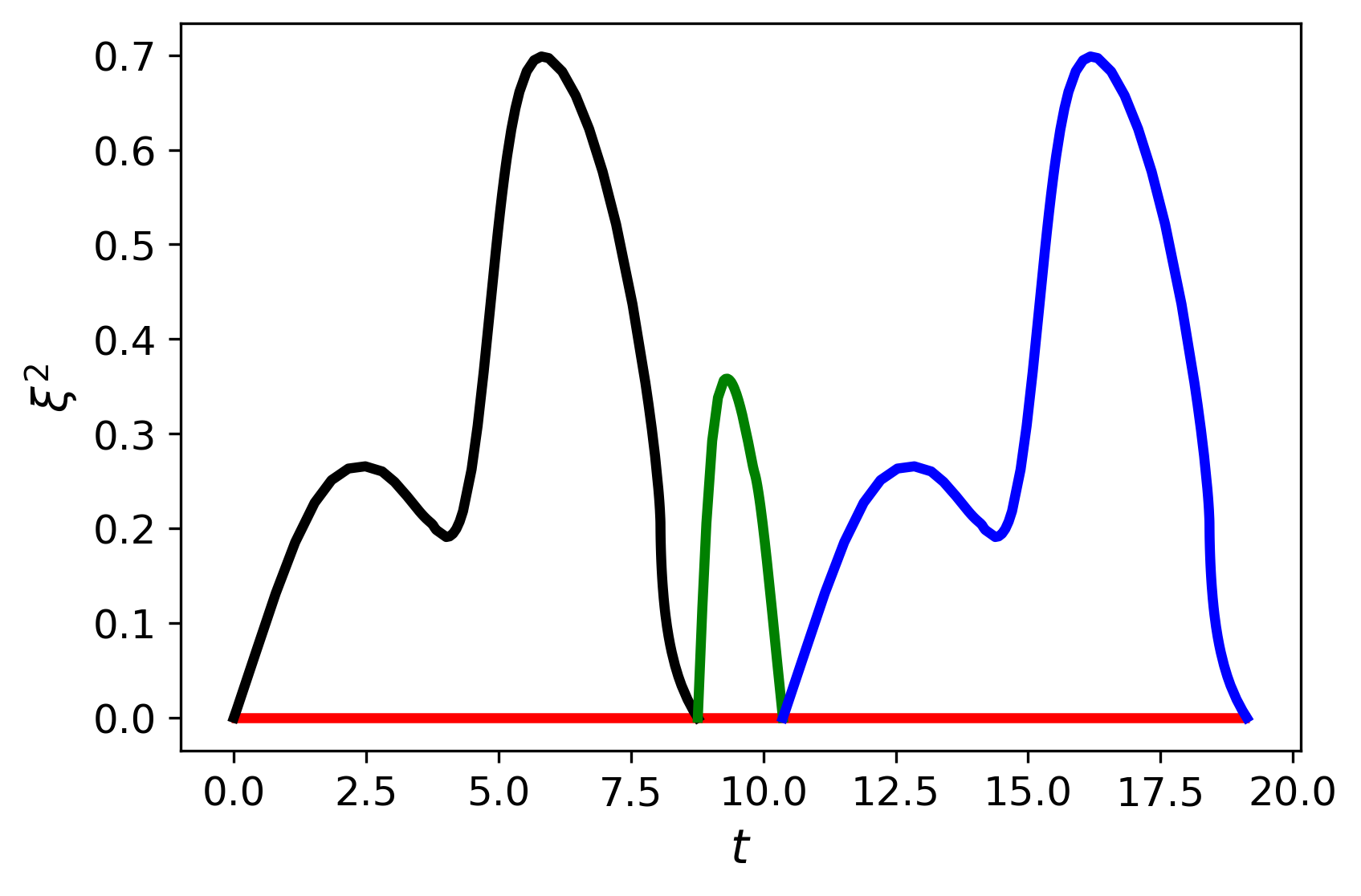}\hfill\includegraphics[scale=0.4]{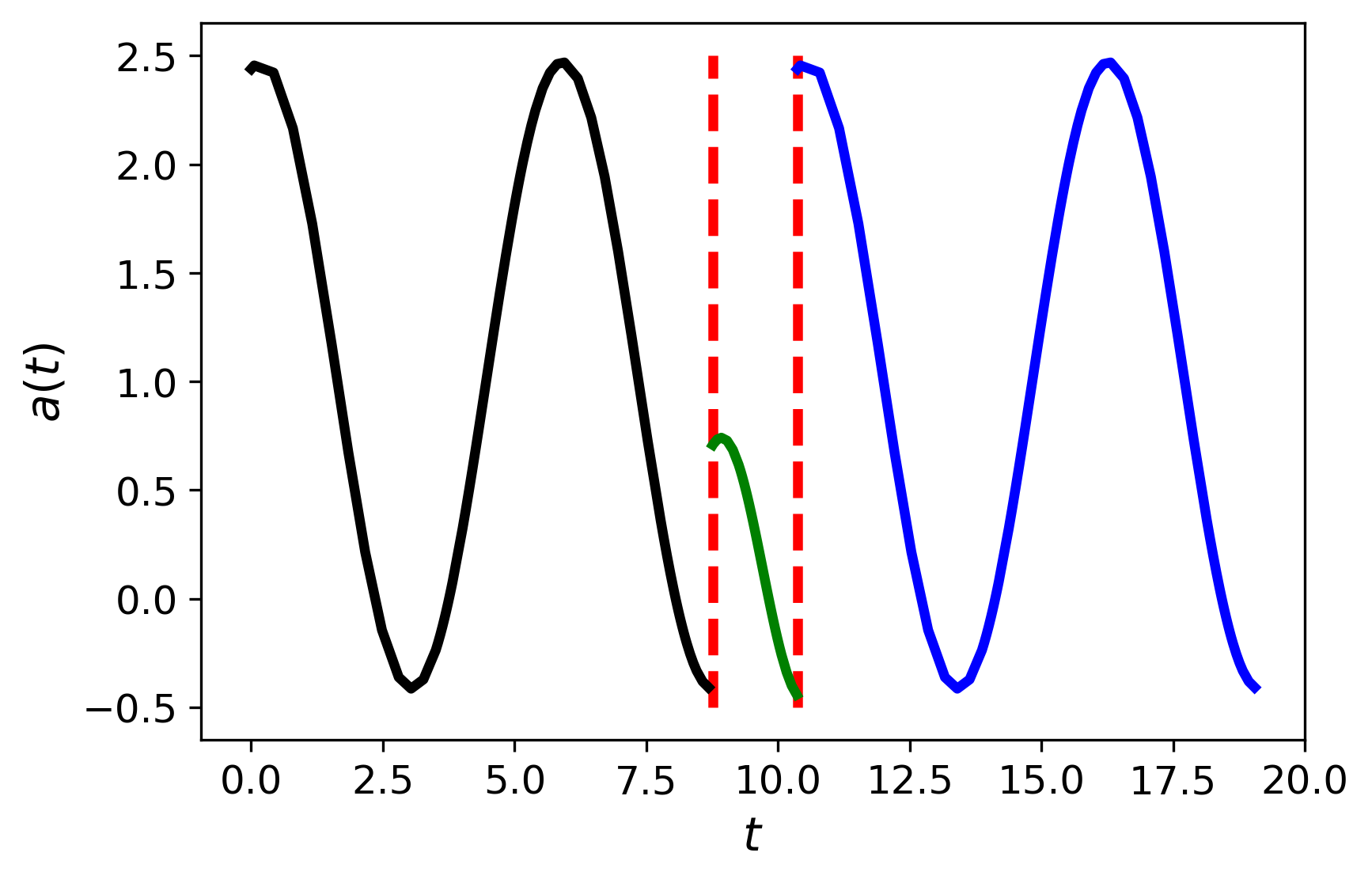}\hfill\includegraphics[scale=0.4]{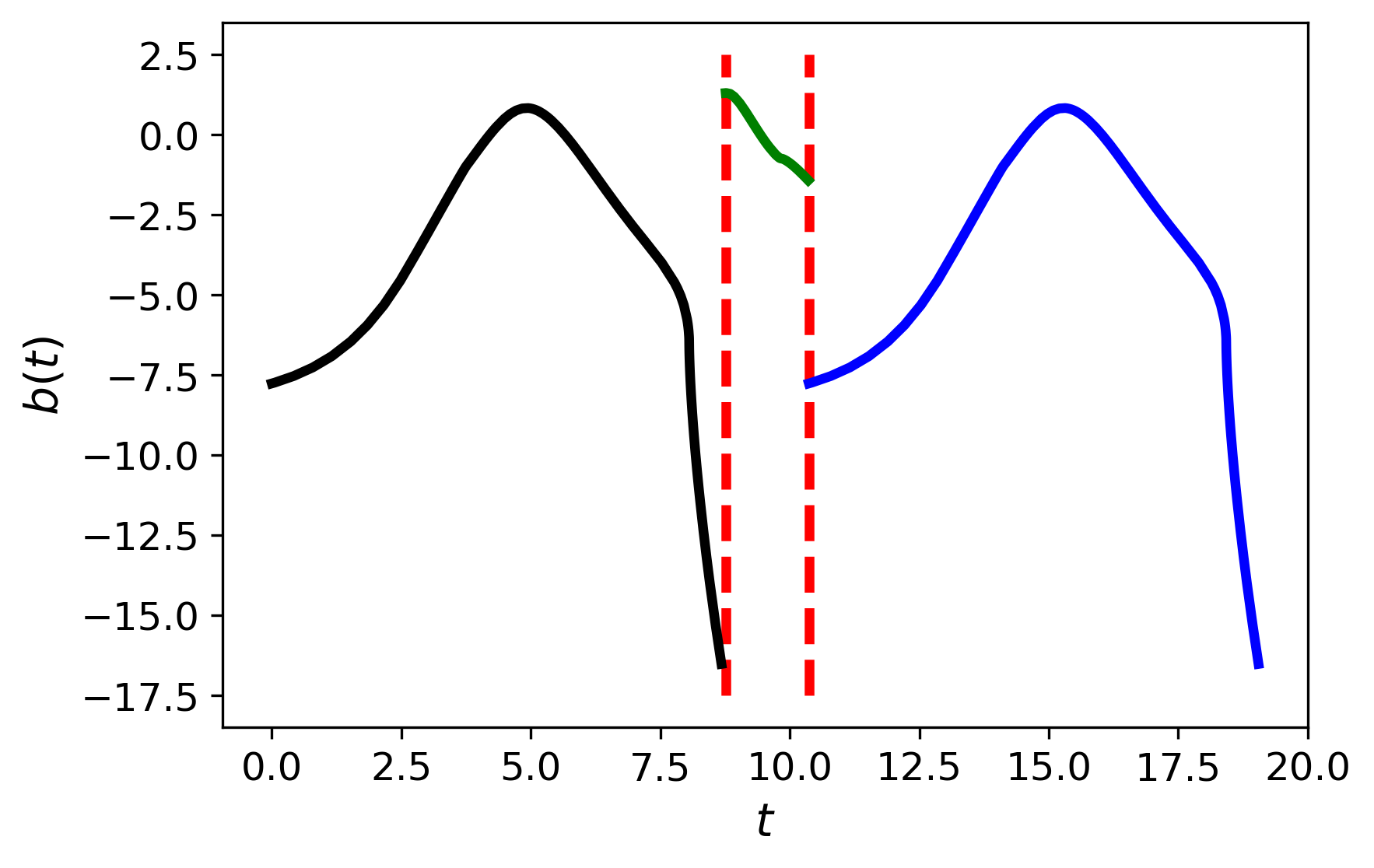}
	\caption{Trajectory of a single leaf (top left), corresponding $\xi^2$ profile (top right) and optimized control functions $a(t)$ (bottom left) and $b(t)$ (bottom right). The black, green, and blue segments distinguish the three arcs that make up the triplet. The vertical red lines indicate the time points where the controls switch between arc types.}
	\label{fig:Triplet}
\end{figure}

\begin{figure} 
	\centering
    \includegraphics[height=3.2cm]{Figures/PointBased with Energy June 4/Arc1} \hfill	\includegraphics[height=3.2cm]{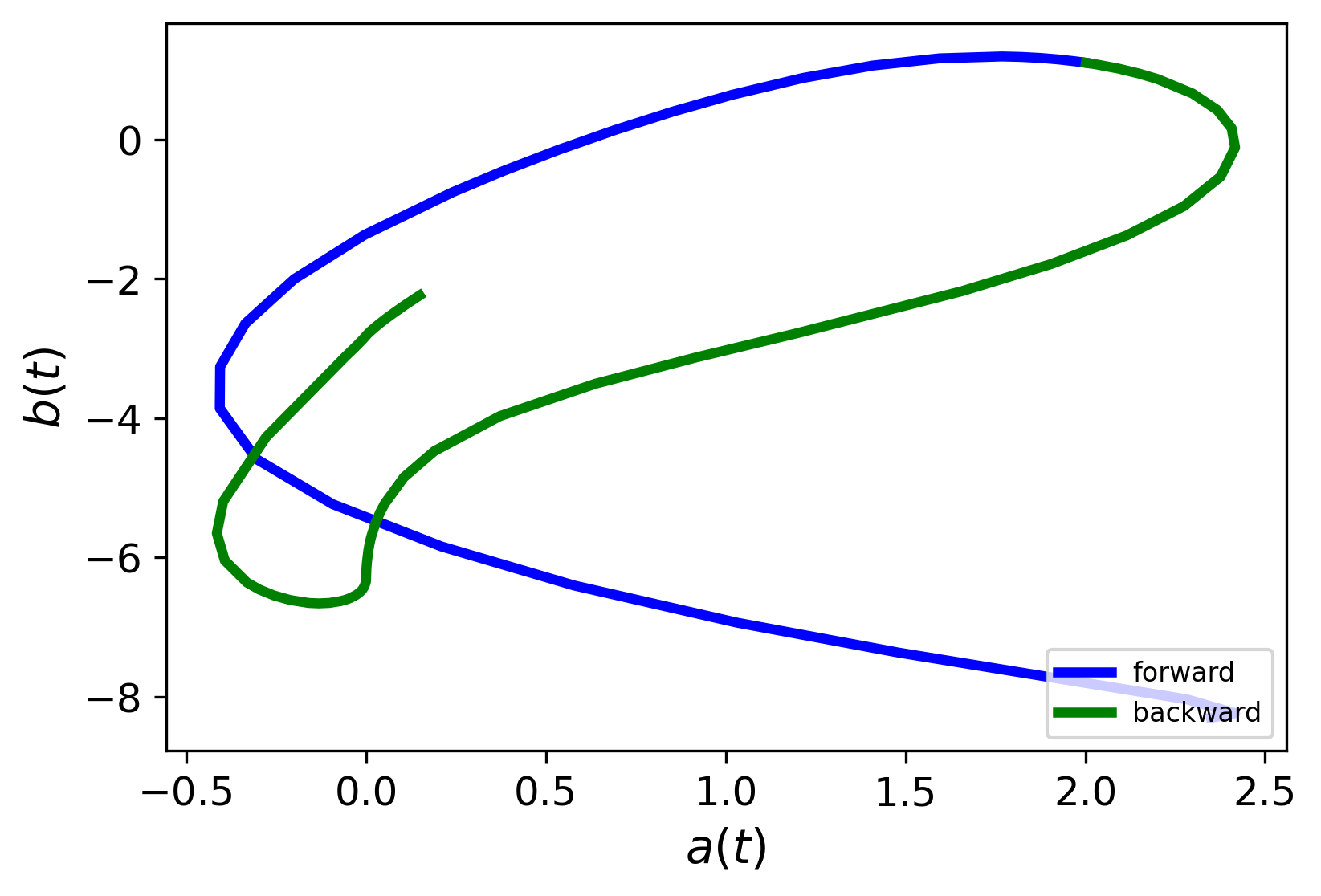}\hfill
    \includegraphics[height=3.2cm]{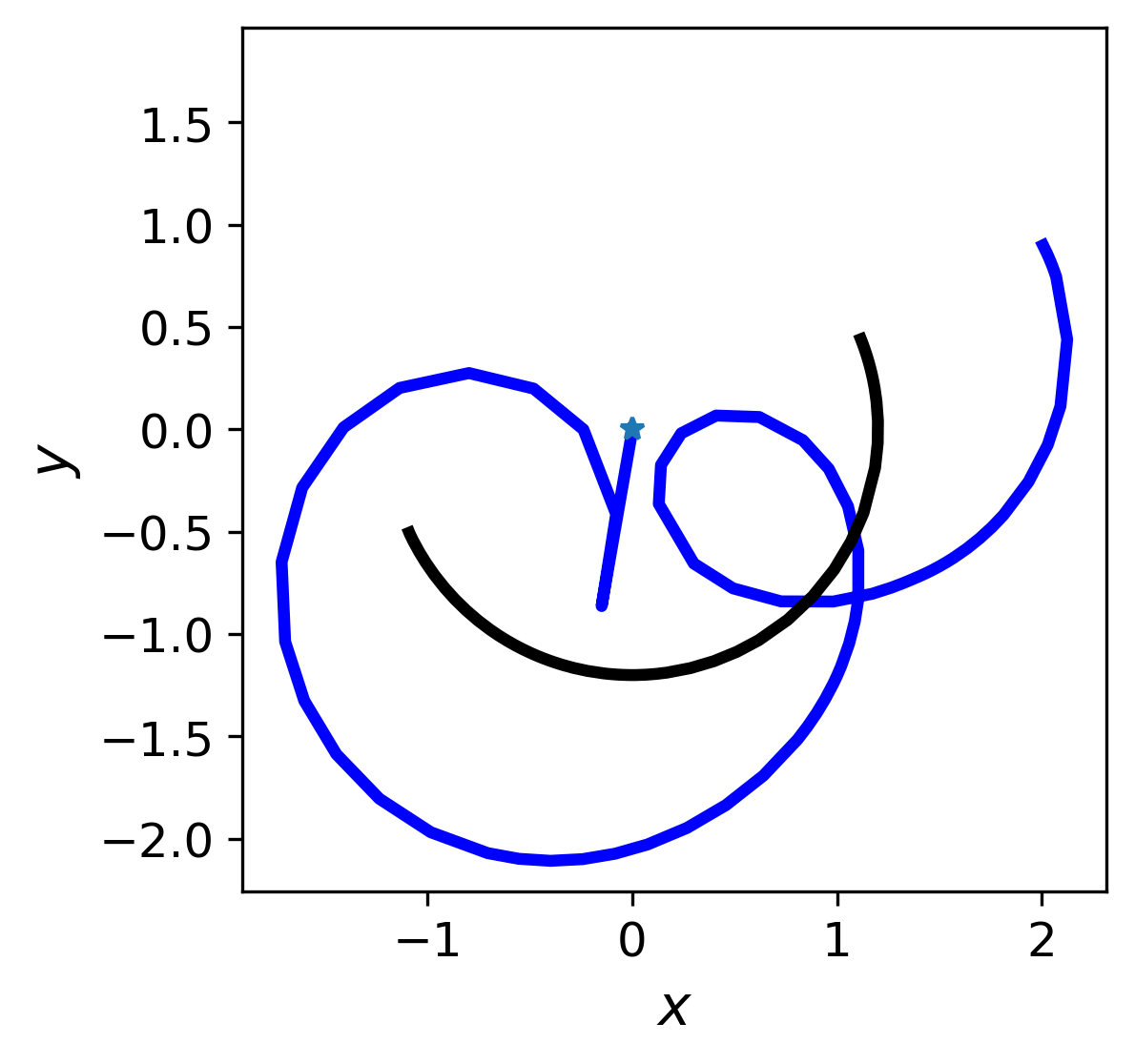}\\
	\includegraphics[height=3.2cm]{Figures/PointBased with Energy June 4/Arc2}\hfill	\includegraphics[height=3.2cm]{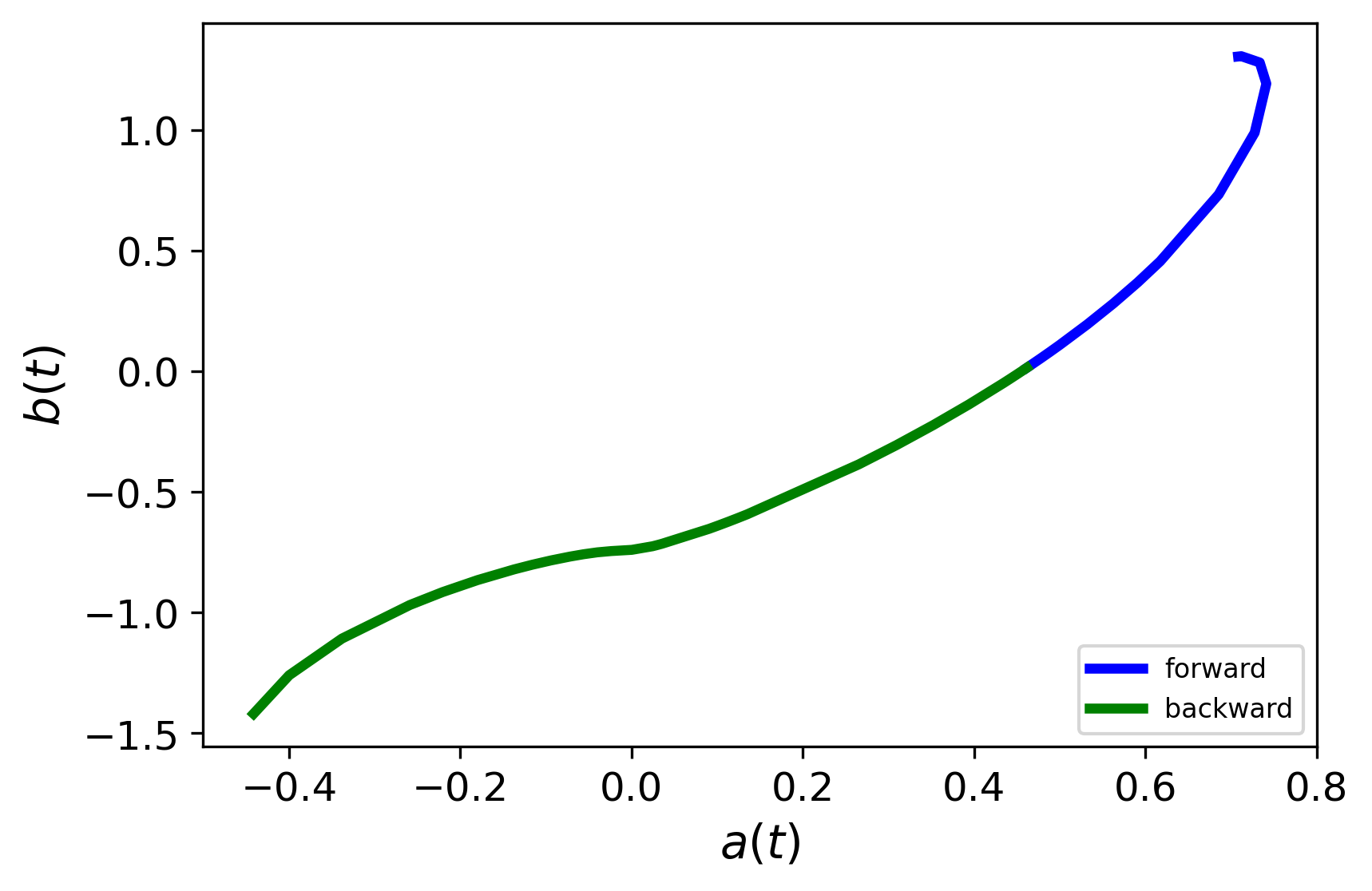}\hfill
	\includegraphics[height=3.2cm]{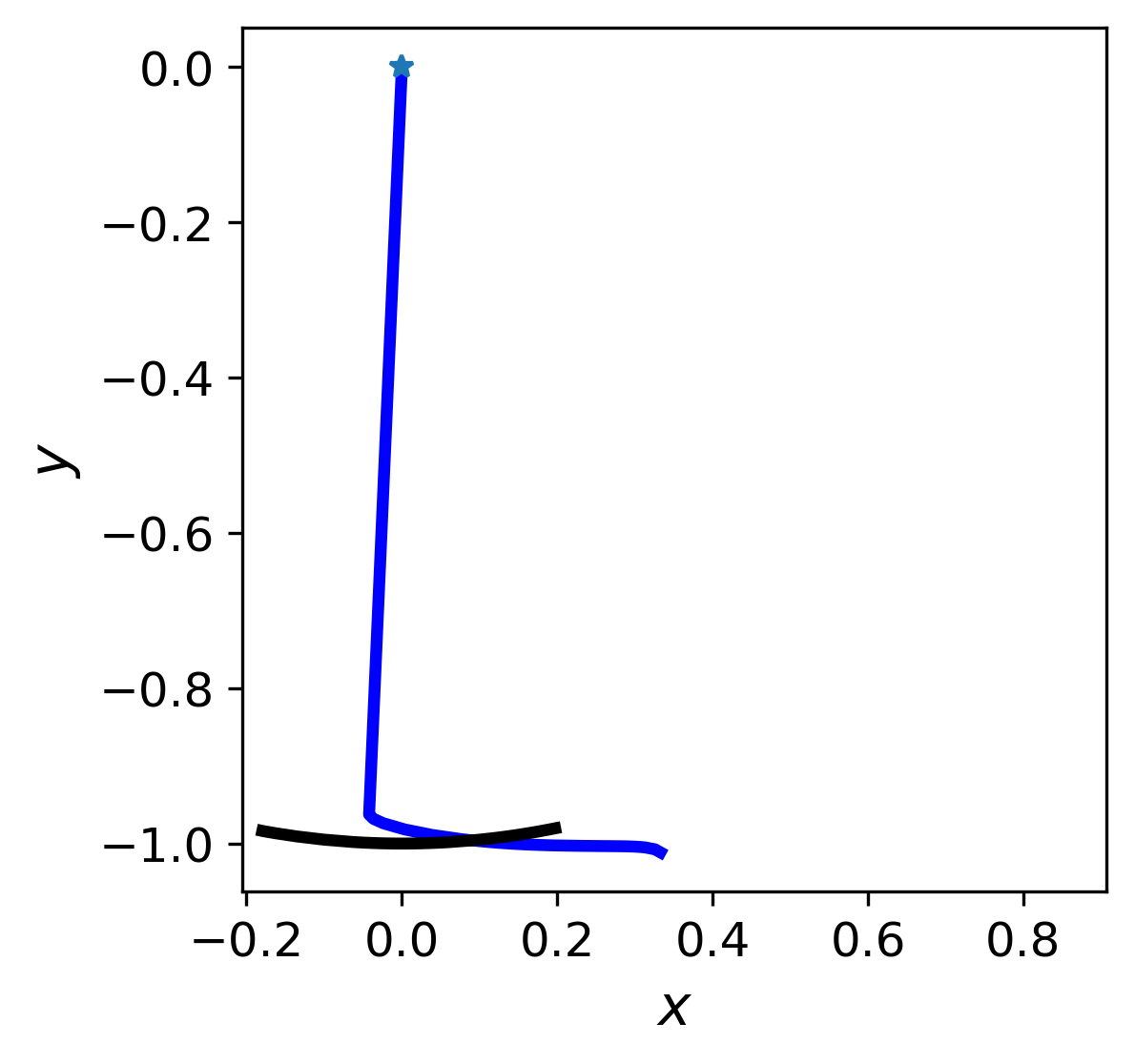}\\
	\includegraphics[height=3.2cm]{Figures/PointBased with Energy June 4/Arc3}\hfill	\includegraphics[height=3.2cm]{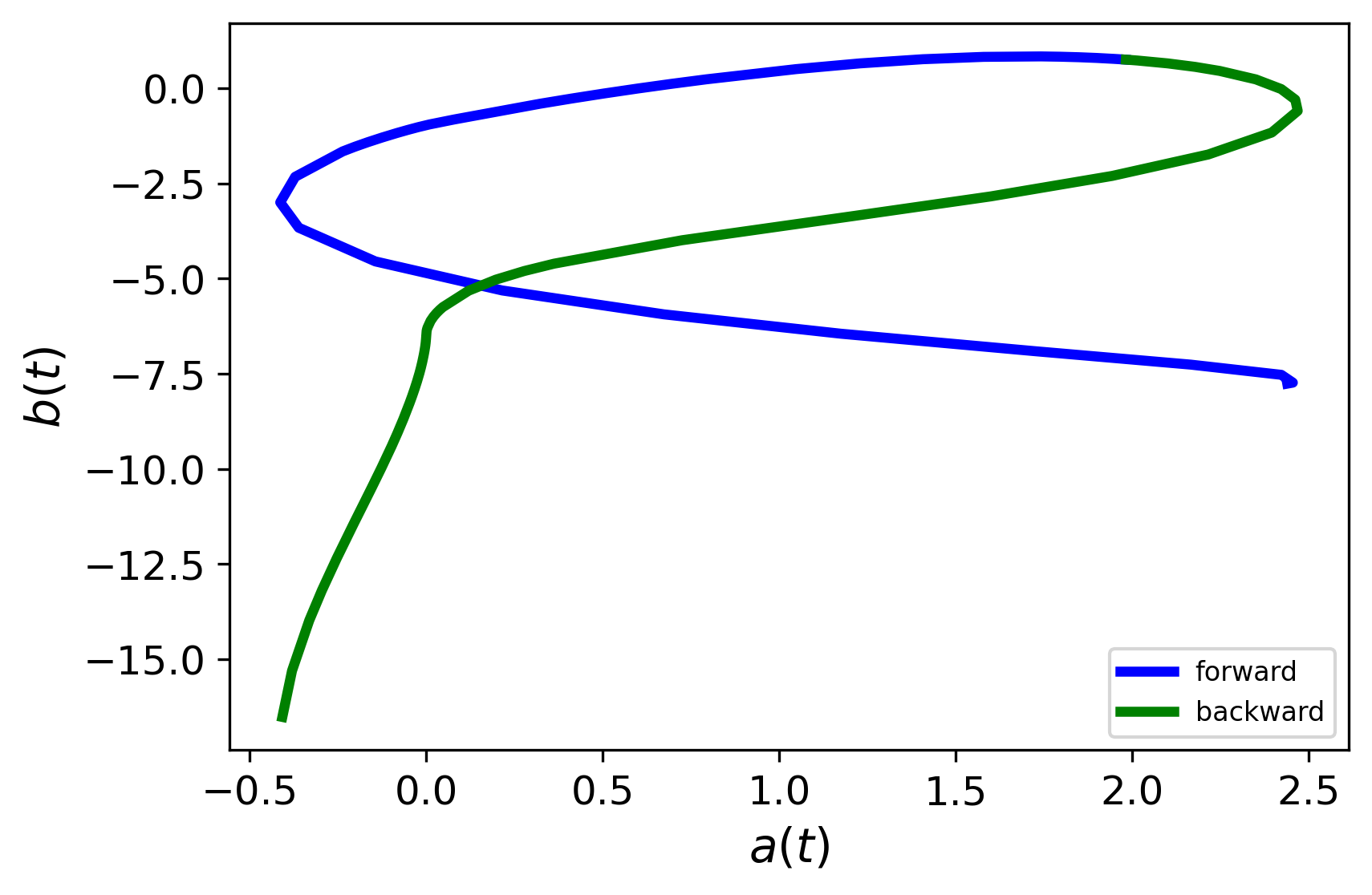}\hfill
	\includegraphics[height=3.2cm]{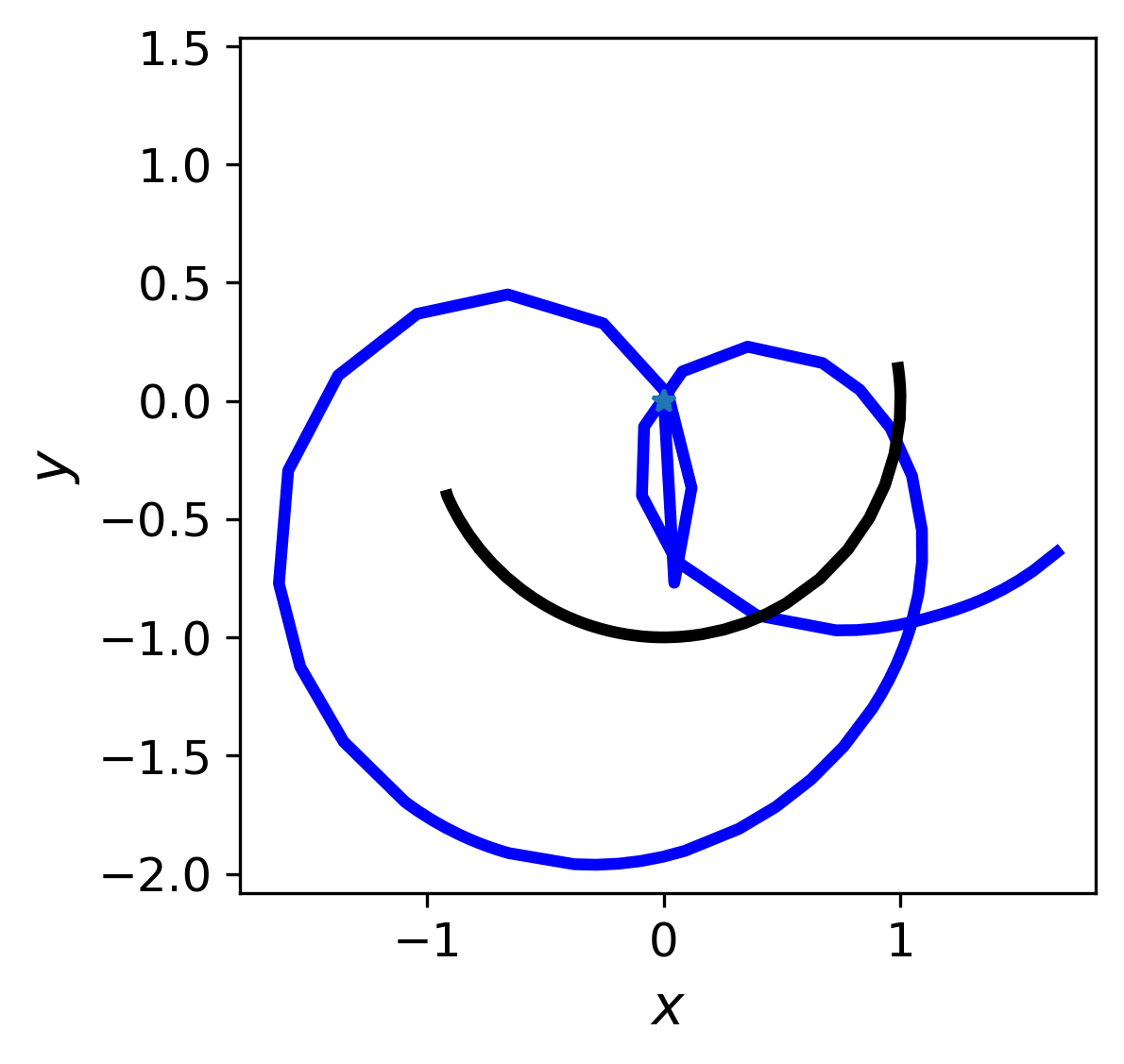}\\
	\caption{Point-optimized skate trajectories (left), control mass trajectories (middle), and overlay of skate and control mass trajectories in spatial frame (right) for arc 1 (top), arc 2 (middle), and arc 3 (bottom). Green and blue indicate the forward- and backward-in-time solutions, respectively. In the far right column, the star indicates the initial point of the control mass.}
	\label{fig:ArcControls}
\end{figure}

\begin{figure} 
\begin{minipage}{0.49\textwidth}
    \begin{flushright}
    \includegraphics[scale=0.4]{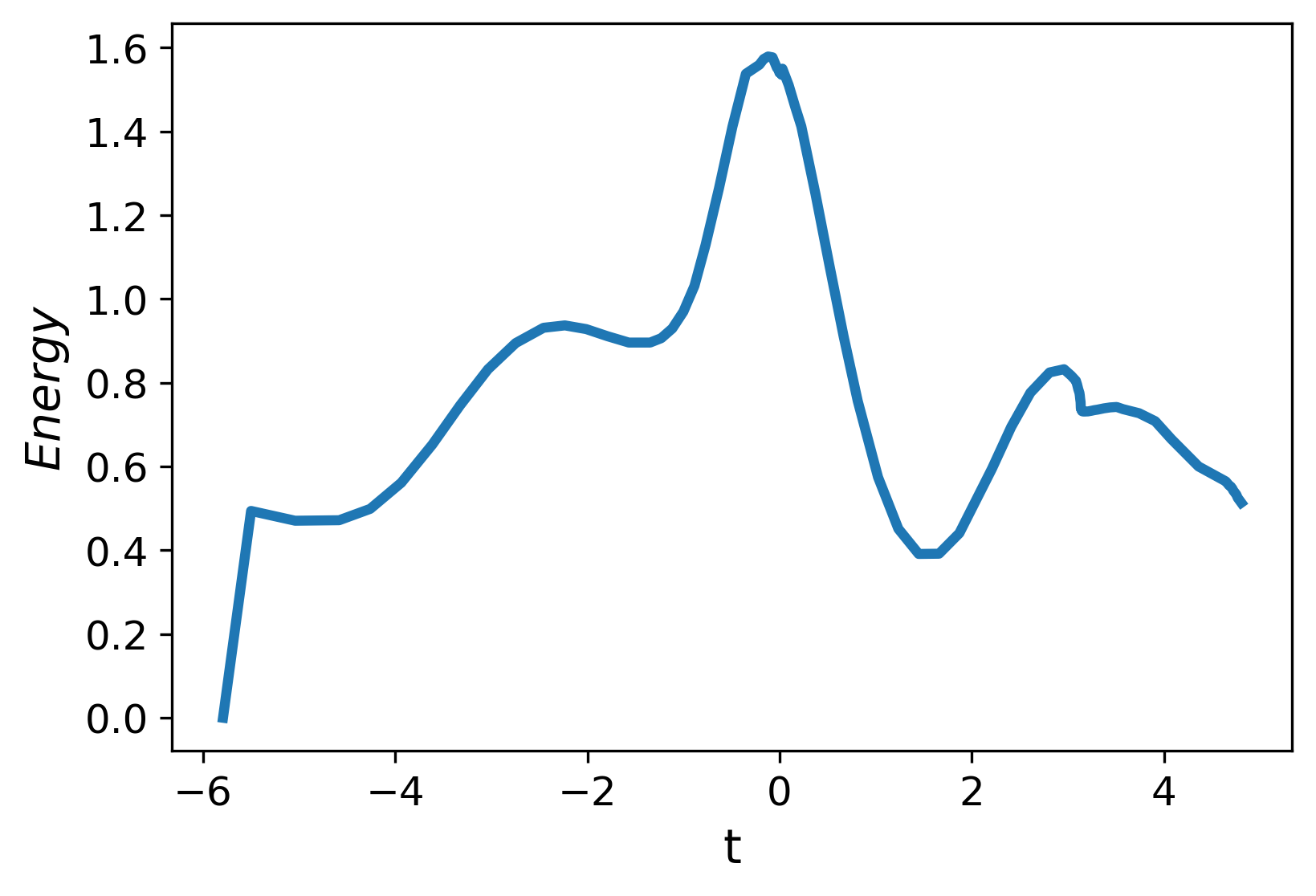}\\
    \includegraphics[scale=0.4]{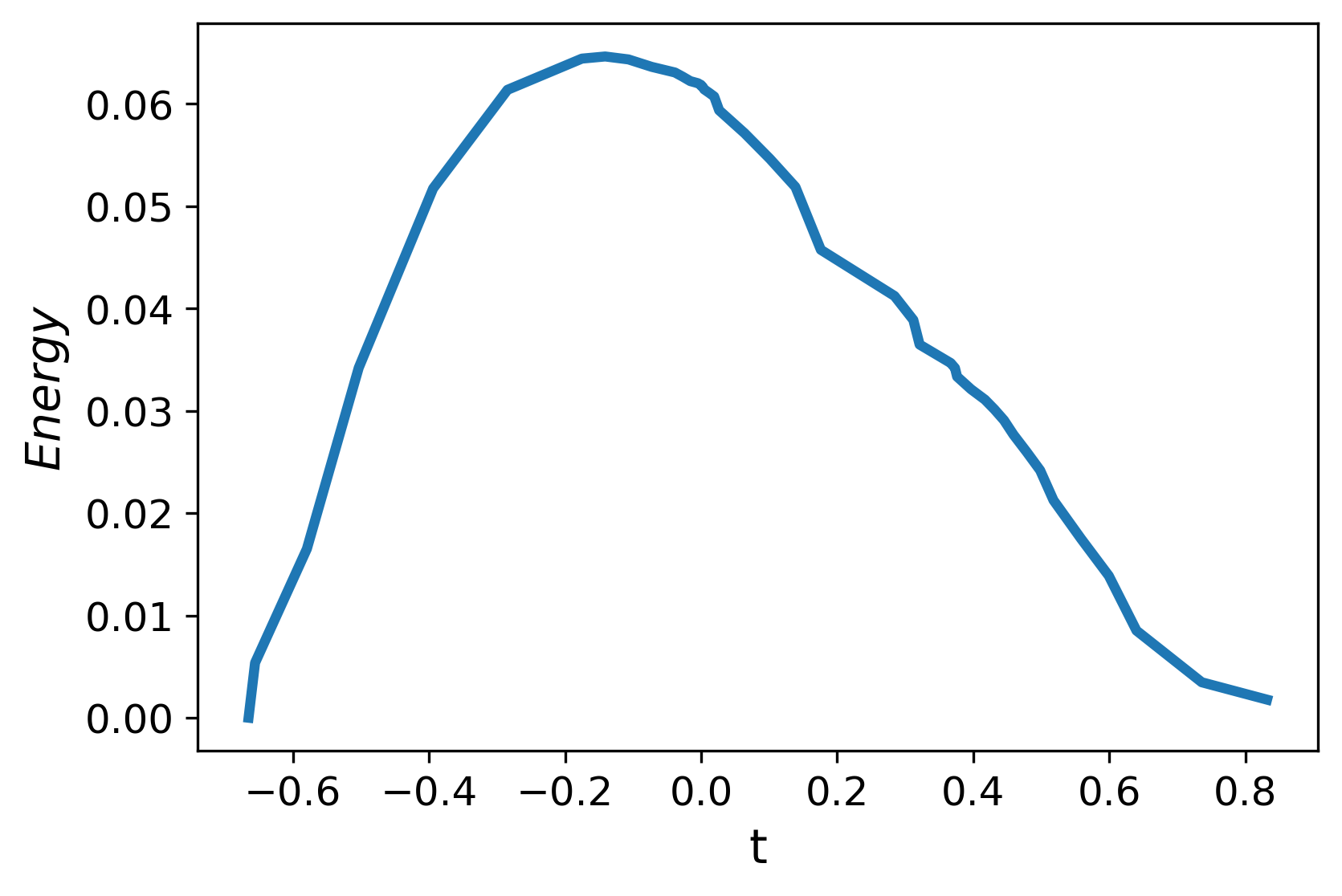}\\
    \includegraphics[scale=0.4]{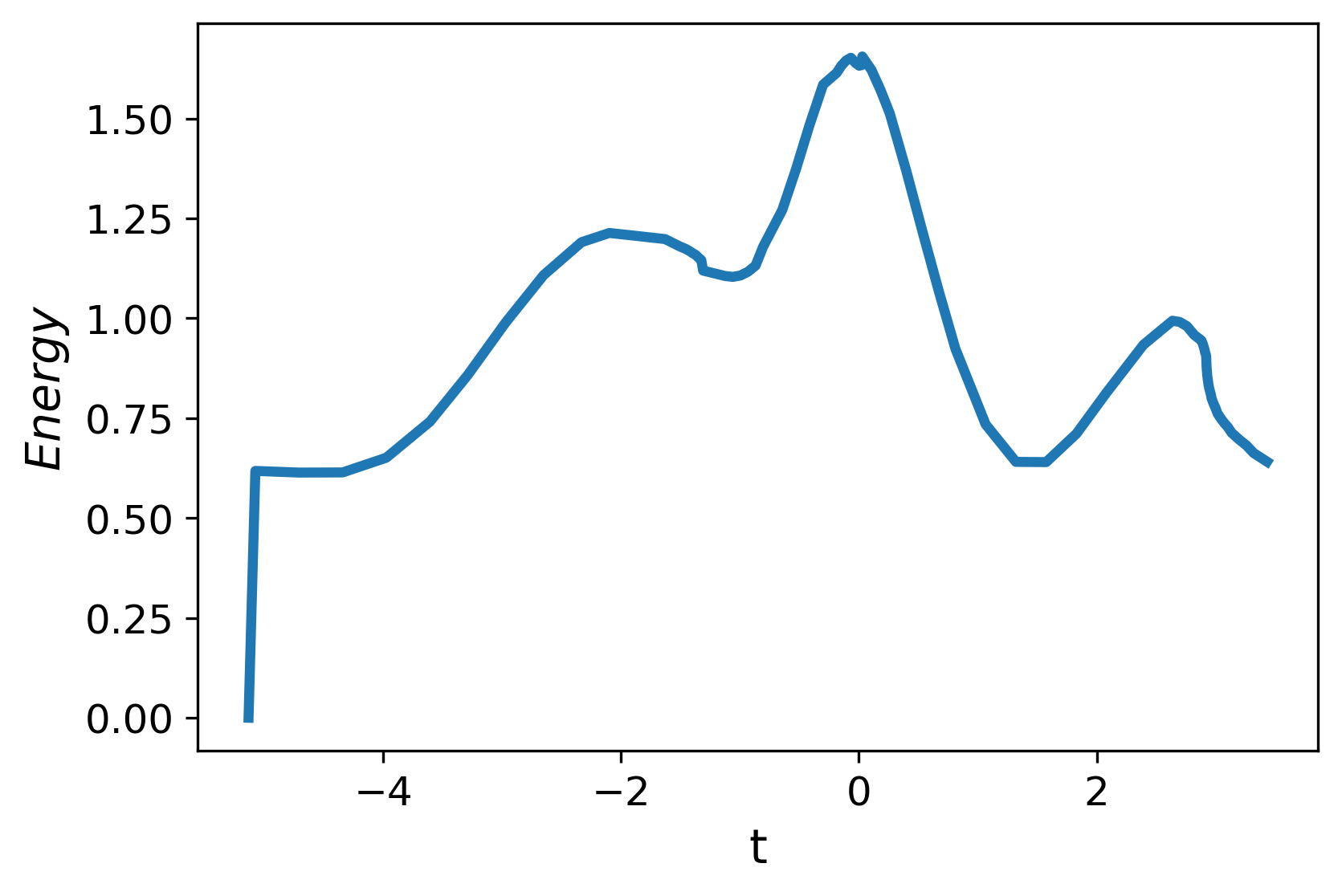}
    \end{flushright}
\end{minipage}
\begin{minipage}{0.49\textwidth}
    \begin{flushright}
    \includegraphics[scale=0.4]{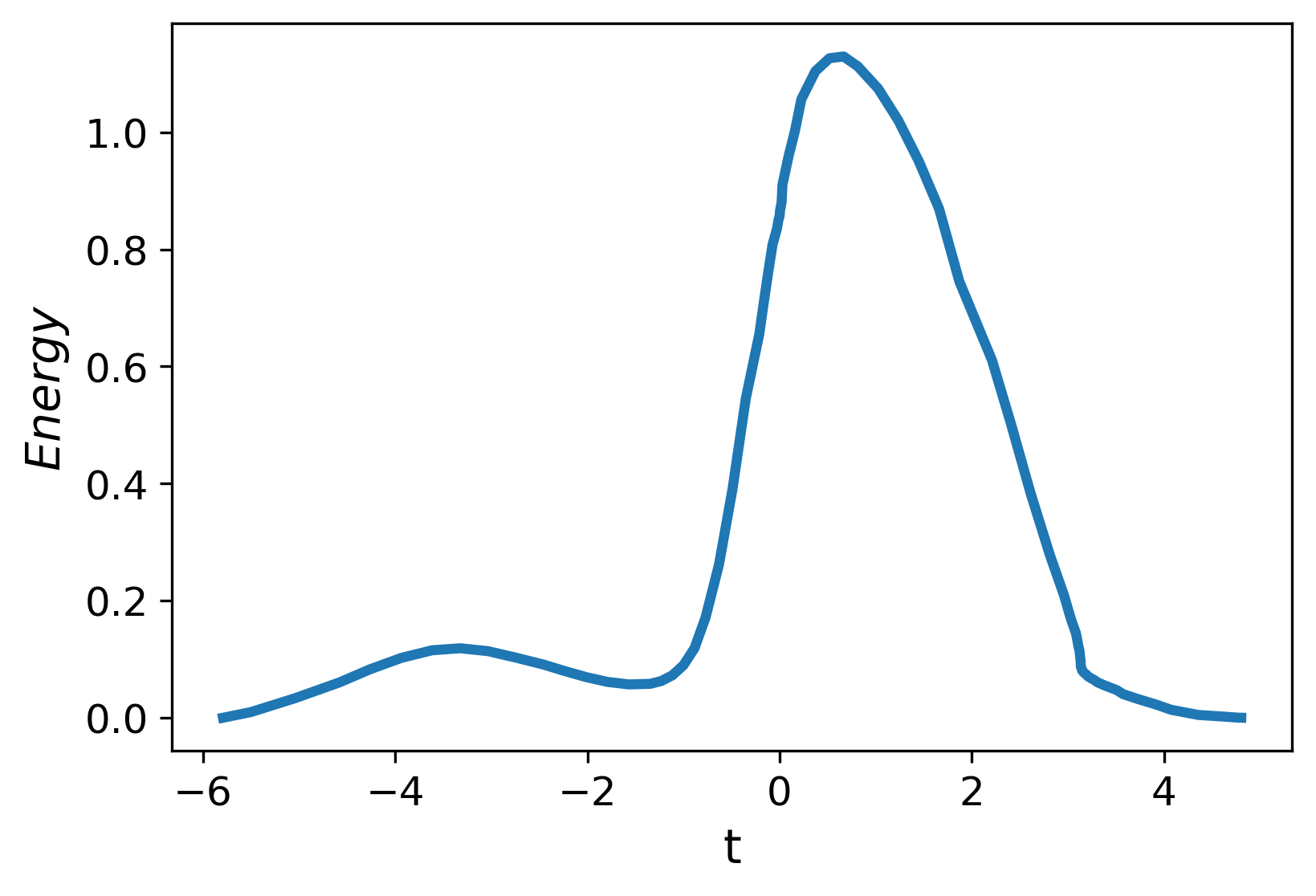}\\
    \includegraphics[scale=0.4]{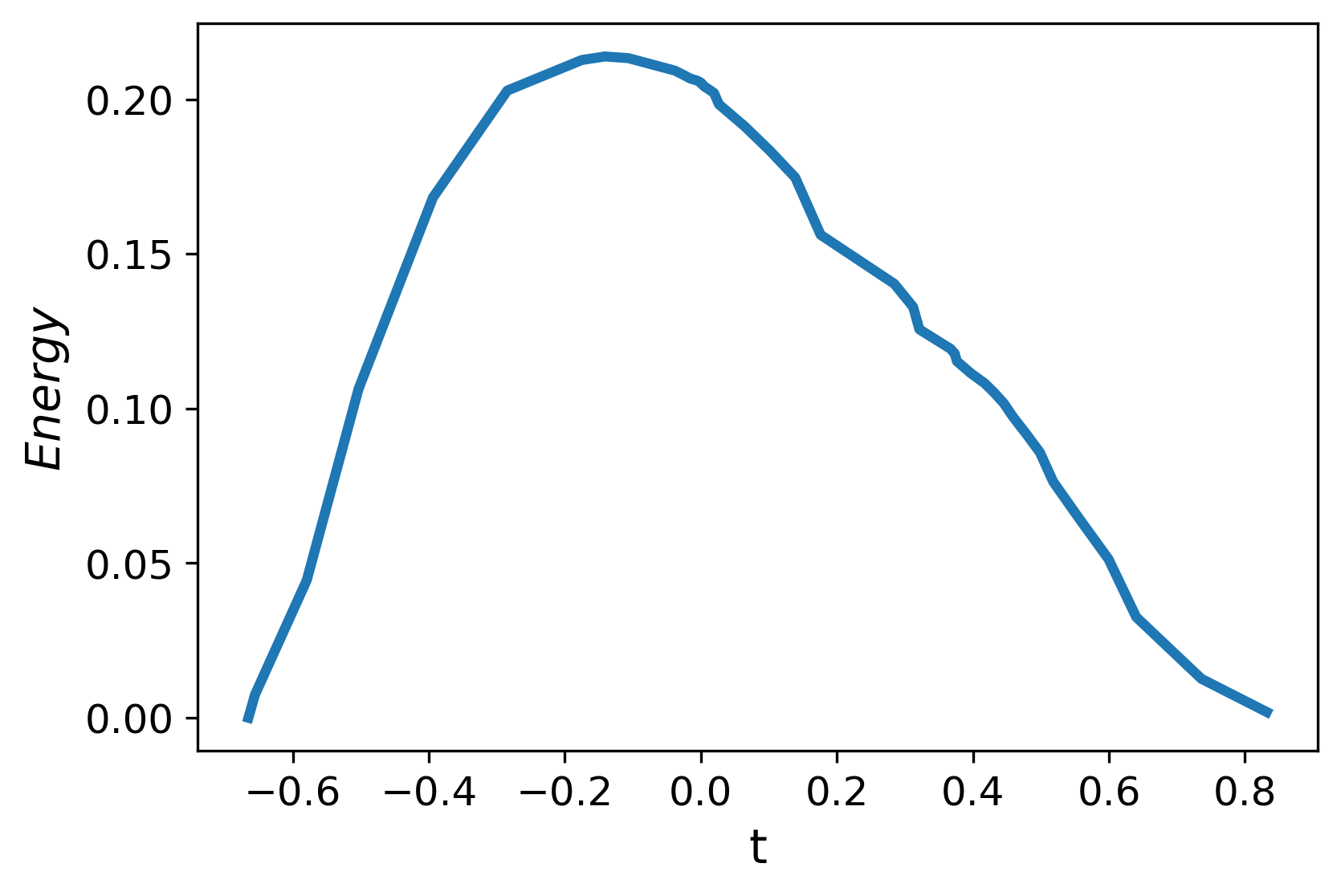}\\
    \includegraphics[scale=0.4]{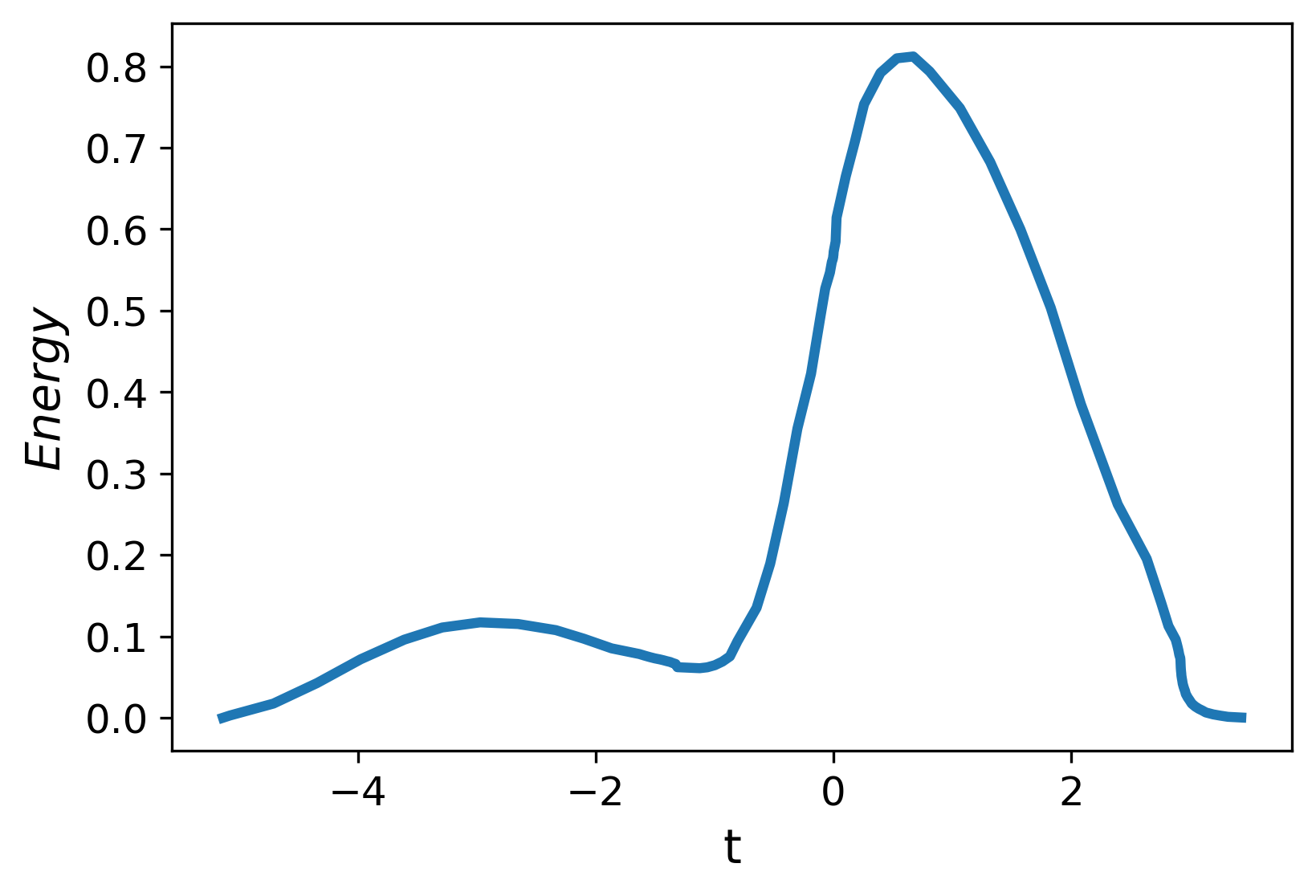}
    \end{flushright}
\end{minipage}
    \caption{Energy over time for control mass (left) and skate (right) for arc 1 (top), arc 2 (middle), and arc 3 (bottom) with point-based optimization.}
    \label{fig: Energy pointbased}
\end{figure}

\begin{figure} 
    \centering
    \includegraphics[width=0.5\textwidth]{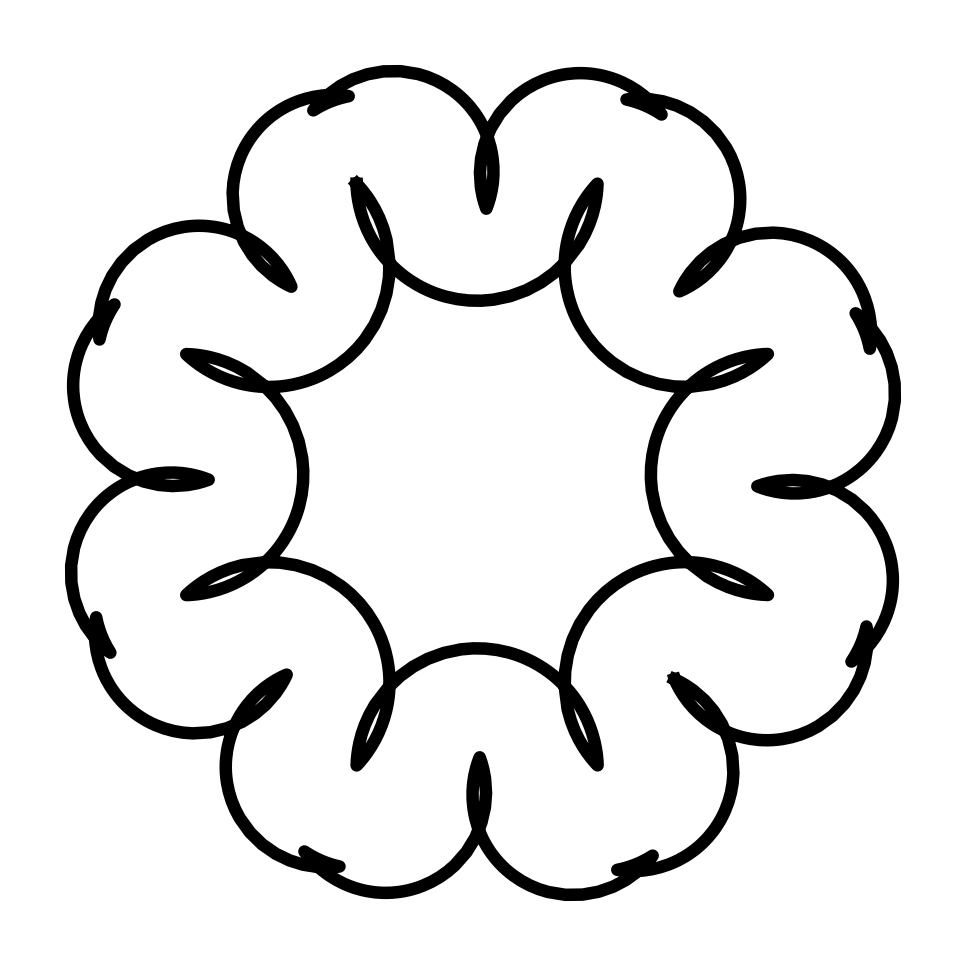}
	\caption{Full pattern reconstruction using point-based optimization. Notice a clear resemblance to Figure~\ref{fig:SimulatedPattern_Length_NewControl}, with slight difference in the length and curvature of the arcs composing the pattern.  }
	\label{fig:SimulatedPattern}
\end{figure}

As shown in Figure \ref{fig:ArcEndPoints} the error from the optimization of the forward-in-time arcs with respect to the target endpoints $(X_{opt}, Y_{opt})$ (indicated by the red dots) varies. The optimization error could be further minimized by choosing optimal simulation parameters, but for the purposes of this work, we were satisfied the given results. Choosing optimal simulation parameters that would yield smaller optimization errors could likely be automated using a machine learning algorithm.
The simulated arc 1, as well as the optimized control functions and $\xi^2$ profile, are shown in Figure \ref{fig:Arc1}.
The forward- and backward-in-time solutions are represented by the green and blue segments, respectively. Similarly, the simulations for arc 2 and arc 3, as well as the optimized control functions and profiles of $\xi^2$, are shown in Figures \ref{fig:Arc2} and \ref{fig:Arc3}. 
The combined triplet simulation, optimized control functions, and corresponding $\xi^2$ profile are shown in Figure \ref{fig:Triplet}.
Figure \ref{fig:ArcControls} shows the trajectory of the position (in body frame) of the added mass $m$ next to the blade trajectory (in spatial frame) for each arc, as well as the overlay of skate and control mass trajectories in the spatial frame. Again, the green indicates the forward-in-time solution, while the blue indicates the backward-in-time solution. Thus, the path starts at the beginning of the blue portion of the curve and ends at the end of the green portion of the curve. 
The total energy of the skate and control mass are shown in Figure \ref{fig: Energy pointbased}. In this case, the energy remains bounded. 
The full pattern produced by transforming (rotation and translation) and repeating the simulated arcs is shown in Figure \ref{fig:SimulatedPattern}, which is very similar to Figure~\ref{fig:SimulatedPattern_Length_NewControl}, but not exactly the same, as the length of the individual arcs is slightly different. Also, in this case, the radius of curvature of the arcs had to be changed slightly to facilitate the convergence, as well as the regularity of solution due to \eqref{sing_cond}. 

\subsection{Length-based Optimization}
\label{sect: LengthBased} 

In this section, we use the more general control function \eqref{eqn: general control function} in combination with the length-based cost function \eqref{eqn:CLengthandSpeed}. The length-based cost function optimizes the control parameters so that the total trajectory (combined backward- and forward-in-time solutions) reaches the desired length. Thus, unlike the point-based cost function, both the forward- and backward-in-time trajectories are solved within the optimization function and are considered when determining the optimal control parameters. 

In the example presented here, the condition on the optimal solution $a(t)$ requiring finite velocity of the control mass \eqref{sing_cond} is violated, and infinite velocity and energy of the control mass are observed. 

The simulation parameters are given in Table \ref{tab:NumericsPars_Length}. The initial control parameter guess, resulting optimized control parameters, and optimization errors are shown in Table \ref{tab:OptimizedPars_Length}. As indicated by the small values for optimization error in Table \ref{tab:OptimizedPars_Length}, the length-based optimization was able to successfully produce arcs that reached the target length. 
\begin{table} 
    \caption{Values used in simulations to create arc 1, arc2, and arc 3 with length-based optimization. In all cases $m=1$, $M=2$, and $I=3$.}
	\label{tab:NumericsPars_Length}
	\centering
	\begin{tabular}{llllll}
		\hline\noalign{\smallskip}
		Arc  & $T$ & $r$ & $L_{opt}$ &  $(\bar{p_1},\bar{p_2}, \bar{\theta},\bar{x},\bar{y},\bar{b})$ & $(a_1, a_2, a_3, \omega)$ \\
		\noalign{\smallskip}\hline\noalign{\smallskip}
		1 & $6$ & $1.2$& $1.1r\pi$ & $(2.5, 3, 0, 0, -1.2, 1.25)$ & $(1, 1, 1, 1)$\\
		2 & $2.1$ & $0.8$&  $0.2r\pi$ & $(0.0015, 0.0025, 0, 0, -0.8, 0.003)$ & $(0.1, 0.1, 0.1, 1)$\\
		3 & $8$& $1.0$ &  $r\pi$ & $(1.5, 2.5, 0, 0, -1, 0.75)$ & $(1, 1, 1, 1)$  \\ 
			\noalign{\smallskip}\hline
	\end{tabular}

\end{table}
\begin{table} 
	\caption{Initial guesses and optimized values for control parameters. }
	\label{tab:OptimizedPars_Length}
	\centering
	\begin{tabular}{llll}
	\hline\noalign{\smallskip}
		Arc  & Guess, $(a_1, a_2, a_3, \omega)$ & Optimized, $(a_1^*, a_2^*, a_3^*, \omega^*)$& Error\\
		\noalign{\smallskip}\hline\noalign{\smallskip}
		1 & $(1, 1, 1, 1)$ & $(1.000, 1.000, 1.000, 1.000)$& $6.916 \times 10^{-4}$\\
		2 & $(0.1, 0.1, 0.1, 1)$& $(0.658, -0.675,  0.428, 1.016)$ & $6.000 \times 10^{-14}$\\
		3 &  $(1, 1, 1, 1)$ & $(1.024, 1.007, 0.996, 1.009)$ & $4.655 \times 10^{-11}$\\ 
		\noalign{\smallskip}\hline
	\end{tabular}

\end{table}

\begin{figure} 
	\includegraphics[scale=0.4]{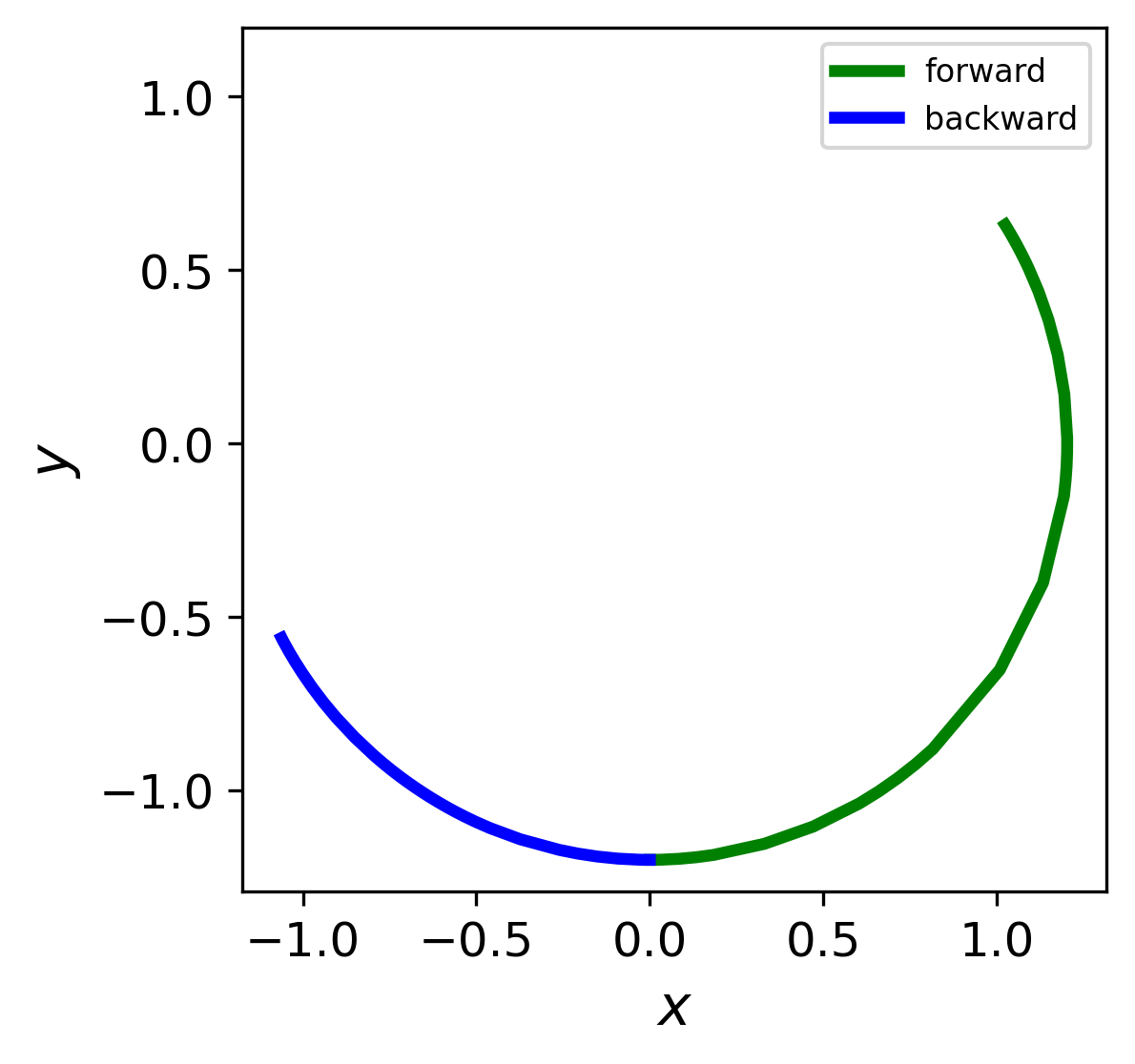} \hfill	\includegraphics[scale=0.4]{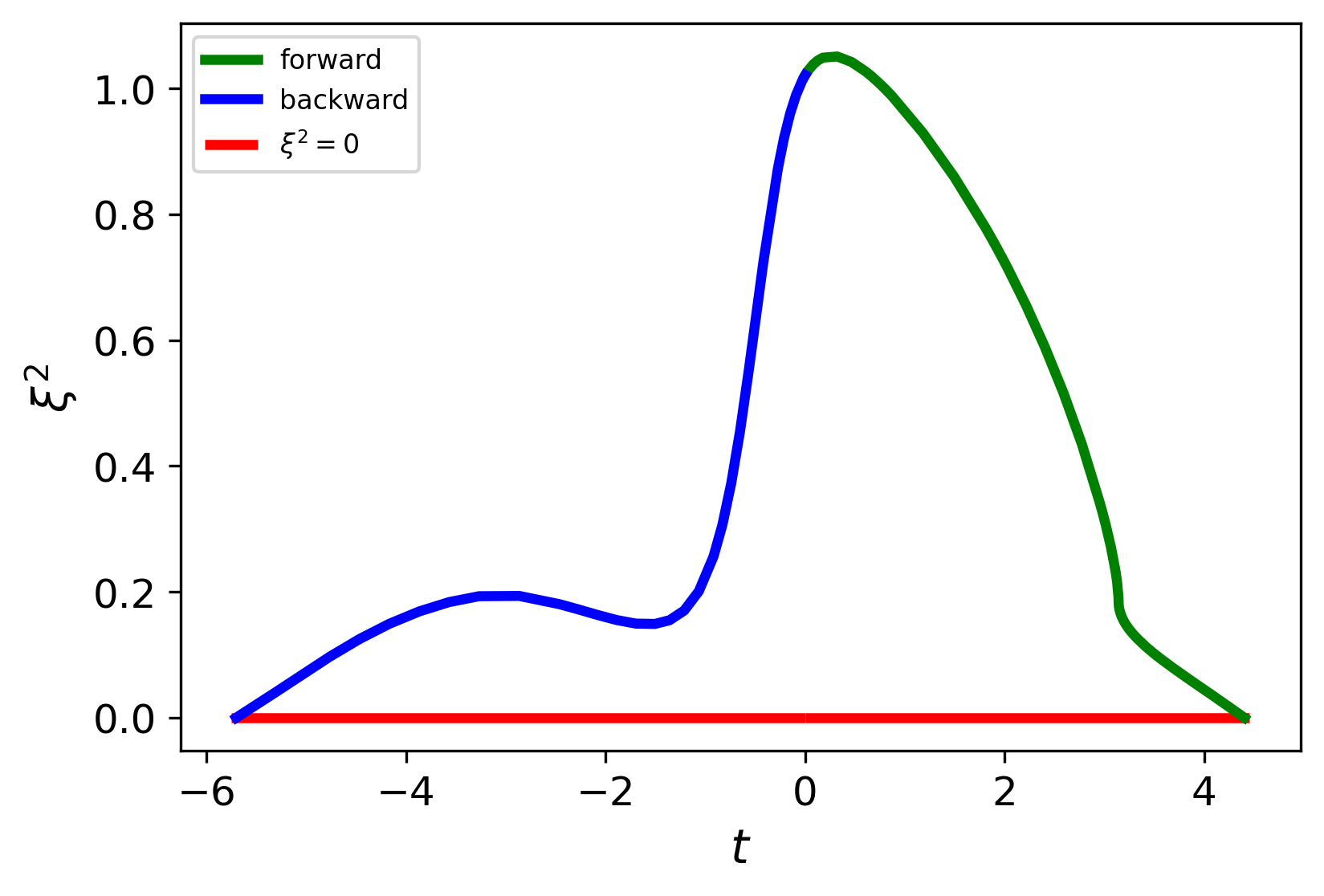}\hfill
	\includegraphics[scale=0.4]{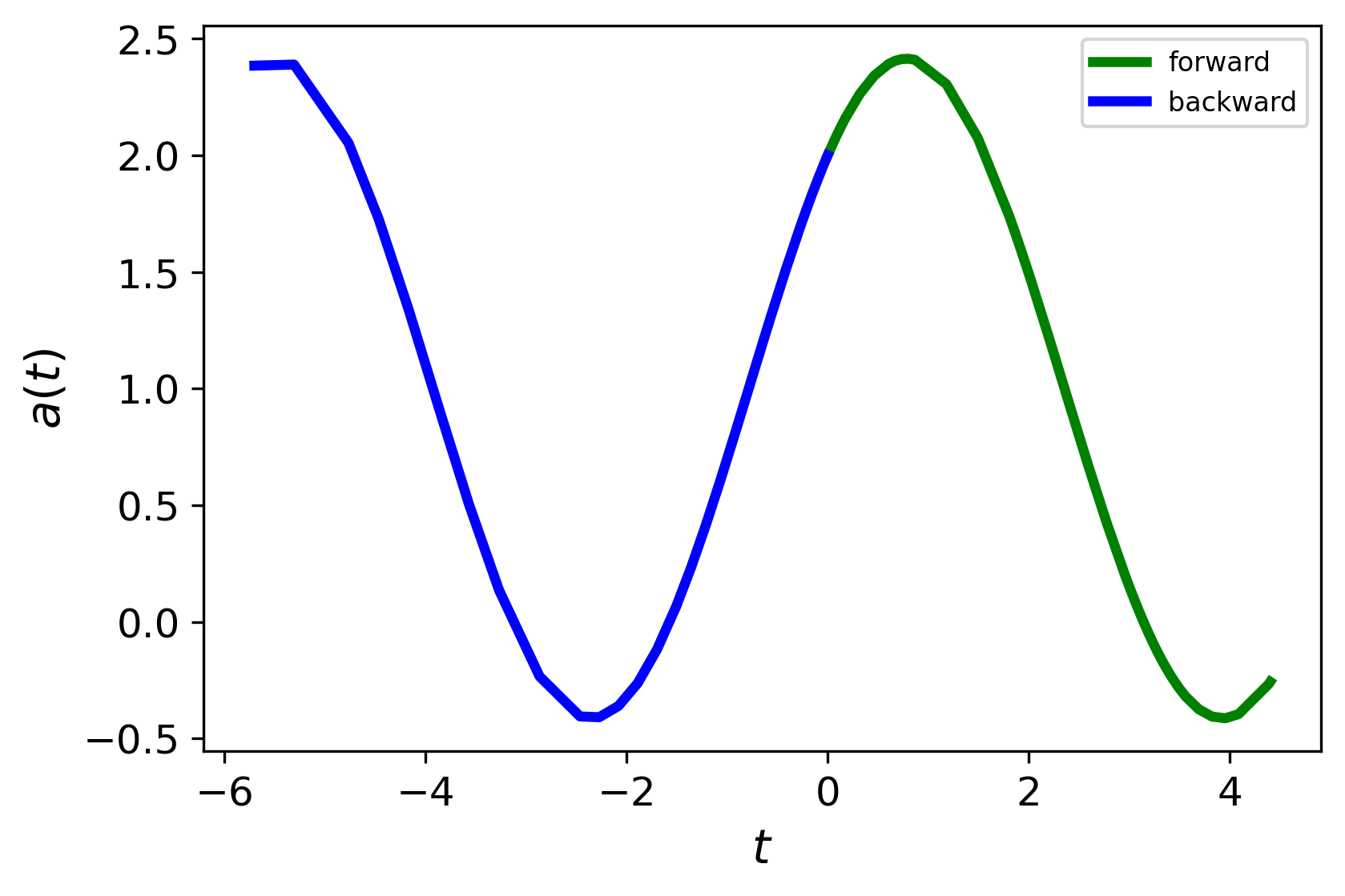}\hfill
	\includegraphics[scale=0.4]{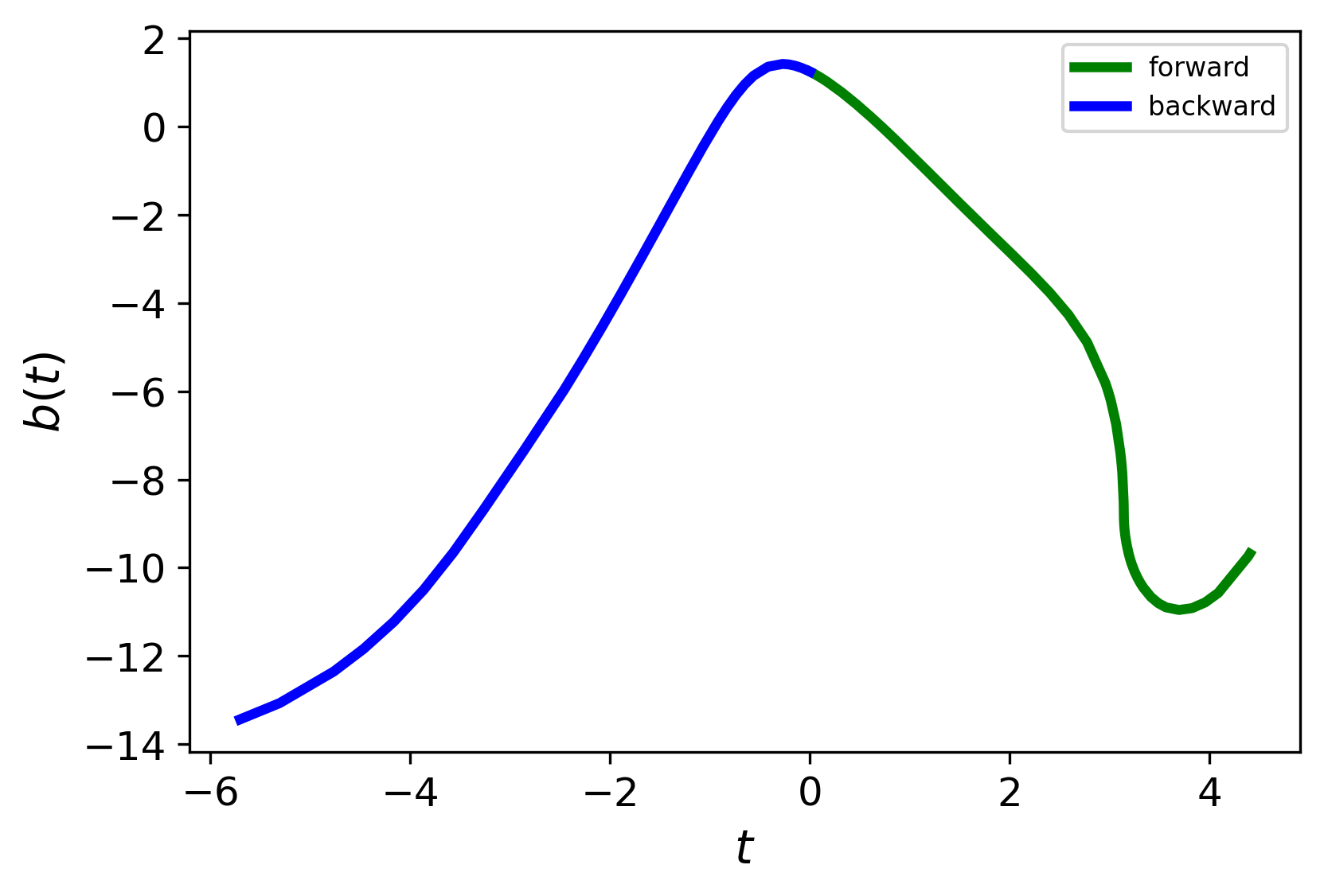}
	\caption{Arc 1 trajectory (top left), $\xi^2$ profile (top right), and optimized control functions, $a(t)$ (bottom left) and $b(t)$ (bottom right). Green and blue indicate the forward- and backward-in-time solutions, respectively.}	\label{fig:Arc1_Length}
\end{figure}
\begin{figure} 
	\includegraphics[scale=0.4]{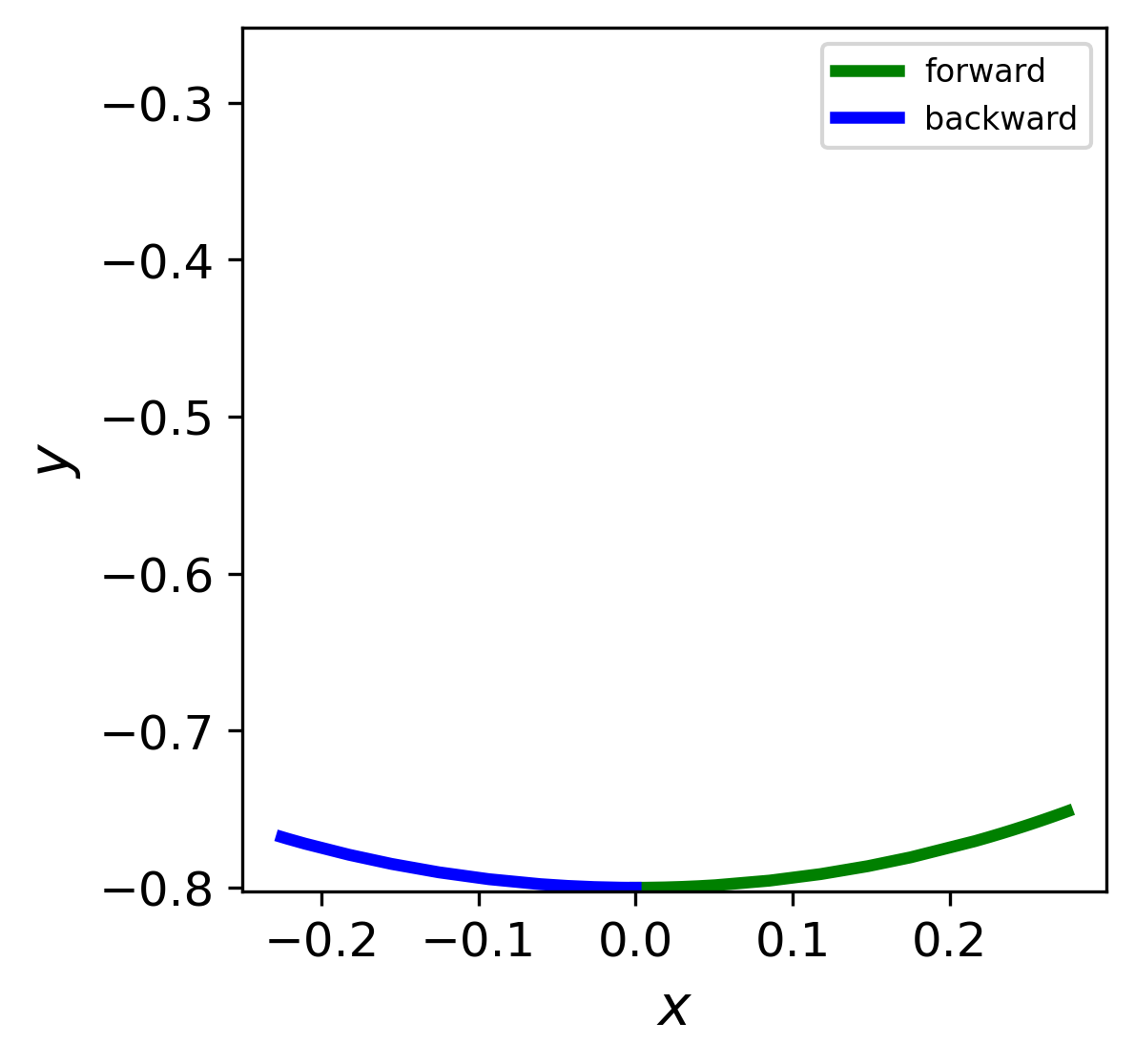}\hfill
	\includegraphics[scale=0.4]{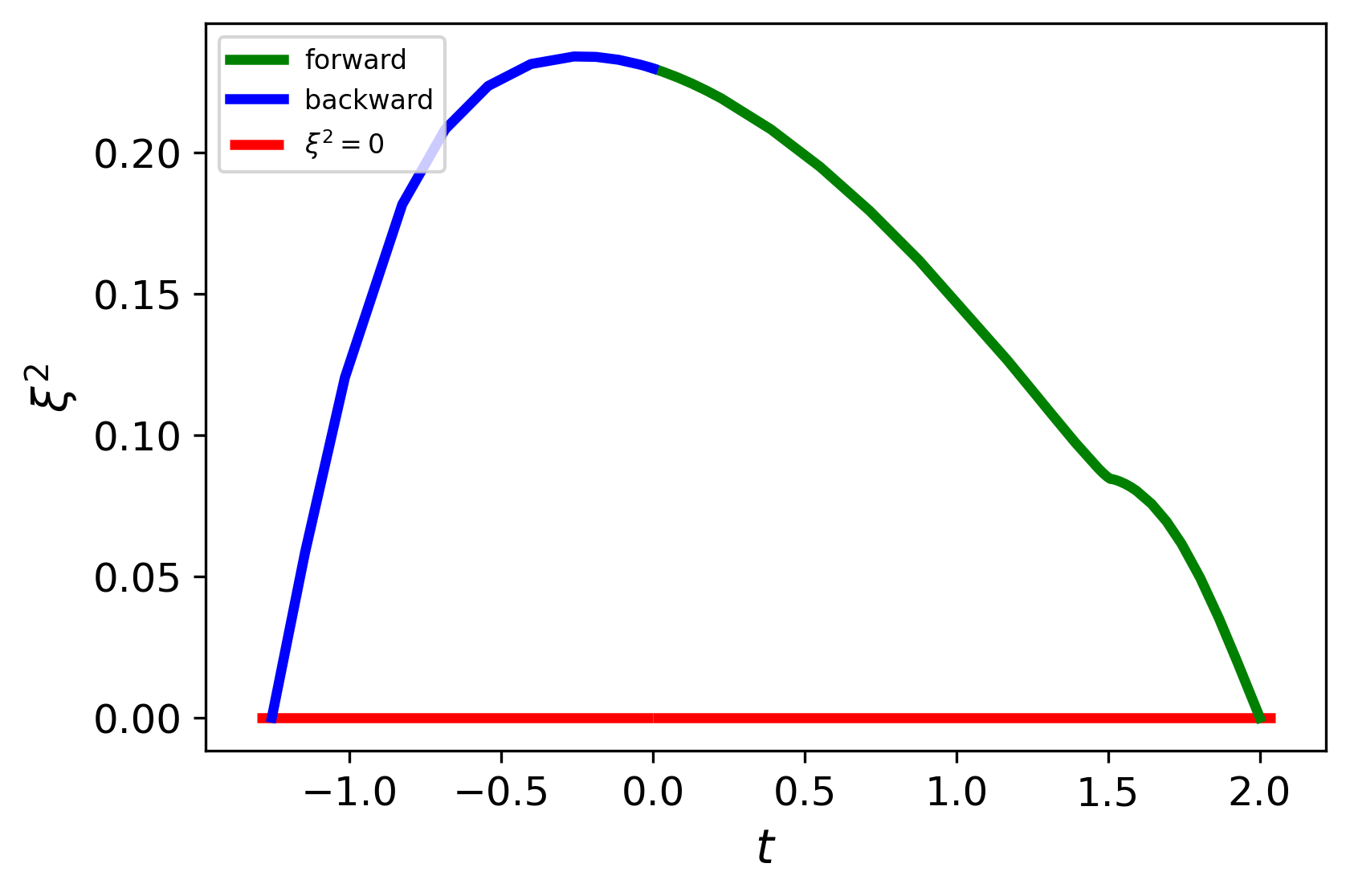}\hfill
	\includegraphics[scale=0.4]{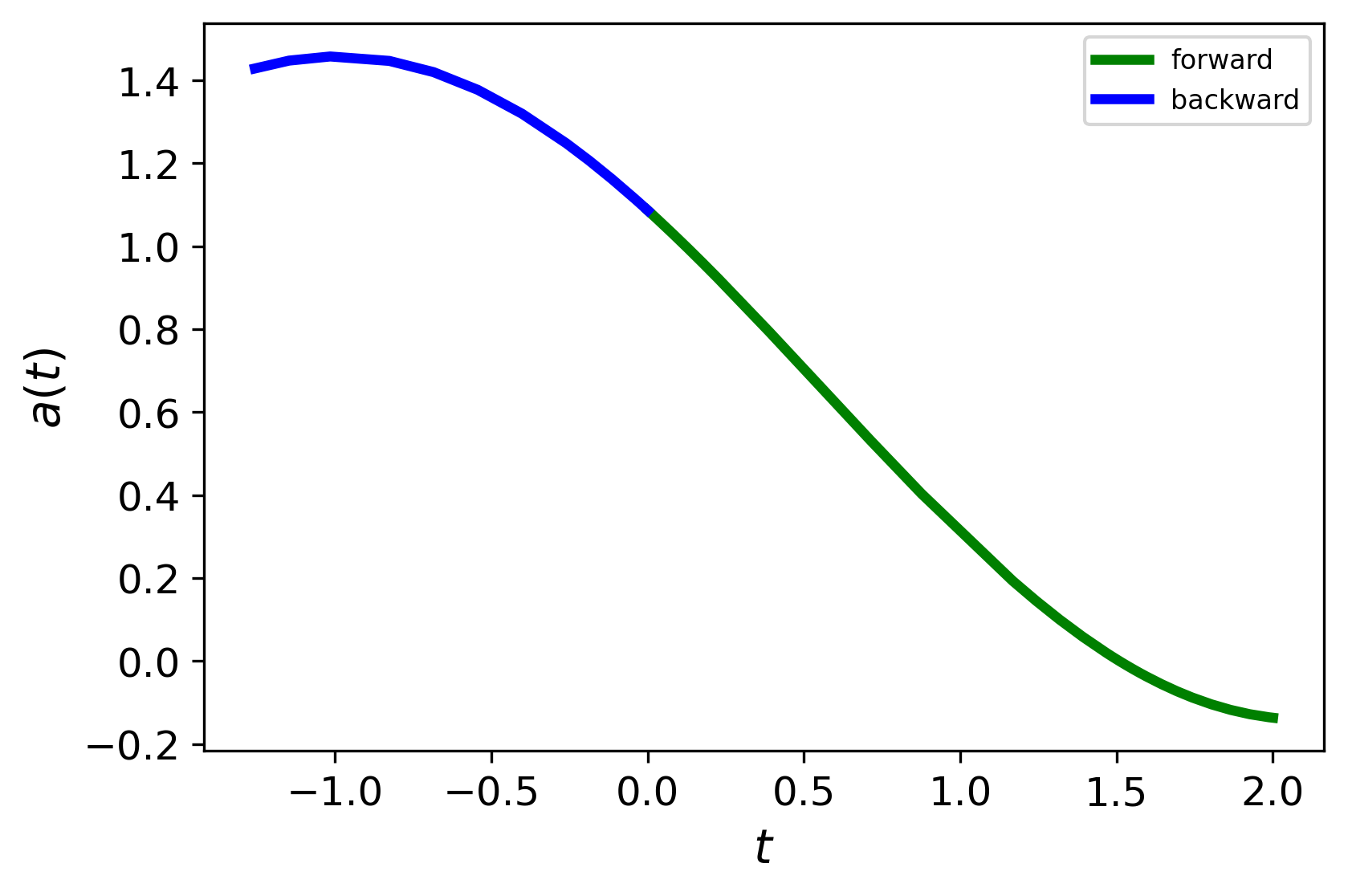}\hfill
	\includegraphics[scale=0.4]{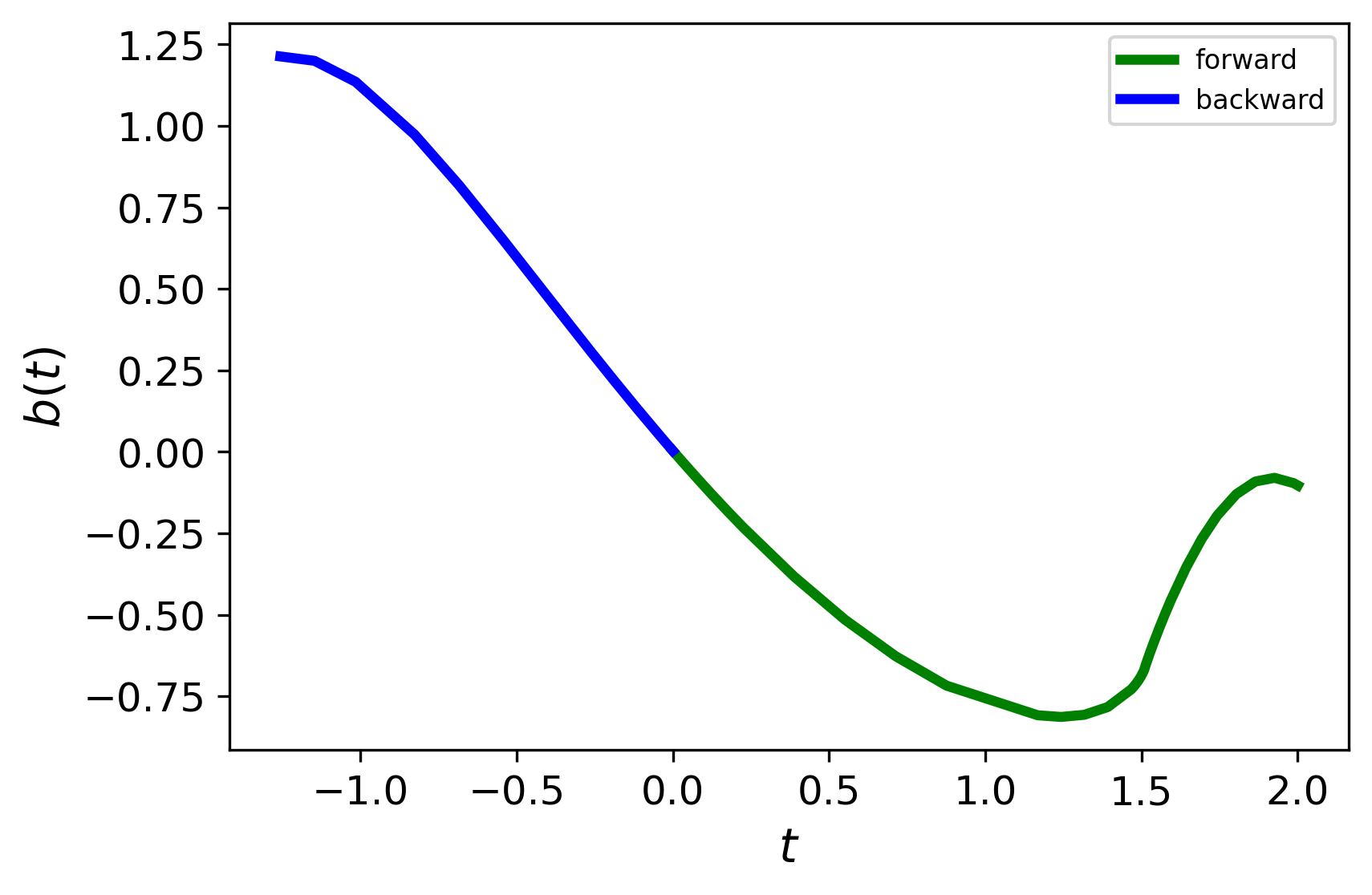}
	\caption{Arc 2 trajectory (top left), $\xi^2$ profile (top right), and optimized control functions, $a(t)$ (bottom left) and $b(t)$ (bottom right). Green and blue indicate the forward- and backward-in-time solutions, respectively.}	\label{fig:Arc2_Length}
\end{figure}
\begin{figure} 
	\includegraphics[scale=0.4]{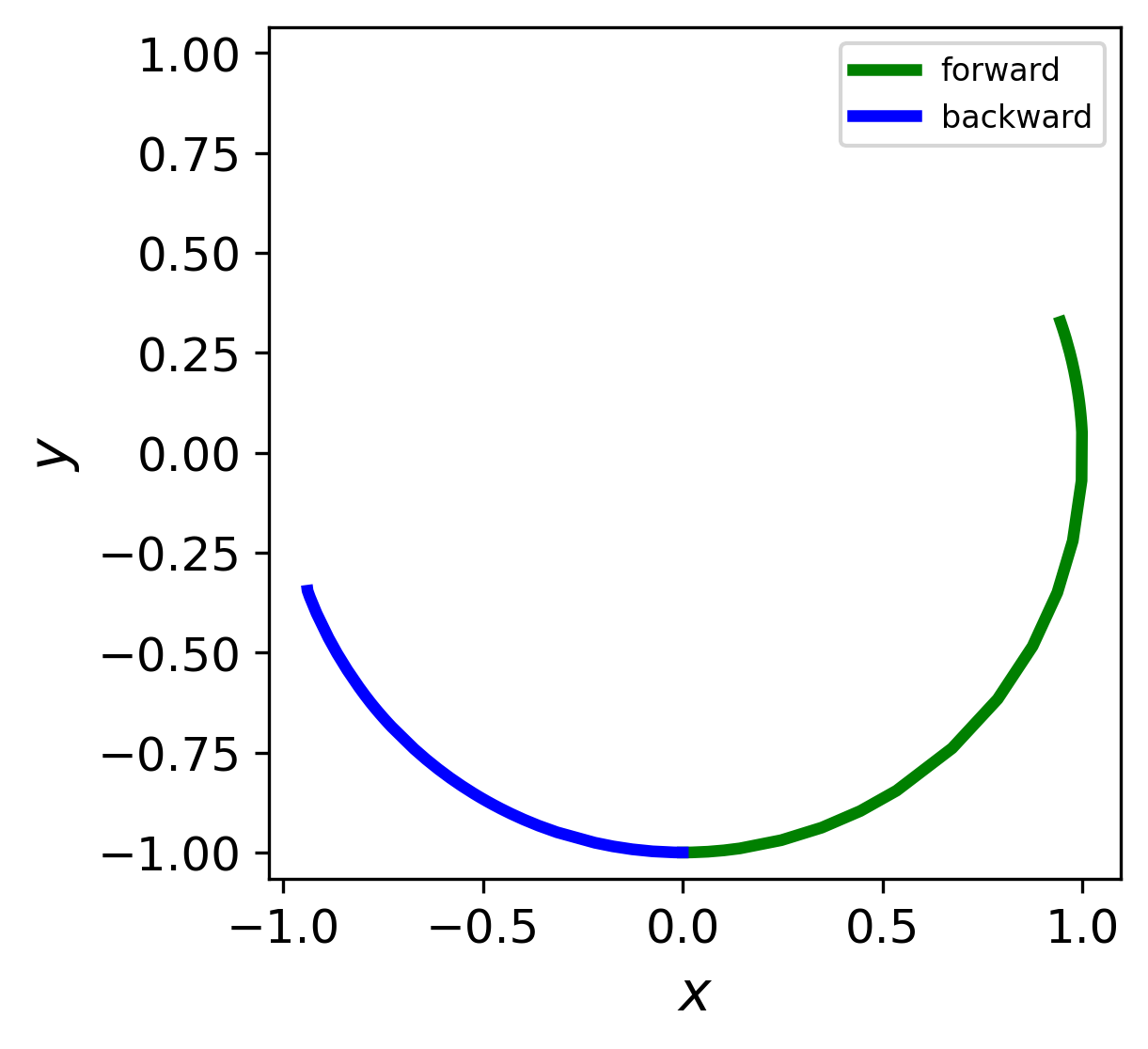}\hfill	\includegraphics[scale=0.4]{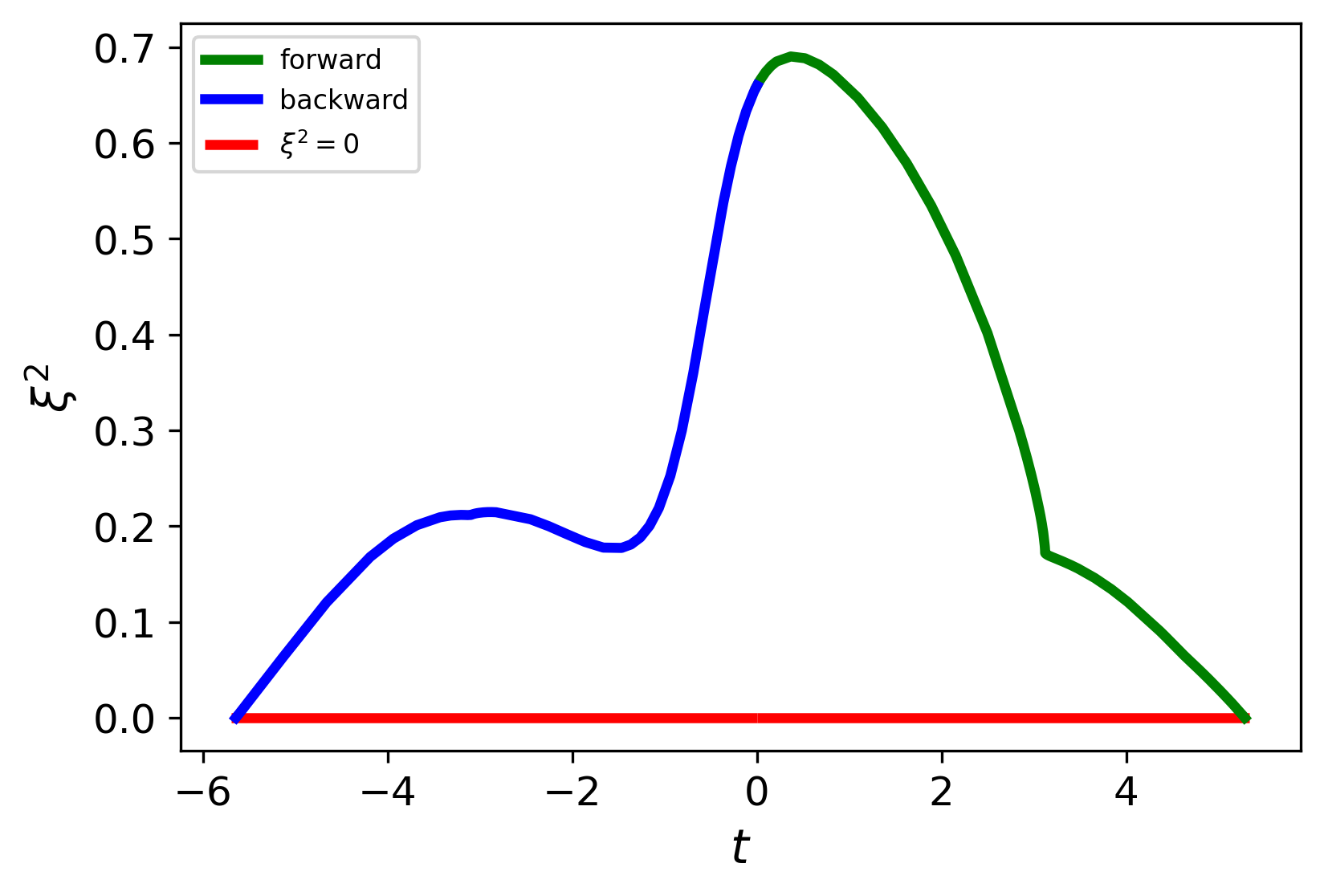}\hfill	\includegraphics[scale=0.4]{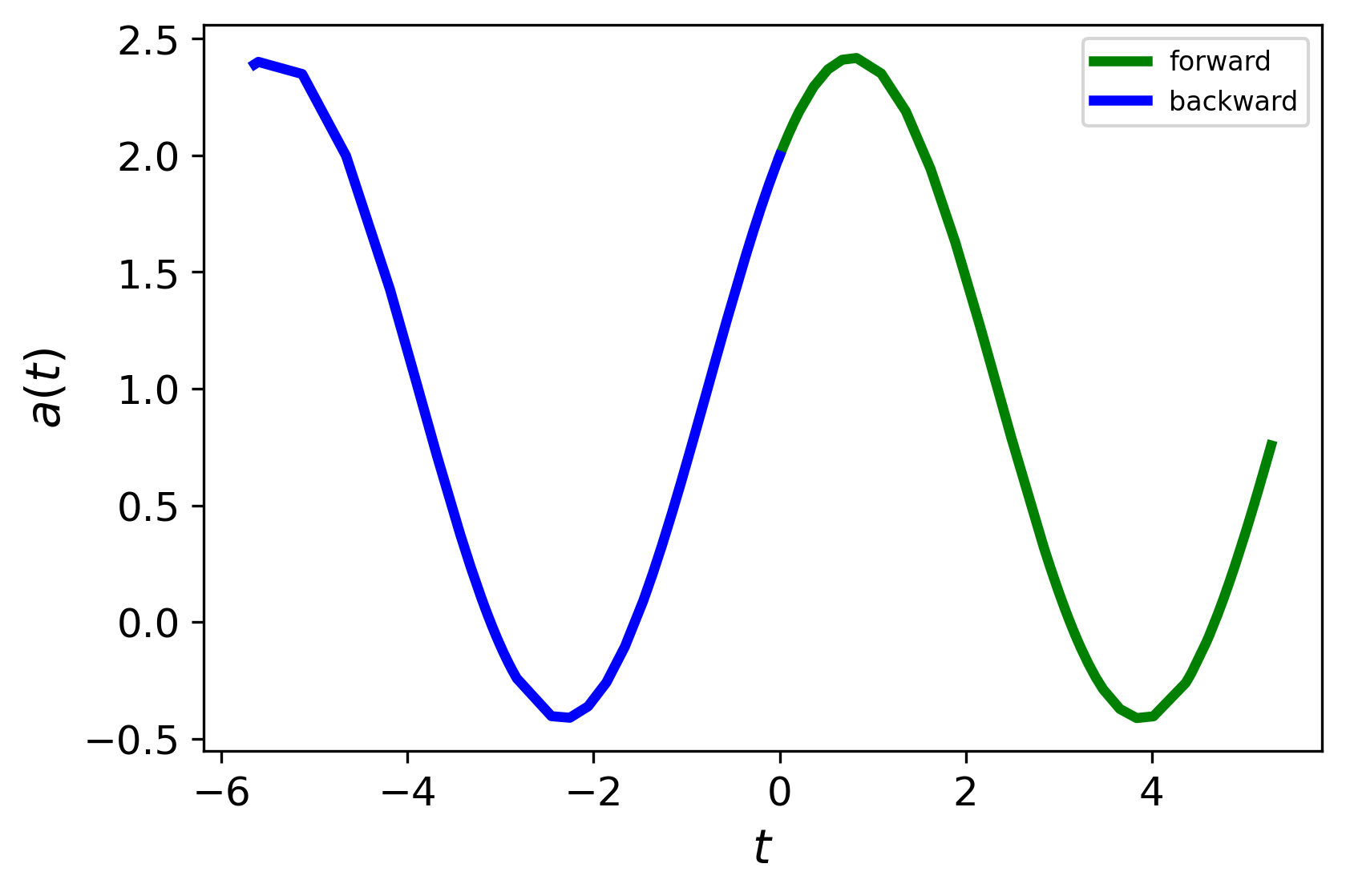}\hfill
	\includegraphics[scale=0.4]{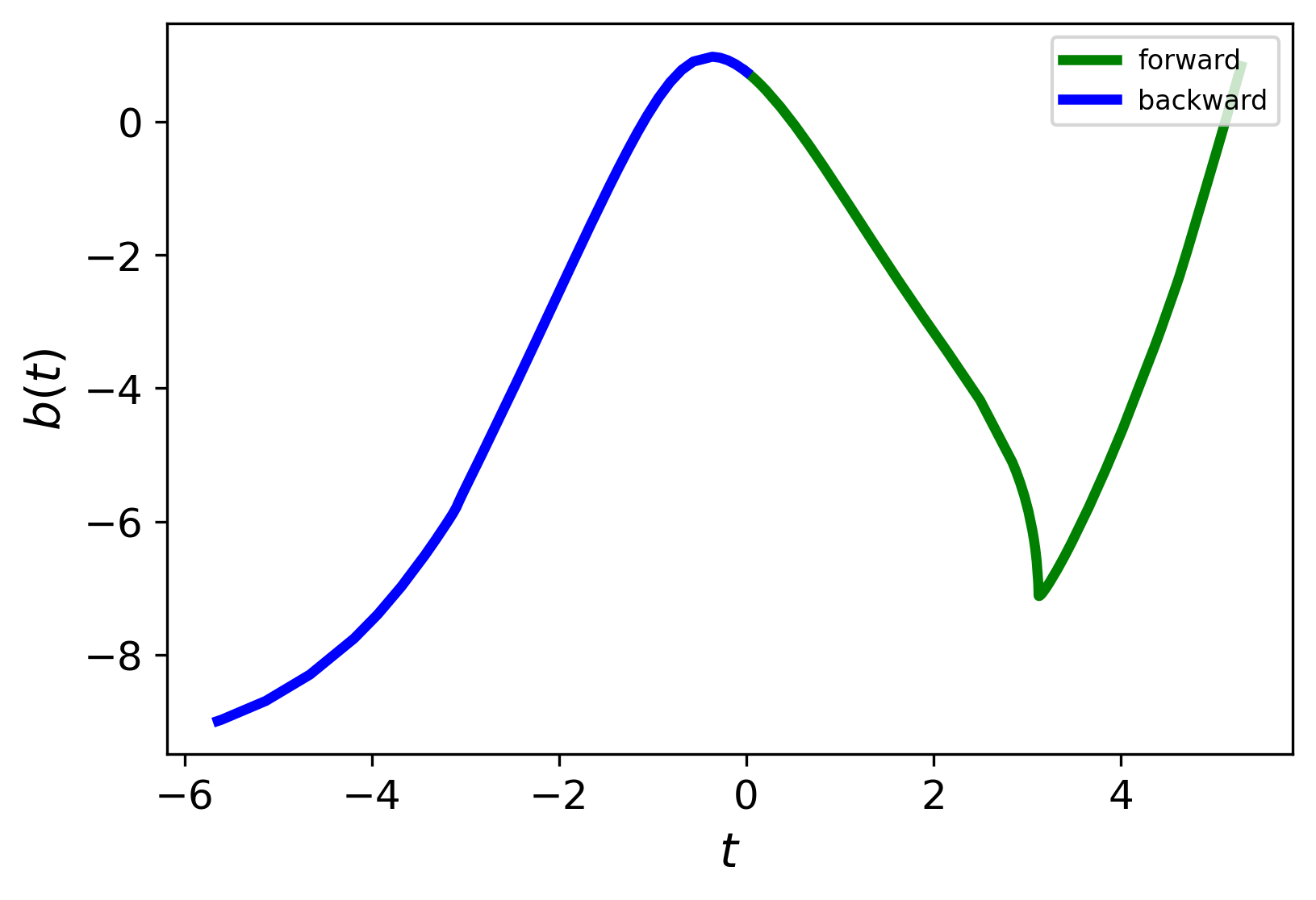}
	\caption{Arc 3 trajectory (top left), $\xi^2$ profile (top right), and optimized control functions, $a(t)$ (bottom left) and $b(t)$ (bottom right). Green and blue indicate the forward- and backward-in-time solutions, respectively.}	\label{fig:Arc3_Length}
\end{figure}
\begin{figure} 
	\centering 
	\includegraphics[scale=0.4]{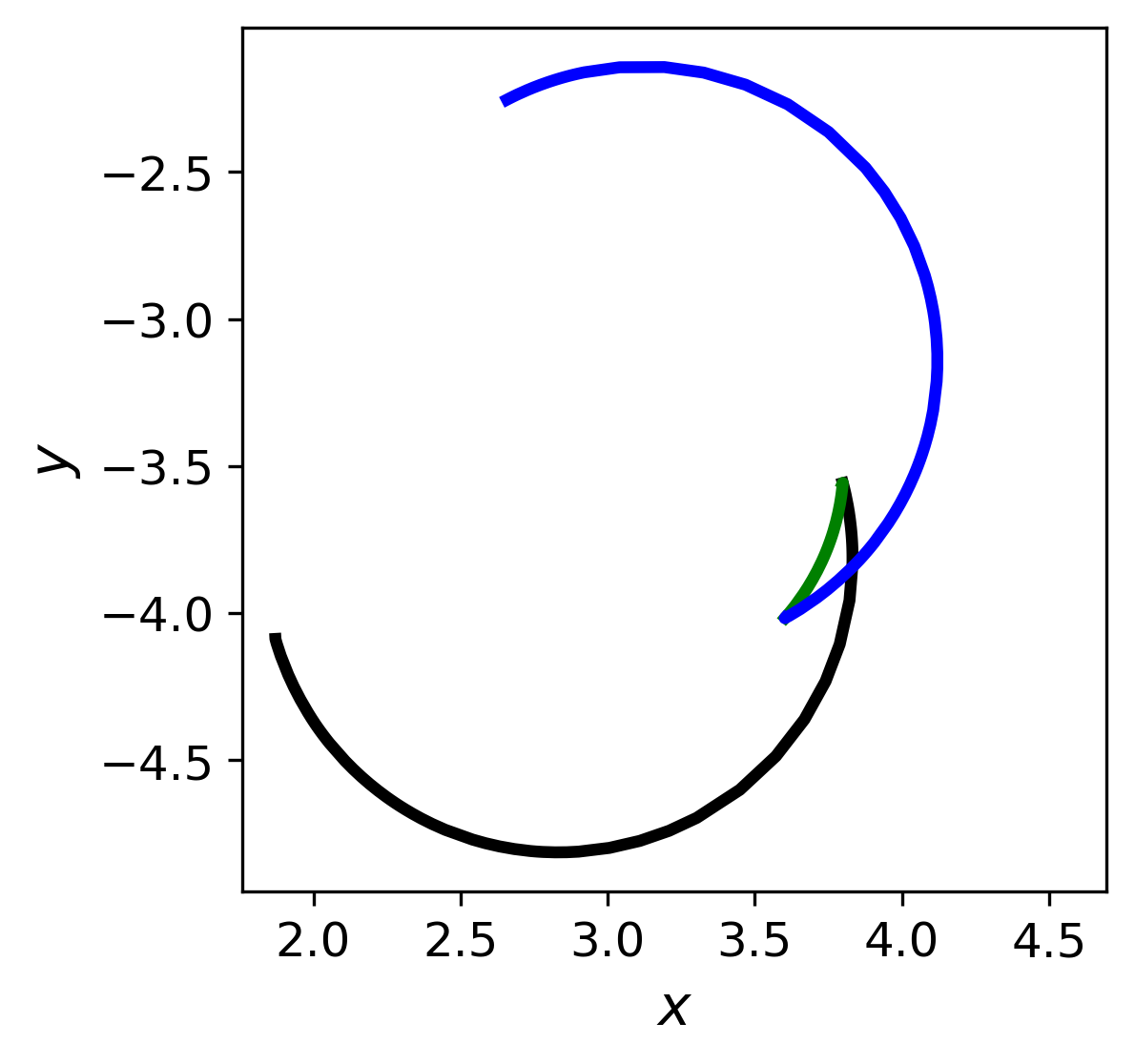}\hfill
	\includegraphics[scale=0.4]{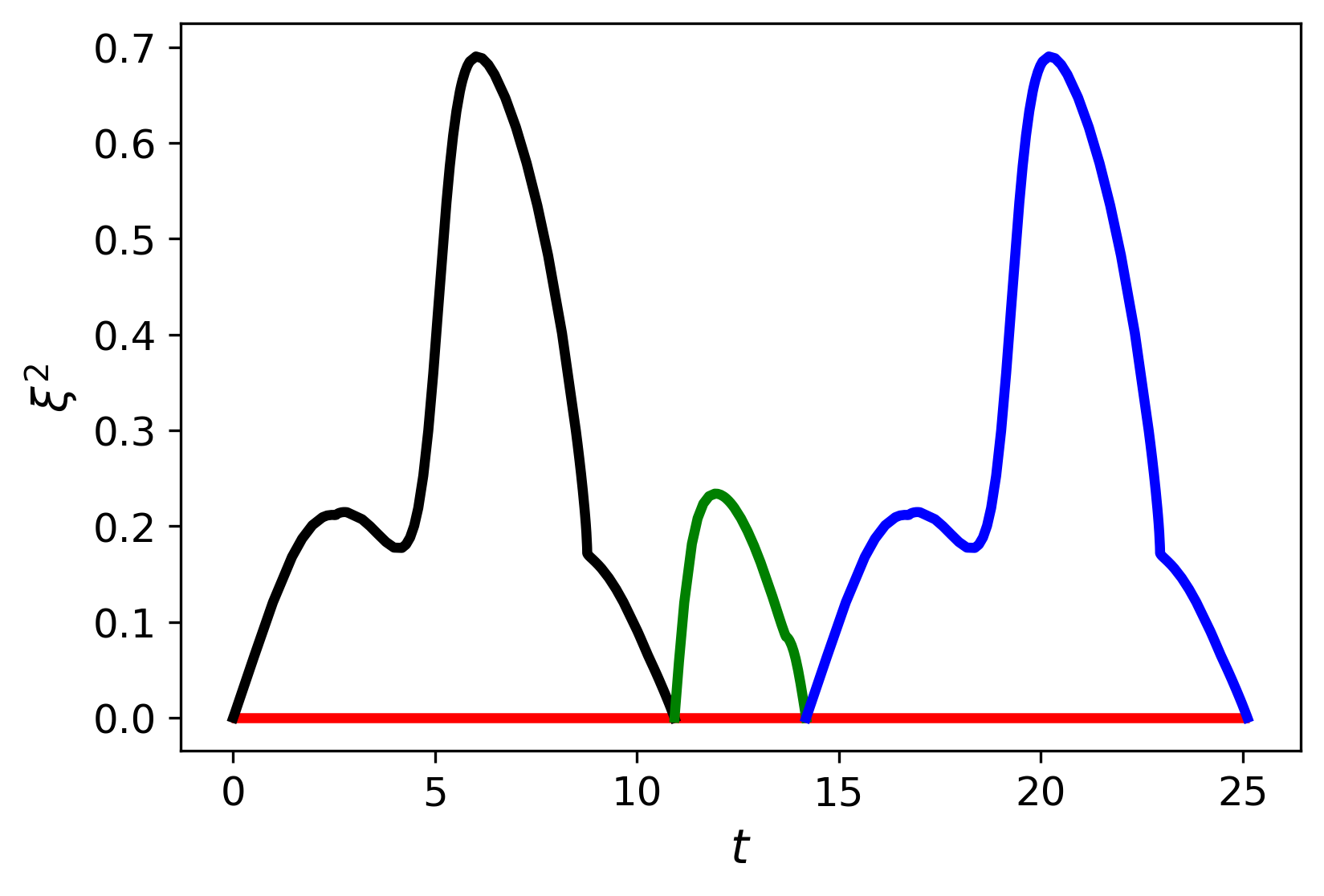}\\
	\includegraphics[scale=0.4]{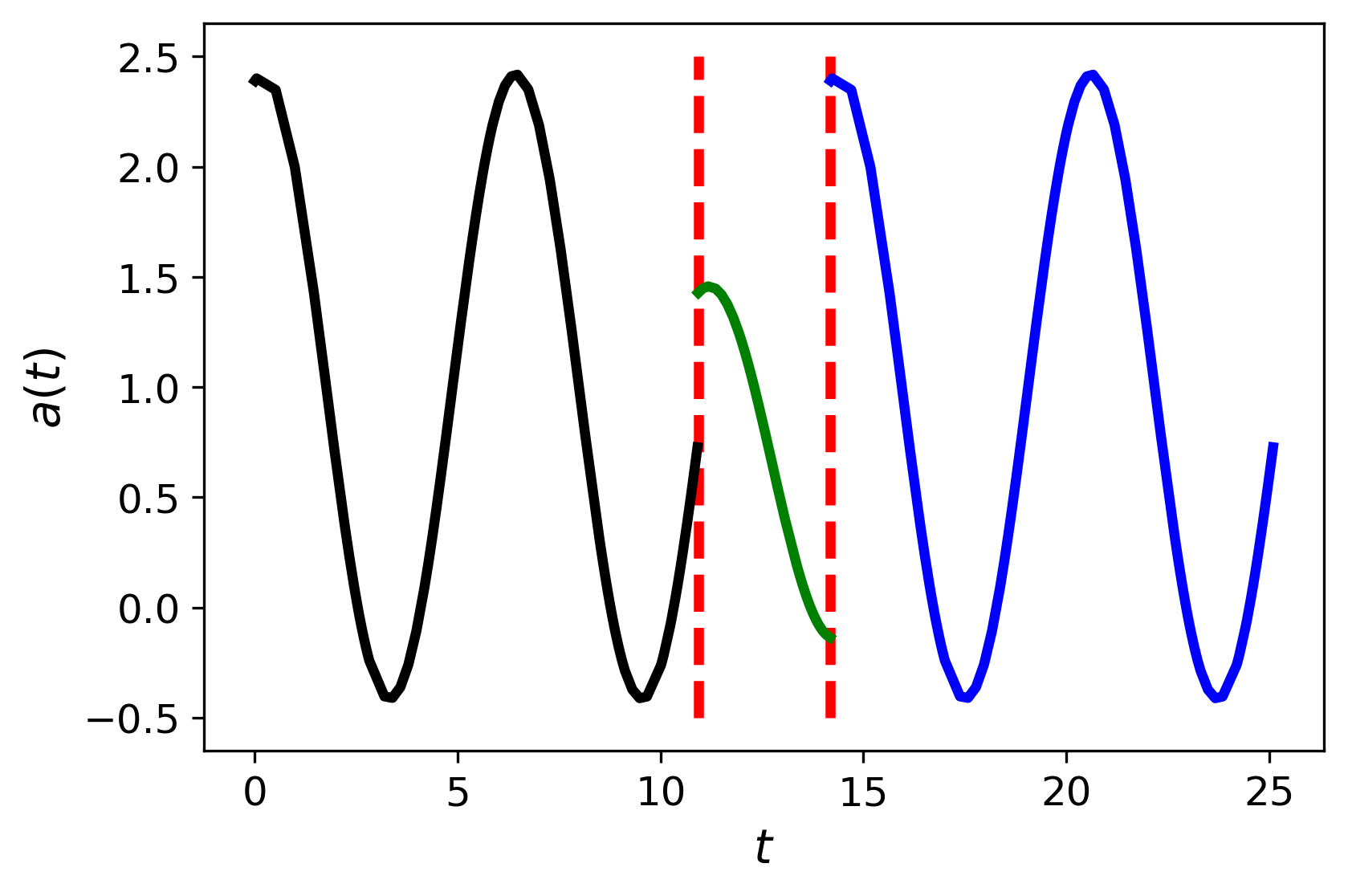}\hfill
	\includegraphics[scale=0.4]{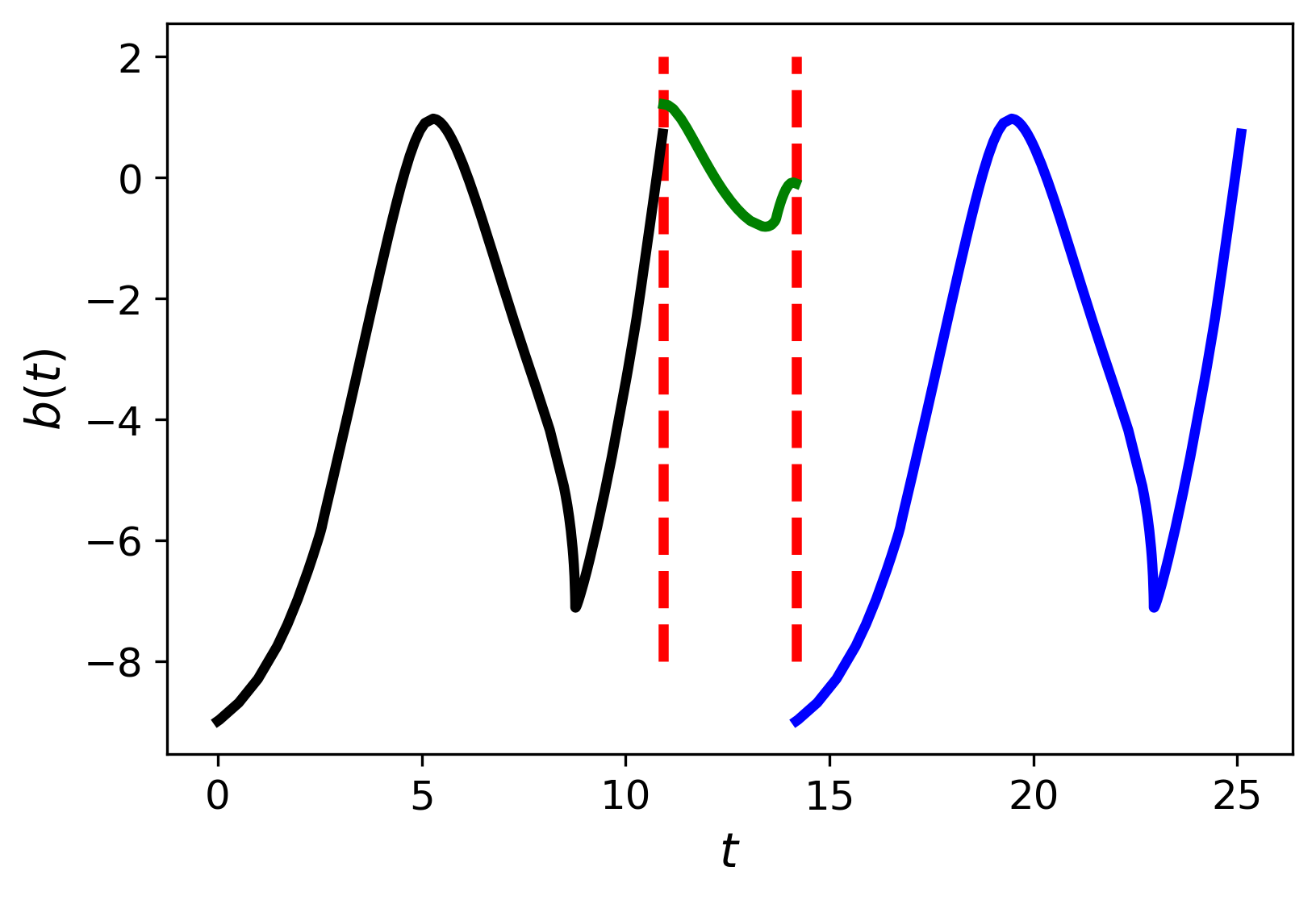}
	\caption{Trajectory of a single leaf (top left), corresponding $\xi^2$ profile (top right) and optimized control functions $a(t)$ (bottom left) and $b(t)$ (bottom right). The black, green, and blue segments distinguish the three arcs that make up the triplet. The vertical red lines indicate the time points where the controls switch between arc types.}
	\label{fig:Triplet_Length}
\end{figure}
\begin{figure} 
	\centering
	\includegraphics[height=3.2cm]{Figures/LengthBased with Energy June 4/Arc1}	\hfill\includegraphics[height=3.2cm]{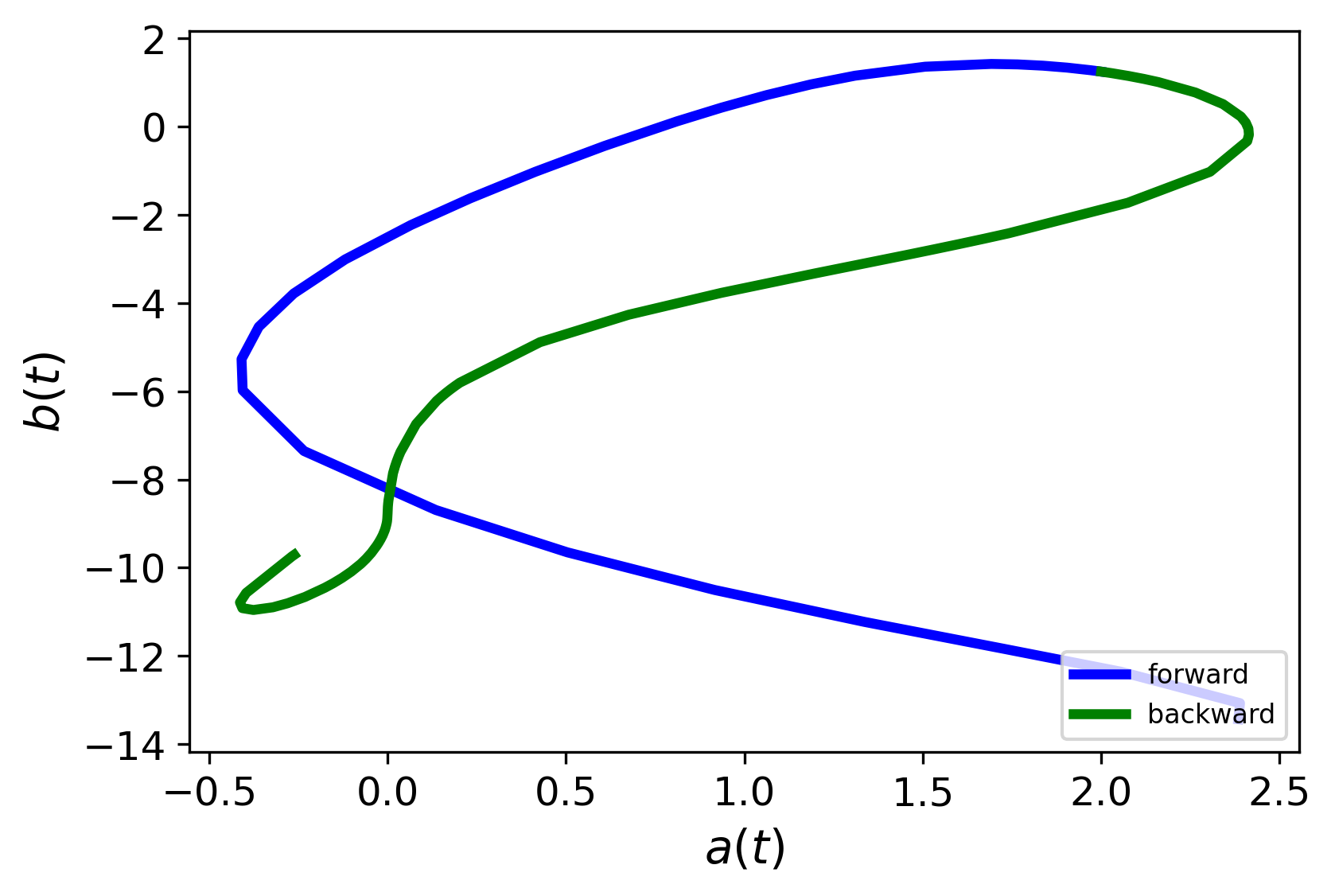}\hfill
	\includegraphics[height=3.2cm]{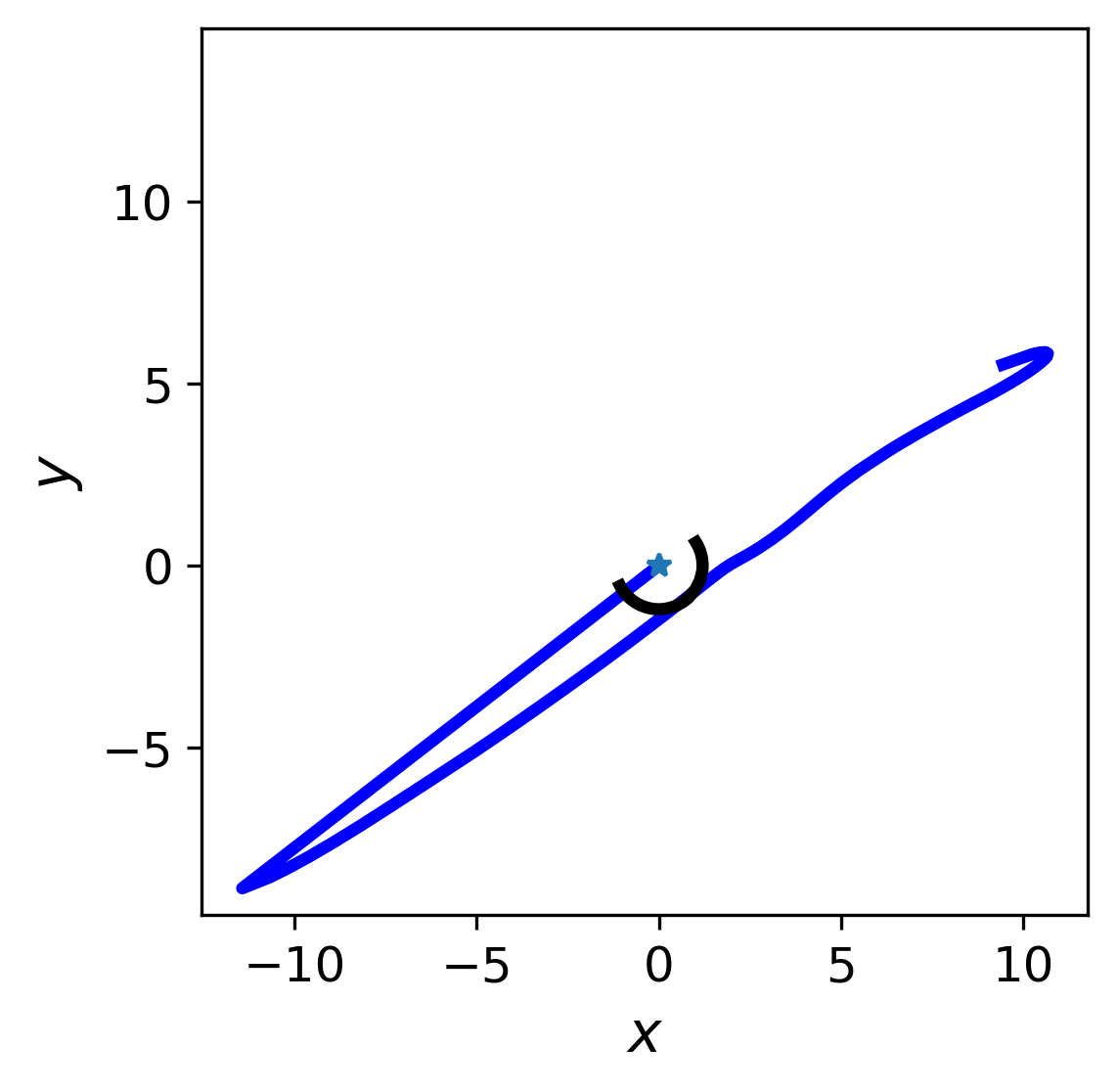}\\
	\includegraphics[height=3.2cm]{Figures/LengthBased with Energy June 4/Arc2}\hfill\includegraphics[height=3.2cm]{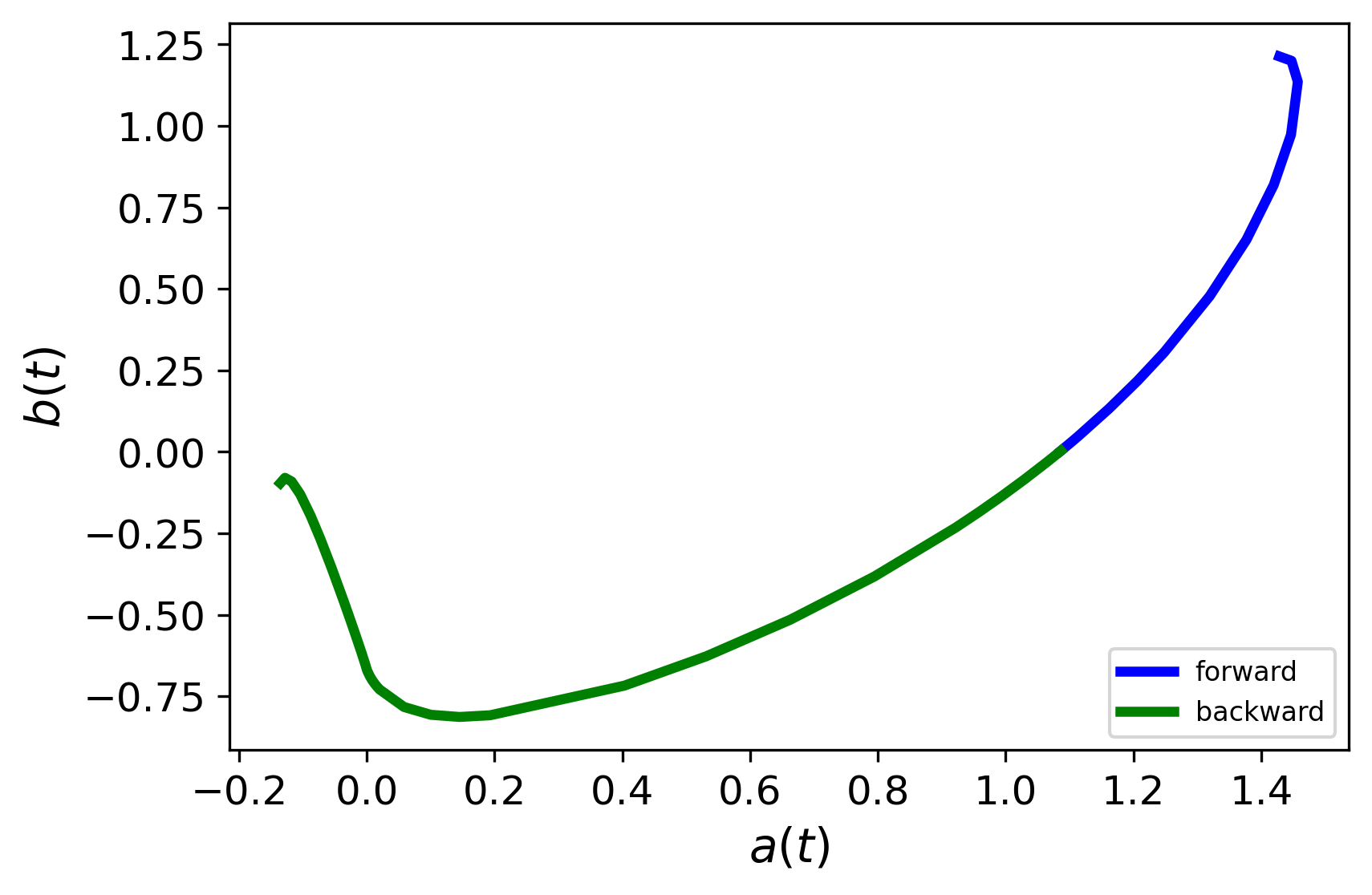}\hfill\includegraphics[height=3.2cm]{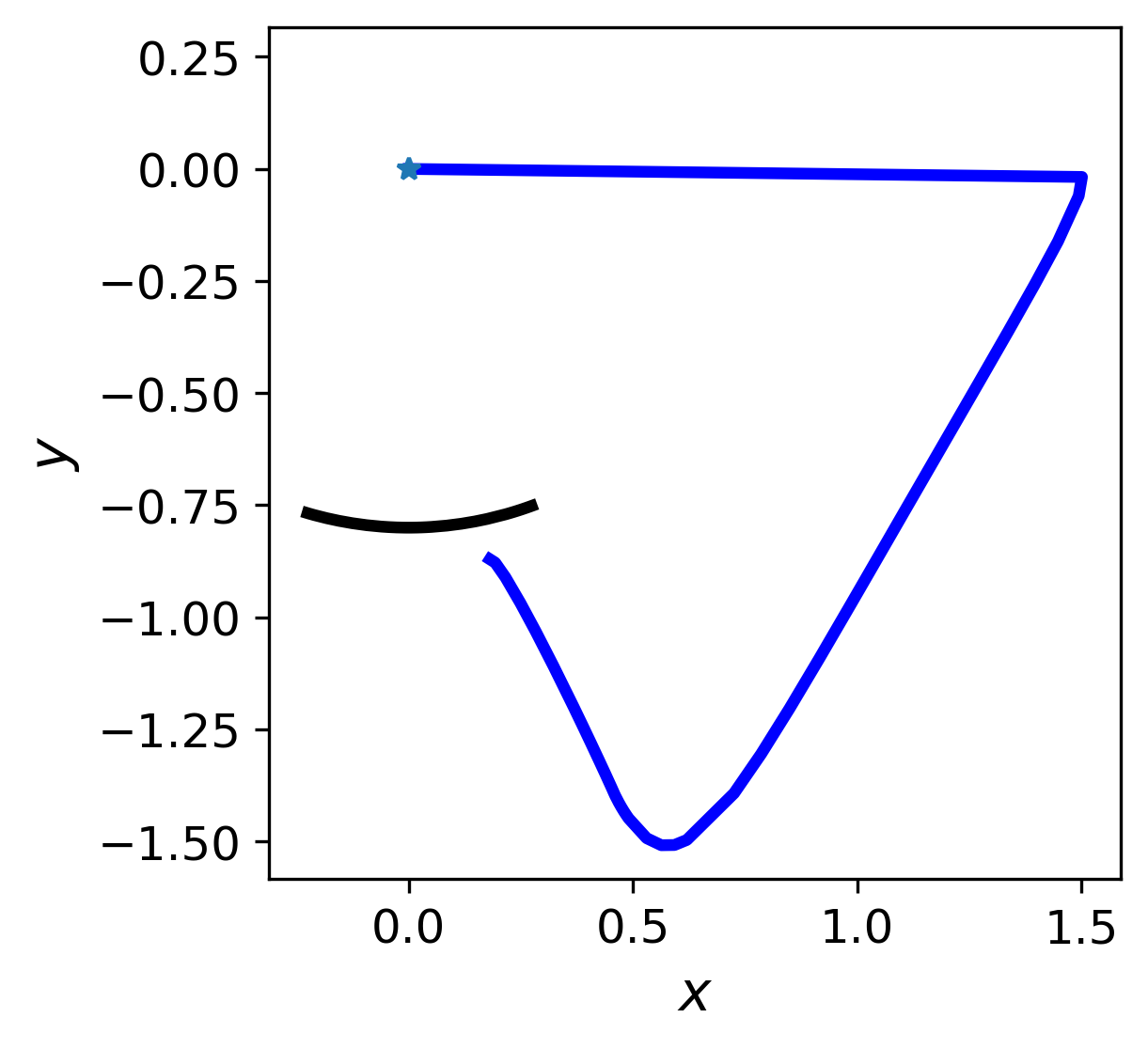}\\
	\includegraphics[height=3.2cm]{Figures/LengthBased with Energy June 4/Arc3}	\hfill\includegraphics[height=3.2cm]{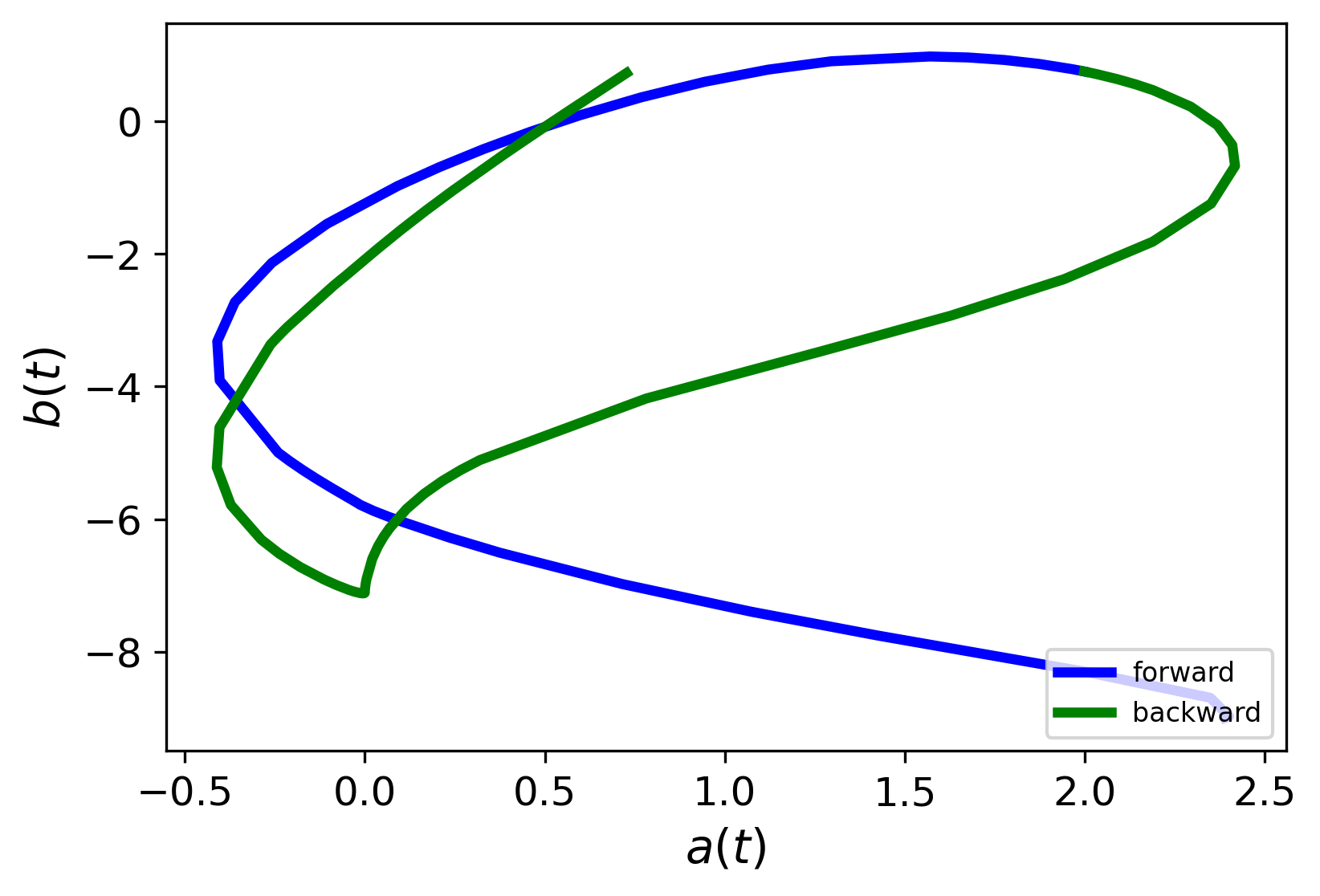}\hfill
	\includegraphics[height=3.2cm]{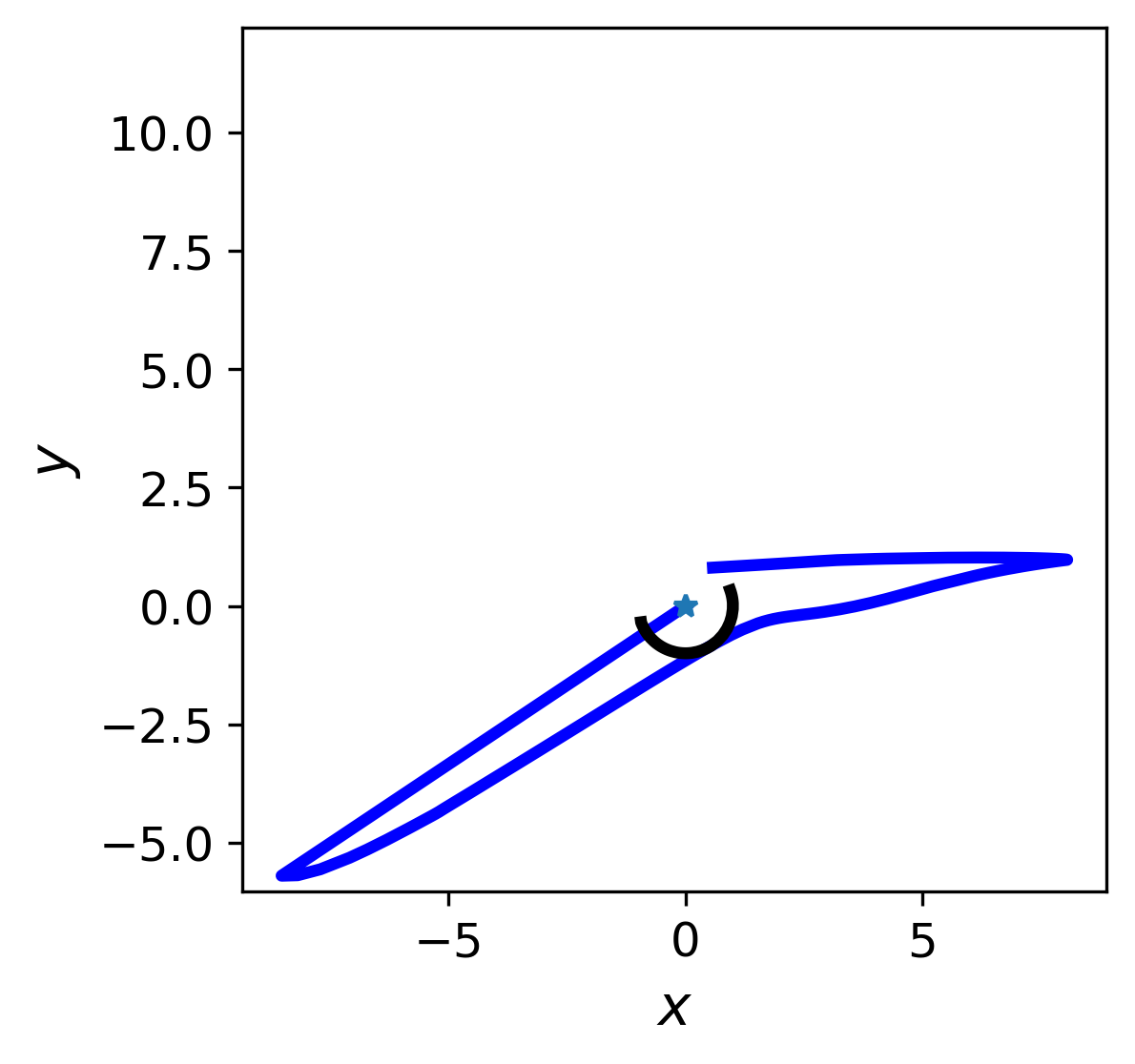}
	\caption{Length-optimized skate trajectories (left), control mass trajectories (middle), and overlay of skate and control mass trajectories in spatial frame (right) for arc 1 (top), arc 2 (middle), and arc 3 (bottom). Green and blue indicate the forward- and backward-in-time solutions, respectively. In the far right column, the star indicates the initial point of the control mass.}
	\label{fig:ArcControls_Length}
\end{figure}

The resulting length-optimized arcs, along with the corresponding $\xi^2$ profiles and optimized control functions, are shown in Figures \ref{fig:Arc1_Length}, \ref{fig:Arc2_Length}, and \ref{fig:Arc3_Length}. 
As shown in Figure \ref{fig:Triplet_Length}, the $\xi^2$ profiles for each leaf of the outer curve vanish at each connection point between the arcs, satisfying the physical requirement.
Figure \ref{fig:ArcControls_Length} shows the trajectory of the position (in body frame) of the control mass $m$, the blade trajectory (in spatial frame), and the overlay of the skate and control mass trajectories (in spatial frame) for each length-optimized arc.

\begin{figure} 
\begin{minipage}{0.49\textwidth}
    \begin{flushright}
        \includegraphics[scale=0.4]{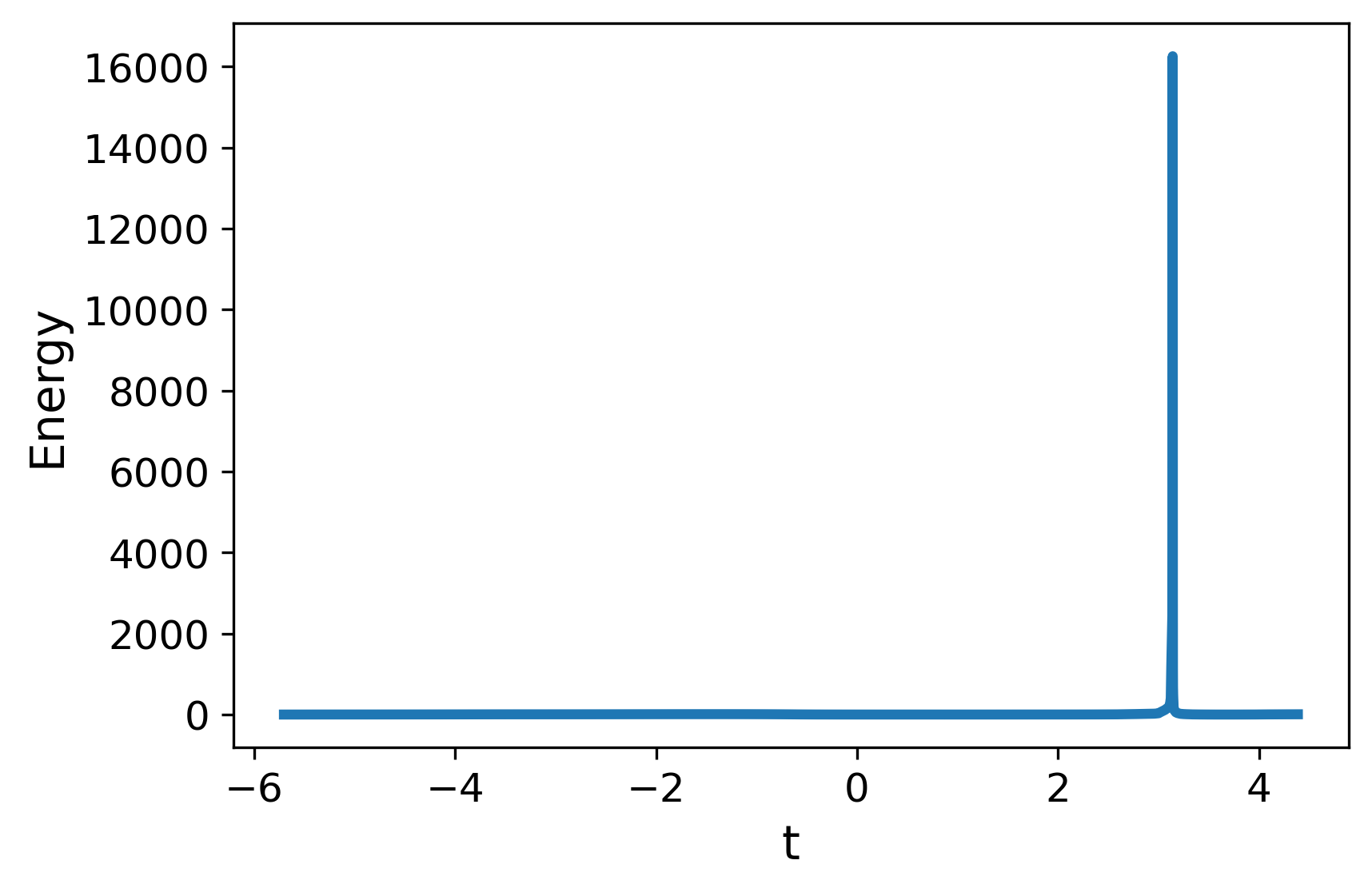}\\
         \includegraphics[scale=0.4]{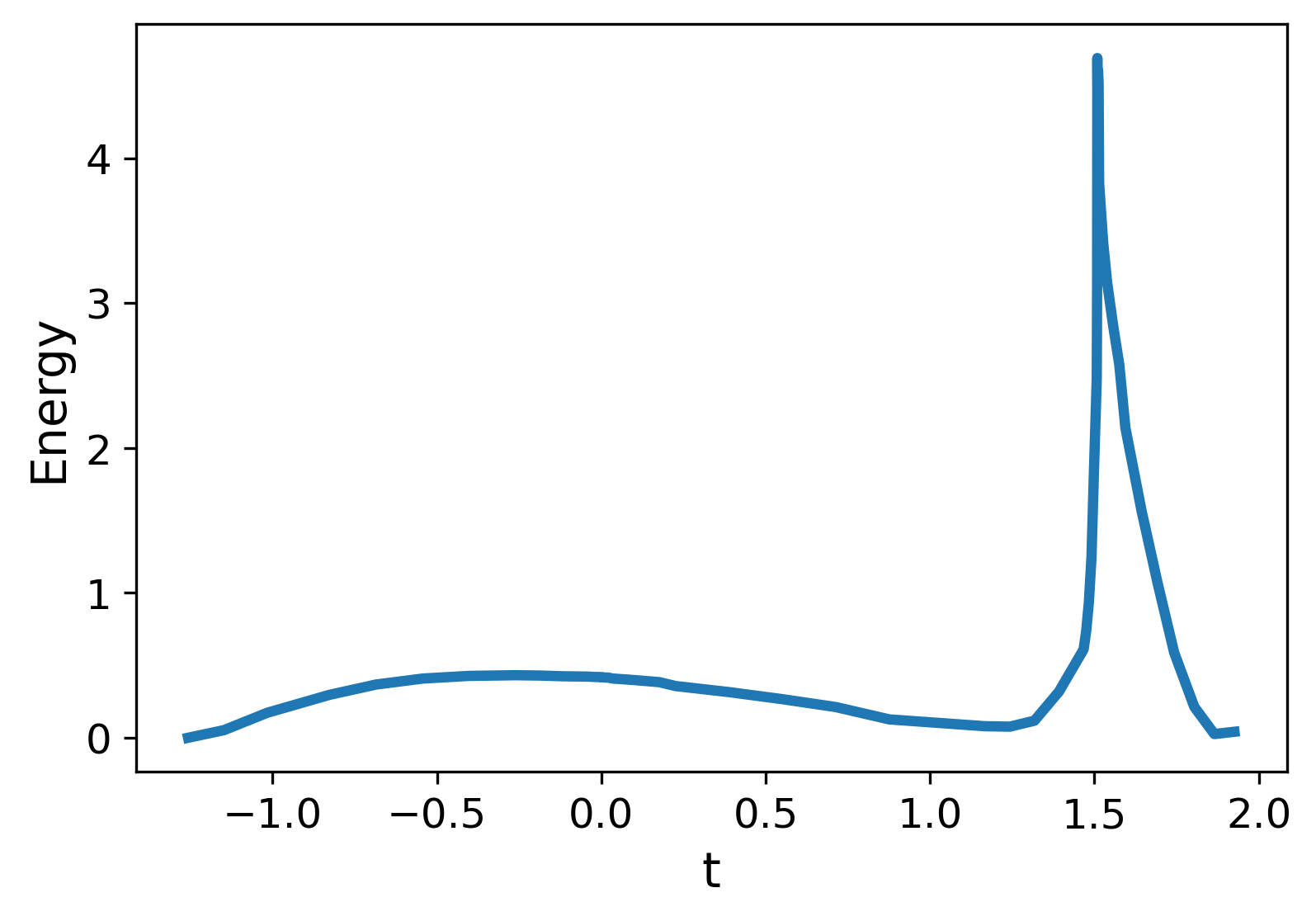}\\
         \includegraphics[scale=0.4]{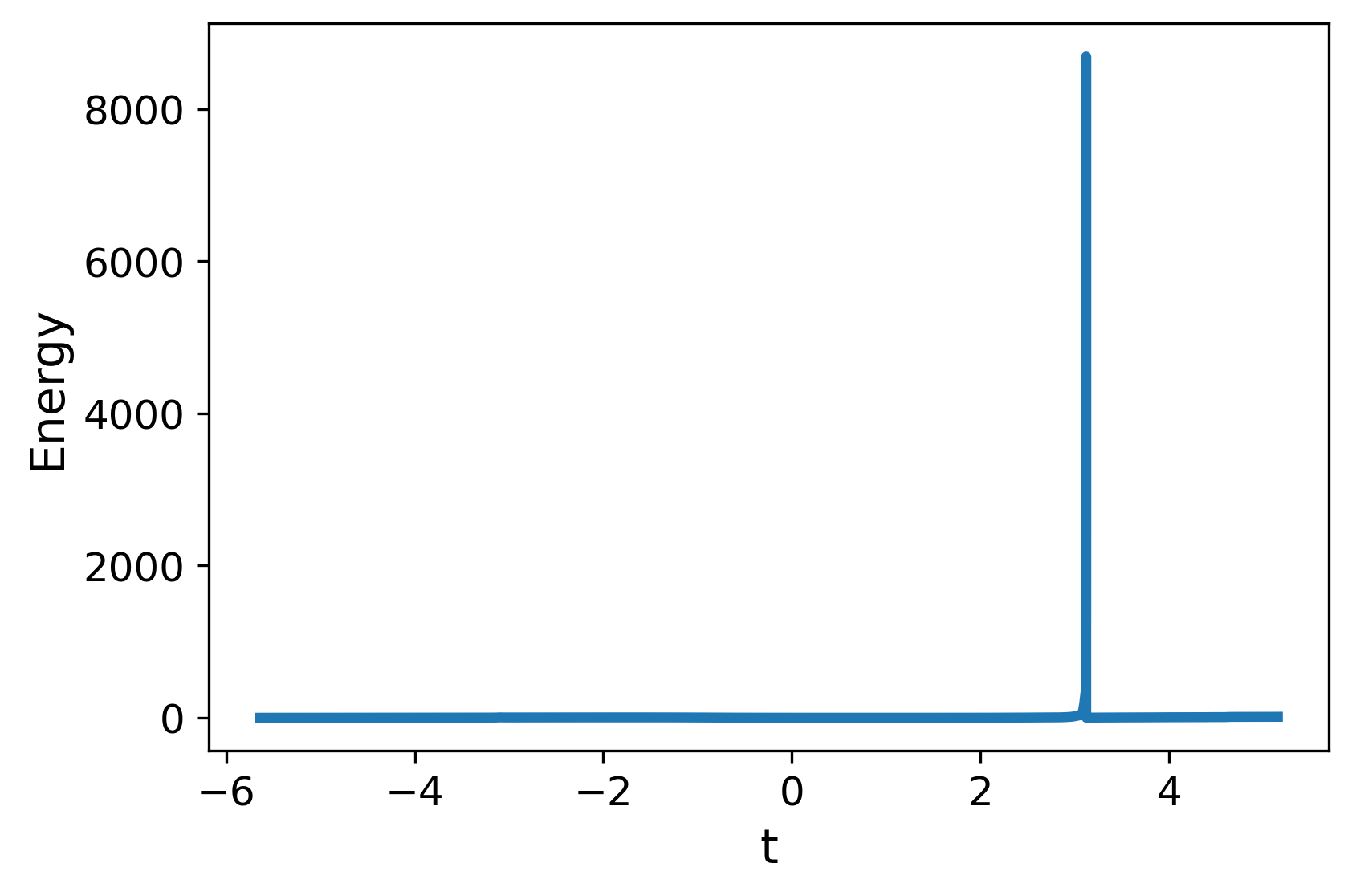}
    \end{flushright}
\end{minipage}
\begin{minipage}{0.49\textwidth}
    \begin{flushright}
        \includegraphics[scale=0.4]{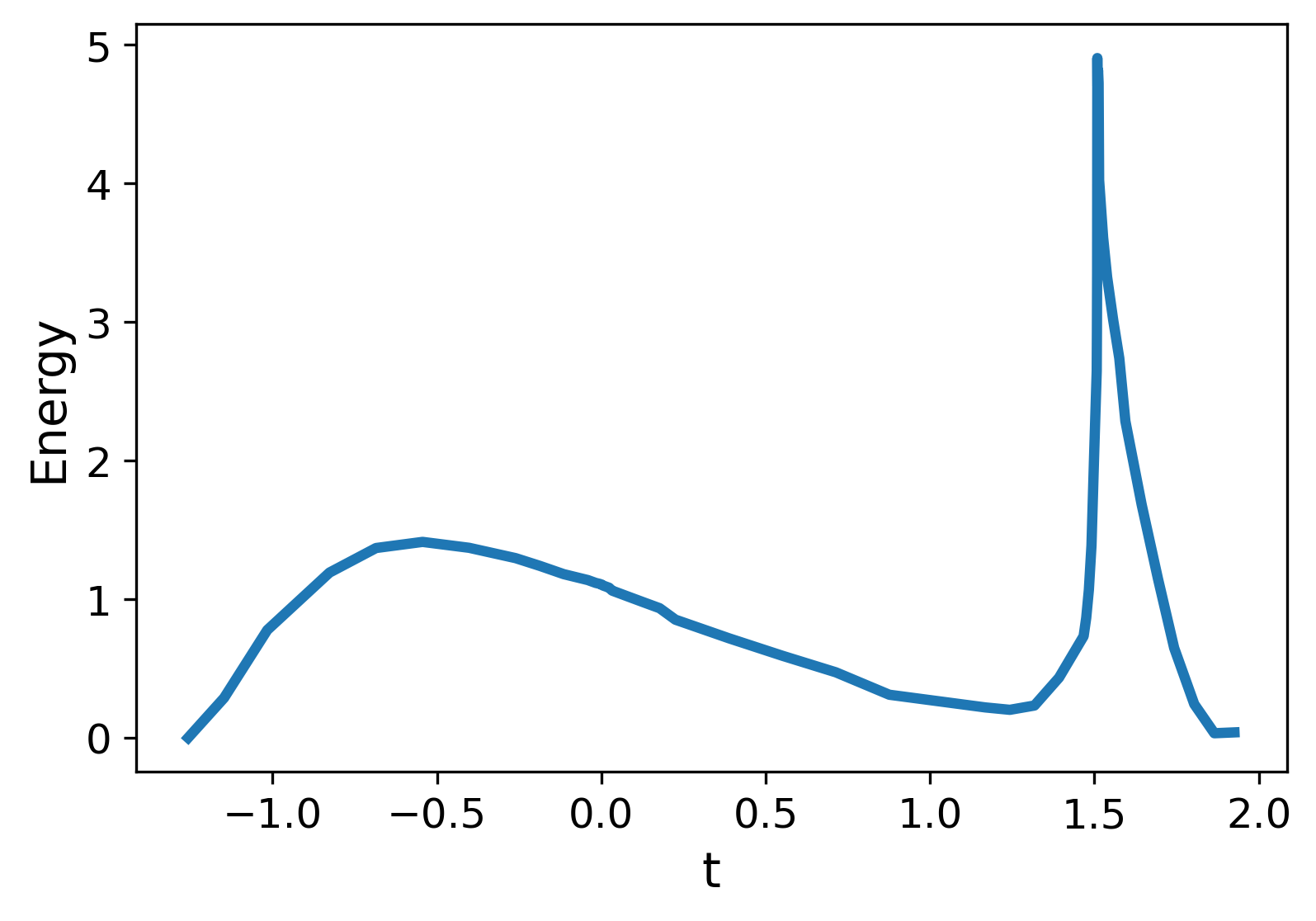}\\ \includegraphics[scale=0.4]{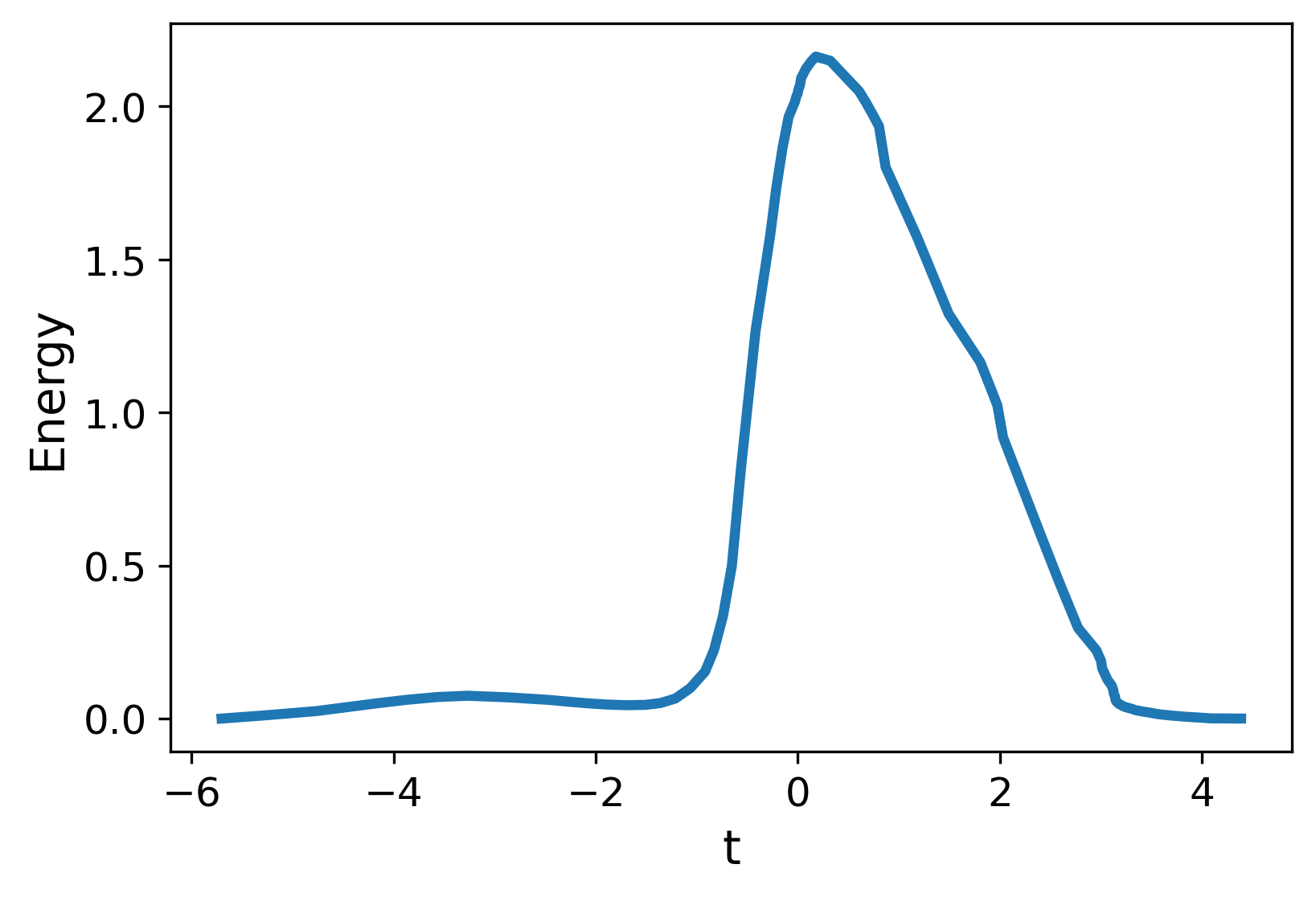}\\
        \includegraphics[scale=0.4]{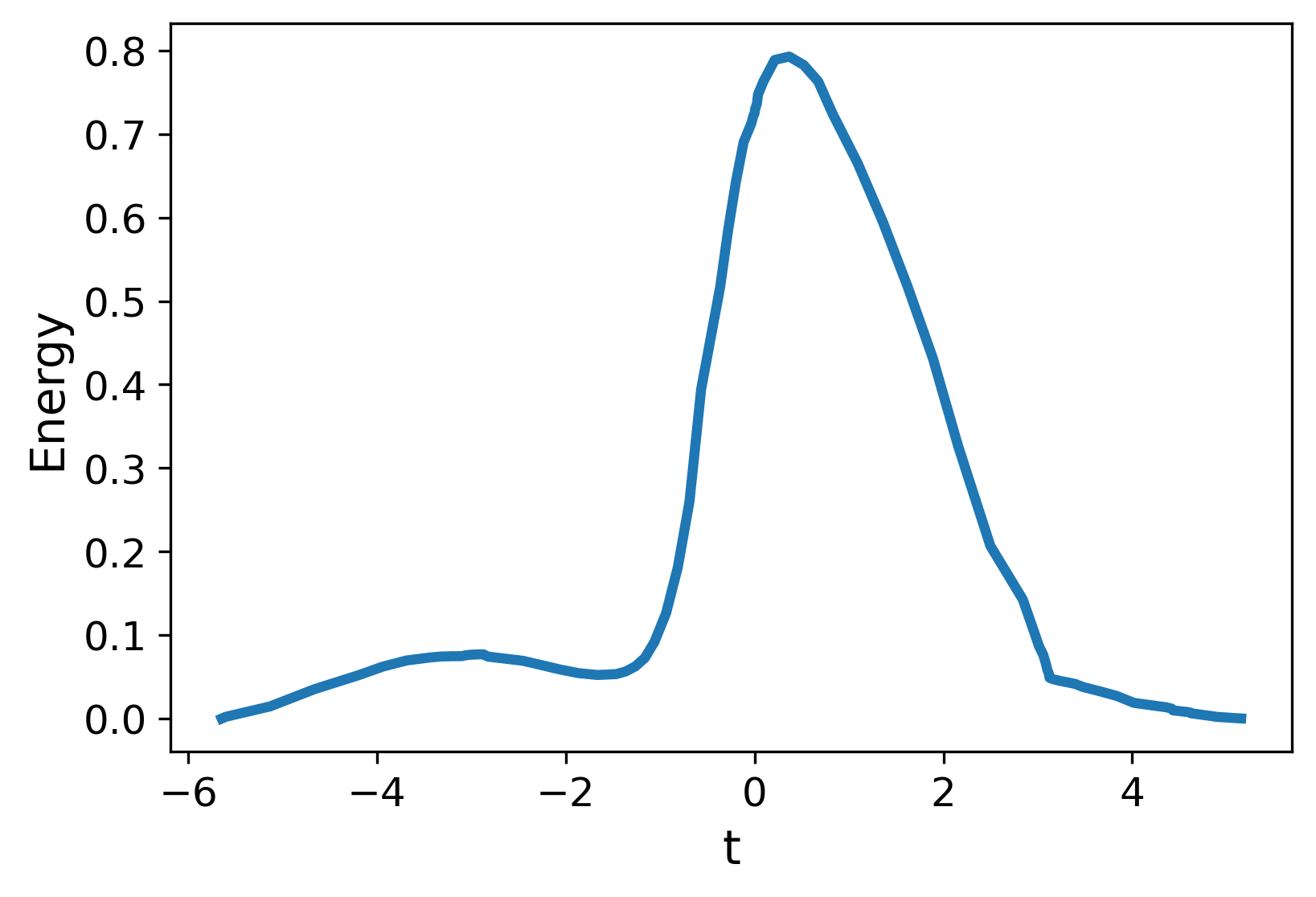}
    \end{flushright}
\end{minipage}
    \caption{Energy over time for control mass (left) and skate (right) for arc 1 (top), arc 2 (middle), and arc 3 (bottom) with length-based optimization.}
    \label{fig: Energy lengthbased}
\end{figure}

 The total energy of the skate and control mass are shown in Figure \ref{fig: Energy lengthbased}. In this case, the energy of the control mass spikes significantly for arc 1 and arc 3, which is not physically realistic for a skater. This behavior led to the use of the more restricted control functions used in the main portion of this paper. 

\begin{figure} 
	\centering
	\includegraphics[width=0.5\textwidth]{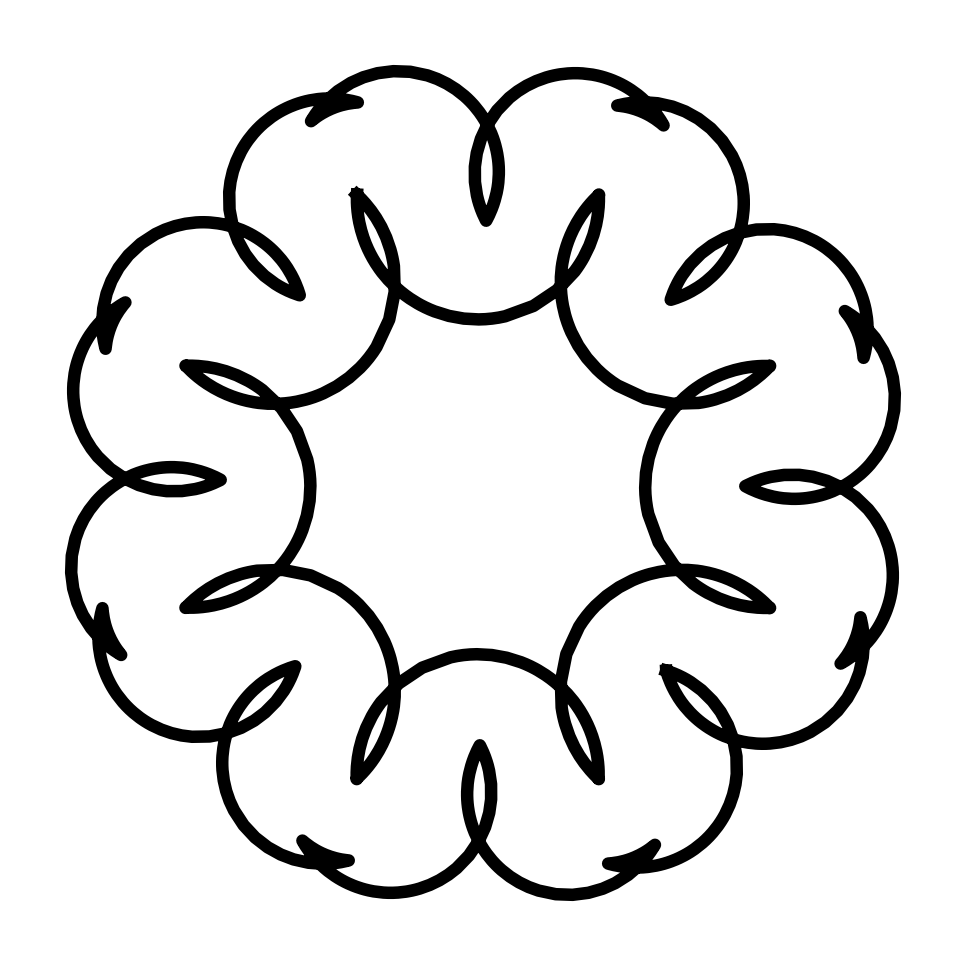}
	\caption{Full pattern reconstruction using length-based optimization.}
	\label{fig:SimulatedPattern_Length}
\end{figure}

The full pattern produced by transforming and repeating the length-optimized arcs is shown in Figure \ref{fig:SimulatedPattern_Length}. Again, notice the similarity with the pattern in Figure~\ref{fig:SimulatedPattern_Length_NewControl}, although the length of the arcs are slightly different from the pattern shown on Figure~\ref{fig:SimulatedPattern_Length_NewControl}.

\end{document}